\documentclass{article}
\usepackage{amsmath}[1996/11/01]
\usepackage{amssymb,amsthm,amsxtra}

\def\[#1\]{\begin{equation}#1\end{equation}}
\makeatletter
\def\beq{%
   \relax\ifmmode
      \@badmath
   \else
      \ifvmode
         \nointerlineskip
         \makebox[.6\linewidth]%
      \fi
      $$
   \fi
}
\def\eeq{%
   \relax\ifmmode
      \ifinner
         \@badmath
      \else
         $$
      \fi
   \else
      \@badmath
   \fi
   \ignorespaces
}

\def\enddisplaymath{\eeq\global\@ignoretrue}
\makeatother

\newtheorem{thm}{Theorem}
\newtheorem{cor}[thm]{Corollary}
\newtheorem{lem}[thm]{Lemma}
\newtheorem{prop}[thm]{Proposition}
\newtheorem{conj}{Conjecture}
\theoremstyle{remark}
\newtheorem*{rem}{Remark}
\newtheorem{rems}{Remark}[thm]

\theoremstyle{definition}
\newtheorem{defn}{Definition}

\numberwithin{equation}{section}
\numberwithin{thm}{section}

\newcommand{\Q}{\mathbb Q}
\newcommand{\R}{\mathbb R}
\newcommand{\Z}{\mathbb Z}
\newcommand{\F}{\mathbb F}
\newcommand{\T}{\mathbb T}

\newcommand{\binomQ}{\genfrac[]{0pt}{}}
\newcommand{\binomI}{\genfrac\{\}{0pt}{}}

\begin{document}

\title{$BC_n$-symmetric polynomials}
\author{Eric M. Rains\footnote{AT\&T Labs -- Research; Presently at:
    Department of Mathematics, University of California, Davis}}

\date{February 6, 2004}
\maketitle

\begin{abstract}
We consider two important families of $BC_n$-symmetric polynomials, namely
Okounkov's interpolation polynomials and Koornwinder's orthogonal
polynomials.  We give a family of difference equations satisfied by the
former, as well as generalizations of the branching rule and Pieri
identity, leading to a number of multivariate $q$-analogues of classical
hypergeometric transformations.  For the latter, we give new proofs of
Macdonald's conjectures, as well as new identities, including an inverse
binomial formula and several branching rule and connection coefficient
identities.  We also derive families of ordinary symmetric functions that
reduce to the interpolation and Koornwinder polynomials upon appropriate
specialization.  As an application, we consider a number of new integral
conjectures associated to classical symmetric spaces.

\end{abstract}

\tableofcontents

\section{Introduction}

In the theory of multivariate orthogonal polynomials, the family of
Koornwinder polynomials \cite{KoornwinderTH:1992} has special significance,
reducing to the Askey-Wilson polynomials \cite{AskeyR/WilsonJ:1985} in the
univariate case, as well as to the Macdonald polynomials for the classical
root systems in appropriate limits \cite{vanDiejen:1995}.  In the
present work, we prove a number of new results on these polynomials, as
well as giving new proofs of Macdonald's ``conjectures'' for these
polynomials (originally proved in \cite{vanDiejen:1996,SahiS:1999}).

The proofs in the literature of the Macdonald conjectures involve heavy use
of ``double affine Hecke algebra'' machinery, and are thus rather far
removed from the standard treatment of the Askey-Wilson case.  In contrast,
we will be taking a more classical approach; just as Askey-Wilson
polynomials are naturally represented and studied as hypergeometric
polynomials, we take as our starting point the multivariate analogue of
this representation, namely Okounkov's ``binomial formula''
\cite{OkounkovA:1998}.  To be precise, we define Koornwinder polynomials
via the binomial formula (or, equivalently, via the ``evaluation'' and
``duality'' portions of the Macdonald conjectures, together called
``evaluation symmetry'' below), and show that the resulting polynomials are
orthogonal with respect to both the Koornwinder \cite{KoornwinderTH:1992}
and multivariate $q$-Racah \cite{vanDiejenJF/StokmanJV:1998} inner
products.  In addition to the new results on Koornwinder polynomials that
this approach allows us to prove, we will see in a follow-up paper
\cite{bctheta} that the hypergeometric approach generalizes to the case of
elliptic hypergeometric functions, including an elliptic analogue
(biorthogonal abelian functions) of Koornwinder polynomials.  (The elliptic
case is also treated in \cite{xforms}, using a combination of difference
and (contour) integral operators.)

The binomial formula expands Koornwinder polynomials in terms of Okounkov's
$BC_n$-symmetric interpolation polynomials \cite{OkounkovA:1998}, which we
must therefore study first.  Our treatment of these polynomials is based on
a certain family of commuting $q$-difference operators; the joint
eigenfunctions $\bar{P}^{*(n)}_\lambda(;q,t,s)$ satisfy the ``extra
vanishing property''
\[
\bar{P}^{*(n)}_\lambda(\mu;q,t,s):=
\bar{P}^{*(n)}_\lambda(q^{\mu_i} t^{n-i} s;q,t,s)
=
0
\quad\lambda\not\subset\mu,
\]
and thus agree with Okounkov's polynomials, up to rescaling.  The resulting
difference equations lead to a new proof of the branching rule for
interpolation polynomials; more generally, we obtain a ``bulk'' version of
the branching rule (in which some of the variables are set to $v$, $tv$,
$t^2 v$,\dots).  Analogously, we also obtain a bulk Pieri identity,
containing both the known $e$-type Pieri identity and a new $g$-type Pieri
identity as special cases.  One important consequence of these bulk
identities is a formula for connection coefficients between interpolation
polynomials for different values of $s$.

One of the interpolation polynomials in the binomial formula shows up
as a ``binomial coefficient'':
\[
\binomQ{\lambda}{\mu}_{q,t,s}
:=
\frac{\bar{P}^{*(n)}_\mu(\lambda;q,t,s t^{1-n})}
{\bar{P}^{*(n)}_\mu(\mu;q,t,s t^{1-n})}.
\]
If we evaluate the bulk Pieri identity at a partition, the resulting sum
can be expressed in terms of binomial coefficients, and turns out to be a
multivariate $q$-analogue of the Saalsch\"utz summation formula for a
${}_3F_2$ hypergeometric series.  The symmetries of this generalized
$q$-Saalsch\"utz formula lead to a formula for the inverse of the matrix
of binomial coefficients, as well as the duality symmetry
\[
\binomQ{\lambda}{\mu}_{q,t,s}
=
\binomQ{\lambda'}{\mu'}_{t,q,1/\sqrt{qt}s}.\footnote{A number of symmetries
of Koornwinder polynomials have been called ``duality'' in the literature;
in the present work, we reserve the term for the symmetry under conjugation
of partitions, e.g., the action of Macdonald's involution on symmetric functions.}
\]
Furthermore, the formula is general enough that a number of arguments
\cite[Chapter 2]{GasperG/RahmanM:1990} lift from the univariate case; in
particular, we obtain a multivariate analogue of Watson's transformation
between a very-well-poised, terminating ${}_8\phi_7$ and a balanced,
terminating ${}_4\phi_3$.

As we mentioned, we take evaluation symmetry as the defining property of
Koornwinder polynomials; this then immediately implies a version of
Okounkov's binomial formula.  Applying the difference equations to the
interpolation polynomials in this formula gives a number of new identities,
notably special cases of branching rules and connection coefficients.
Moreover, we find that the corresponding difference operators act nicely on
Koornwinder polynomials, and have nice adjoints with respect to the
Koornwinder inner product; these facts combine to show that our Koornwinder
polynomials are orthogonal with respect to the correct inner product, and
thus the Koornwinder polynomials as usually defined satisfy evaluation
symmetry.  A different approach works in the $q$-Racah case; here, our
generalized hypergeometric transformations can be used to directly lift the
univariate proof \cite[Section 7.2]{GasperG/RahmanM:1990}.  Other
consequences of our theory include an inverse to the binomial formula
(expanding an interpolation polynomial in Koornwinder polynomials), a
connection coefficient formula for Koornwinder polynomials, and a
generalization of the Nasrallah-Rahman integral representation of a
very-well-poised ${}_8\phi_7$ \cite[Section 6.3]{GasperG/RahmanM:1990}.

Most of the formulas alluded to above are expressed in terms of ordinary
Macdonald polynomials; more precisely, in terms of principal
specializations of skew Macdonald polynomials.  This suggests that there
may be further connections between interpolation polynomials and ordinary
symmetric functions.  Indeed, it turns out that there are two
four-parameter families of symmetric functions that reduce to the
three-parameter family of interpolation polynomials upon appropriate
specialization.  Similarly, via the binomial formula, one obtains two
seven-parameter families that reduce to Koornwinder polynomials.  In
addition to some curious symmetries of these lifted Koornwinder
polynomials, having no counterparts for integer $n$, we also obtain a
refinement of the fact that the leading terms of interpolation and
Koornwinder polynomials are Macdonald polynomials, namely that the
$BC_n$-symmetric polynomials can be expressed in the triangular form
\[
\sum_{\mu\subset\lambda} c_{\lambda\mu} P_\mu(x_1^{\pm 1},x_2^{\pm 1},\dots
x_n^{\pm 1};q,t).
\]
We also consider the limit $n\to\infty$ of the $n$-variable Koornwinder
polynomials.

As an application of our results, we consider a pair of ``vanishing
conjectures''.  These are analogues for Macdonald polynomials of the fact
that the integral of a Schur function over the orthogonal (or symplectic)
group vanishes unless the corresponding irreducible representation of
$U(n)$ has a vector fixed by the orthogonal group.  More precisely, we
conjecture that the integral of $P_\lambda(;q,t)$ over a suitable
Koornwinder weight vanishes unless all parts of $\lambda$ are even; dually,
the integral over a different Koornwinder weight vanishes unless all parts
of $\lambda'$ are even.  In addition to the Schur case $q=t$, we can prove
a several other special cases, most notably that the first conjecture holds
if either $\lambda_1\le 4$ or $\ell(\lambda)\le 1$.  The latter case
involves a quadratic transformation of basic hypergeometric series, and
thus our conjectures can be interpreted as multivariate analogues of
quadratic transformations.  We also give a number of analogous conjectures
related to other classical symmetric spaces, obtained by various {\em ad
  hoc} methods.  (In fact, most of the conjectures of section
\ref{sec:vanish} have since been proved; see \cite{quadtrans}.  More
precisely, that paper proves Conjectures \ref{conj:vanish_U/O_T},
\ref{conj:vanish_U_grass2}, and \ref{conj:vanish_Sp_grass} using (single)
affine Hecke algebra techniques; Conjectures \ref{conj:vanish_U/Sp_T} and
\ref{conj:vanish_O_grass} are not amenable to that approach, but are
equivalent to Conjectures \ref{conj:vanish_U/O_T} and
\ref{conj:vanish_Sp_grass} as shown below.)

The paper is organized as follows.  Section \ref{sec:sfnotation} introduces
the notations for partitions and (ordinary) Macdonald polynomials we will
need below (for basic hypergeometric series, see
\cite{GasperG/RahmanM:1990}), as well as proving a few transformation
results for principal specializations of skew Macdonald polynomials.  We
then introduce the interpolation polynomials in Section \ref{sec:interp},
proving their main properties; then Section \ref{sec:hyperg} treats the
corresponding hypergeometric transformations.  Section \ref{sec:koorn}
introduces Koornwinder polynomials and proves the main results discussed
above (in particular, the Macdonald conjectures).  Sections
\ref{sec:sfinterp} and
\ref{sec:sfkoorn} then consider the corresponding families of symmetric
functions, with special attention to the case $n\to\infty$ of lifted
Koornwinder polynomials.  Finally, Section \ref{sec:vanish} states the various
vanishing conjectures, and proves the known special cases.

\subsection*{Acknowledgements}
Throughout this work, Peter Forrester has been a very helpful sounding
board and guide to the literature; many thanks are thus due.  Appreciation
is also due to the organizers of the conference on Applications of Macdonald
Polynomials at the Newton Institute in April 2001, which is where the idea
of looking for vanishing integrals like those of Section \ref{sec:vanish} 
first arose.

\section{Notation}\label{sec:sfnotation}

For partitions and Macdonald polynomials, we use the notations of
\cite{MacdonaldIG:1995}, Chapters I and VI, supplemented as follows.
First, we will denote the number of parts of $\lambda$ by $\ell(\lambda)$
(being somewhat clearer than $\lambda'_1$).  Next, we extend the dominance
partial order $\le$ to partitions not of the same size in the obvious way;
note that conjugation is no longer order-reversing ($311\ge 22$, even
though both partitions are self-conjugate).  We also define relations
$\prec$ and $\succ$ such that $\kappa\prec\lambda$ (equivalently
$\lambda\succ\kappa$) for two partitions iff $\lambda/\kappa$ is a vertical
strip; that is, $\kappa_i\le \lambda_i\le \kappa_i+1$ for all $i$.  We also
introduce some transformations of partitions.  For a partition $\lambda$,
we define partitions $2\lambda$ and $\lambda^2$ by
$(2\lambda)_i=2\lambda_i$, $(\lambda^2)_i=\lambda_{\lceil i/2\rceil}$.
Finally, if $\ell(\lambda)\le n$, $m^n+\lambda$ denotes the partition such
that $(m^n+\lambda)_i=m+\lambda_i$; if also $\lambda_1\le m$, $m^n-\lambda$
denotes the partition such that $(m^n-\lambda)_i = m-\lambda_{n+1-i}$. Note
that $(m^n-\lambda)'=n^m-\lambda'$.

The ring $\Lambda$ of symmetric functions will, unless otherwise noted,
have coefficients in $\F=\Q(q,t)$; we also consider the completion
$\hat\Lambda$ with respect to the natural grading (i.e., a sequence $f_k\to
0$ iff the degree of the lowest degree term of $f_k$ tends to $\infty$).
There are thus three natural types of basis: homogeneous bases,
inhomogeneous bases of $\Lambda$, and inhomogeneous bases of $\Lambda'$.
In the latter two cases, we require that the leading terms (the nonzero
homogeneous components of largest/smallest degree) form a homogeneous
basis.  If the polynomials $f_\lambda$ form a basis, we write $[f_\lambda]
g$ for the coefficient of $f_\lambda$ in the expansion of $g$ in terms of
that basis; this notation is mildly abusive, but the specific basis meant
should always be clear.  A similar notation applies to $BC_n$-symmetric
polynomials.  We will also need a notation for plethystic substitution in
symmetric functions; if $c_k$ is a sequence of elements of some
$\F$-algebra and $f\in\Lambda$, then $f([c_k])$ denotes the image of $f$
under the homomorphism $p_k\mapsto c_k$.

The formulas in the sequel will frequently involve products of the
form
\[
\prod_{(i,j)\in \lambda} f(i,j),
\]
where $(i,j)\in \lambda$ means that $1\le i$ and $1\le j\le
\lambda'_i$.  We will also have frequent recourse to $q$-symbols,
\[
(a;q)_n := \prod_{0\le i<n} (1-a q^i).
\]
If $n$ is omitted, the product is over all $i\ge 0$; also, $(x_1,x_2,\dots
x_m;q)_n:=\prod_{1\le i\le m} (x_i;q)_n$.  (A number of useful $q$-symbol
identities are given in \cite[Appendix I]{GasperG/RahmanM:1990}.)  The
following identity will be useful:

\begin{lem}
Let $a$, $b$, $c$, $d$, $e$ be arbitrary quantities.  Then
\begin{align}
\prod_{(i,j)\in \lambda} (1-a^{\lambda_i} b^j c^{\lambda'_j} d^i e)
&=
\prod_{1\le i\le j\le l} (a^{\lambda_i} b^{\lambda_{j+1}+1} c^j d^i e;b)_{\lambda_j-\lambda_{j+1}}\\
&=
\prod_{1\le i\le j\le l} (a^{\lambda_i} b^{\lambda_j} c^j d^i
e;b^{-1})_{\lambda_j-\lambda_{j+1}},
\end{align}
for any $l\ge \ell(\lambda)$.  If $|b|<1$, then
\[
\prod_{(i,j)\in\lambda} (1-a^{\lambda_i} b^j c^{\lambda'_j} d^i e)
=
\prod_{1\le i\le l}
\frac{(a^{\lambda_i} b c^l d^i e;b)}
{((ab)^{\lambda_i} b (cd)^i e;b)}
\prod_{1\le i<j\le l}
\frac{(a^{\lambda_i} b^{\lambda_j+1} c^{j-1} d^i e;b)}
{(a^{\lambda_i} b^{\lambda_j+1} c^j d^i e;b)},
\]
while if $|b|>1$, then
\[
\prod_{(i,j)\in\lambda} (1-a^{\lambda_i} b^j c^{\lambda'_j} d^i e)
=
\prod_{1\le i\le l}
\frac{((ab)^{\lambda_i} (cd)^i e;b^{-1})}
{(a^{\lambda_i} c^l d^i e;b^{-1})}
\prod_{1\le i<j\le l}
\frac{(a^{\lambda_i} b^{\lambda_j} c^j d^i e;b^{-1})}
{(a^{\lambda_i} b^{\lambda_j} c^{j-1} d^i e;b^{-1})}
\]
In particular,
\begin{align}
\prod_{(i,j)\in\lambda} (1-b^j d^i e)
&=
\prod_{1\le i\le l} (b d^i e;b)_{\lambda_i}\\
&=
\prod_{1\le i\le l} (b^{\lambda_i} d^i e;b^{-1})_{\lambda_i}.
\end{align}
\end{lem}

\begin{proof}
We find
{
\allowdisplaybreaks
\begin{align}
\prod_{(i,j)\in\lambda} (1-a^{\lambda_i} b^j c^{\lambda'_j} d^i e)
&=
\prod_{1\le i\le l}
\prod_{1\le j\le \lambda_i}
(1-a^{\lambda_i} b^j c^{\lambda'_j} d^i e)\\
&=
\prod_{1\le i\le k\le l}
\prod_{\substack{1\le j\le \lambda_i\\\lambda'_j=k}}
(1-a^{\lambda_i} b^j c^{\lambda'_j} d^i e)\\
&=
\prod_{1\le i\le k\le l}
\prod_{\lambda_{k+1}<j\le \lambda_k}
(1-a^{\lambda_i} b^j c^{\lambda'_j} d^i e)\\
&=
\prod_{1\le i\le k\le l} (a^{\lambda_i} b^{\lambda_k} c^k d^i e;b^{-1})_{\lambda_k-\lambda_{k+1}}\\
&=
\prod_{1\le i\le k\le l} (a^{\lambda_i} b^{\lambda_{k+1}+1} c^k d^i
e;b)_{\lambda_k-\lambda_{k+1}}.
\end{align}
}
If $|b|<1$, then
\begin{align}
\prod_{1\le i\le j\le l} (a^{\lambda_i} b^{\lambda_{j+1}+1} c^j d^i e;b)_{\lambda_j-\lambda_{j+1}}
&=
\prod_{1\le i\le j\le l}
\frac{(a^{\lambda_i} b^{\lambda_{j+1}+1} c^j d^i e;b)}
{(a^{\lambda_i} b^{\lambda_{j}+1} c^j d^i e;b)}\\
&=
\prod_{1\le i\le l}
\frac{(a^{\lambda_i} b c^l d^i e;b)}
{(a^{\lambda_i} b^{\lambda_i+1} (cd)^i e;b)}
\frac
{\prod_{1\le i\le j<l} (a^{\lambda_i} b^{\lambda_{j+1}+1} c^j d^i e;b)}
{\prod_{1\le i<j\le l} (a^{\lambda_i} b^{\lambda_{j}+1} c^j d^i e;b)}\\
&=
\prod_{1\le i\le l}
\frac{(a^{\lambda_i} b c^l d^i e;b)}
{(a^{\lambda_i} b^{\lambda_i+1} (cd)^i e;b)}
\prod_{1\le i<j\le l}
\frac
{\prod_{1\le i<j\le l} (a^{\lambda_i} b^{\lambda_{j}+1} c^{j-1} d^i e;b)}
{\prod_{1\le i<j\le l} (a^{\lambda_i} b^{\lambda_{j}+1} c^j d^i e;b)},
\end{align}
and similarly if $|b|>1$.
\end{proof}

\begin{rem}
This is essentially the argument used in section VI.6 of Macdonald for
the special case $a=q$, $b=1/q$, $c=t$, $d=1/t$.
\end{rem}

Three special cases are of particular importance; we define
{\allowdisplaybreaks
\begin{align}
C^+_\lambda(x;q,t)&:=\prod_{(i,j)\in \lambda} (1-q^{\lambda_i+j-1}
t^{2-\lambda'_j-i} x)\\*
&\phantom{:}=
\prod_{1\le i\le l}
\frac{(q^{\lambda_i} t^{2-l-i} x;q)}
{(q^{2\lambda_i} t^{2-2i} x;q)}
\prod_{1\le i<j\le l}
\frac{(q^{\lambda_i+\lambda_j} t^{3-i-j} x;q)}
{(q^{\lambda_i+\lambda_j} t^{2-i-j} x;q)},\\
C^-_\lambda(x;q,t)&:=\prod_{(i,j)\in \lambda} (1-q^{\lambda_i-j}
t^{\lambda'_j-i} x)\\*
&\phantom{:}=
\prod_{1\le i\le l}
\frac{(x;q)}
{(q^{\lambda_i} t^{l-i} x;q)}
\prod_{1\le i<j\le l}
\frac{(q^{\lambda_i-\lambda_j} t^{j-i} x;q)}
{(q^{\lambda_i-\lambda_j} t^{j-i-1} x;q)},\\
C^0_\lambda(x;q,t)&:=\prod_{(i,j)\in \lambda} (1-q^{j-1} t^{1-i} x)\\*
&\phantom{:}=
\prod_{1\le i\le l} (t^{1-i} x;q)_{\lambda_i}.
\end{align}
}
Thus, for instance, we have the following expressions for standard
Macdonald polynomial quantities in this notation.
\begin{align}
b_\lambda(q,t)&=\frac{C^-_\lambda(t;q,t)}{C^-_\lambda(q;q,t)}\\
P_\lambda(\left[\frac{1-u^k}{1-t^k}\right];q,t)
&=
\frac{t^{n(\lambda)} C^0_\lambda(u;q,t)}{C^-_\lambda(t;q,t)}\\
\langle P_\lambda,P_\lambda\rangle''_n
&=
\frac{C^0_\lambda(t^n;q,t)C^-_\lambda(q;q,t)}{C^0_\lambda(q t^{n-1};q,t)C^-_\lambda(t;q,t)}.
\end{align}
Again, by convention, multiple arguments before the semicolon indicate
a product; e.g.,
\[
C^+_\lambda(x,y;q,t) = C^+_\lambda(x;q,t)C^+_\lambda(y;q,t).
\]
We list some useful transformations of these quantities.

\begin{lem}
For any partition $\lambda$,
\begin{align}
C^+_\lambda(1/x;1/q,1/t)
&=
(-qx)^{-|\lambda|} t^{3n(\lambda)} q^{-3n(\lambda')} C^+_\lambda(x;q,t)\\
C^-_\lambda(1/x;1/q,1/t)
&=
(-1/x)^{|\lambda|} t^{-n(\lambda)} q^{-n(\lambda')} C^-_\lambda(x;q,t)\\
C^0_\lambda(1/x;1/q,1/t)
&=
(-1/x)^{|\lambda|} t^{n(\lambda)} q^{-n(\lambda')}
C^0_\lambda(x;q,t)
\end{align}
\end{lem}

\begin{lem}
For any partition $\lambda$,
\begin{align}
C^+_{\lambda'}(x;q,t) &= C^+_{\lambda}(qt x;1/t,1/q) \\
C^-_{\lambda'}(x;q,t) &= C^-_{\lambda}(x;t,q) \\
C^0_{\lambda'}(x;q,t) &= C^0_{\lambda}(x;1/t,1/q)
\end{align}
\end{lem}

\begin{lem}
For any integers $m,n\ge 0$ and partition $\lambda$ with $\ell(\lambda)\le
n$,
\begin{align}
\frac{C^+_{m^n+\lambda}(x;q,t)}{C^+_{m^n}(x;q,t)}
&=
\frac{C^0_\lambda(q^{2m} t^{1-n} x;q,t)
C^+_\lambda(q^{2m} x;q,t)}
{C^0_\lambda(q^{m} t^{1-n} x;q,t)}
\\
\frac{C^-_{m^n+\lambda}(x;q,t)}{C^-_{m^n}(x;q,t)}
&=
\frac{C^0_\lambda(q^m t^{n-1} x;q,t) C^-_\lambda(x;q,t)}
{C^0_\lambda(t^{n-1} x;q,t)}
\\
\frac{C^0_{m^n+\lambda}(x;q,t)}{C^0_{m^n}(x;q,t)}
&=
C^0_\lambda(q^m x;q,t).
\end{align}
\end{lem}

\begin{lem}
For any integers $m,n\ge 0$ and partition $\lambda\subset m^n$,
\begin{align}
\frac{C^+_{m^n-\lambda}(x;q,t)}{C^+_{m^n}(x;q,t)}
&=
\frac{
C^+_\lambda(q^{2m-1} t^{-2n+3} x;1/q,1/t)
C^0_\lambda(q^{m-1} t^{2-2n} x,q^{2m-1} t^{2-n} x;1/q,1/t)}
{
C^0_{2\lambda^2}(q^{2m-1} t^{2-2n} x;1/q,1/t)}
\\
\frac{C^-_{m^n-\lambda}(x;q,t)}{C^-_{m^n}(x;q,t)}
&=
\frac{C^-_\lambda(x;q,t)}{C^0_\lambda(t^{n-1} x;q,t) C^0_\lambda(q^{m-1} x;1/q,1/t)}
\\
\frac{C^0_{m^n-\lambda}(x;q,t)}{C^0_{m^n}(x;q,t)}
&=
\frac{1}{C^0_\lambda(q^{m-1} t^{1-n} x;1/q,1/t)}
\end{align}
\end{lem}

\begin{lem}
For any partition $\lambda$,
\begin{align}
C^+_{2\lambda}(x;q,t) &= C^+_\lambda(x,qx;q^2,t)\\
C^-_{2\lambda}(x;q,t) &= C^-_\lambda(x,xq;q^2,t)\\
C^0_{2\lambda}(x;q,t) &= C^0_\lambda(x,xq;q^2,t)\\
C^+_{\lambda^2}(x;q,t) &= C^+_\lambda(x/t,x/t^2;q,t^2)\\
C^-_{\lambda^2}(x;q,t) &= C^-_\lambda(x,xt;q,t^2)\\
C^0_{\lambda^2}(x;q,t) &= C^0_\lambda(x,x/t;q,t^2)
\end{align}
\end{lem}

We will also need analogous results for principal specializations of skew
Macdonald polynomials.

\begin{lem}
For any partitions $\kappa\subset\lambda$,
\begin{align}
P_{\lambda/\kappa}(\left[\frac{1-u^{-k}}{1-t^{-k}}\right];1/q,1/t)
&=
(t/u)^{|\lambda/\kappa|}
P_{\lambda/\kappa}(\left[\frac{1-u^k}{1-t^k}\right];q,t)\\
P_{\lambda'/\kappa'}(\left[\frac{1-u^k}{1-q^k}\right];t,q)
&=
(-u)^{|\lambda/\kappa|}
\frac{b_\lambda(q,t)}{b_\kappa(q,t)}
P_{\lambda/\kappa}(\left[\frac{1-u^{-k}}{1-t^k}\right];q,t).
\end{align}
Let $m,n\ge 0$ be integers.  If $\ell(\lambda)\le n$, then
\[
P_{(m^n+\lambda)/(m^n+\kappa)}(
\left[\frac{1-u^k}{1-t^k}\right];q,t)
=
\frac{C^0_\lambda(q^{m+1} t^{n-1},t^n;q,t)C^0_\kappa(q^m t^n,q t^{n-1};q,t)}
{C^0_\kappa(q^{m+1} t^{n-1},t^n;q,t)C^0_\lambda(q^m t^n,q t^{n-1};q,t)}
P_{\lambda/\kappa}(
\left[\frac{1-u^k}{1-t^k}\right];q,t).
\]
If also $\lambda_1\le m$, then
\[
P_{(m^n-\kappa)/(m^n-\lambda)}(
\left[\frac{1-u^k}{1-t^k}\right];q,t)
=
(q/t)^{|\lambda/\kappa|}
\frac{C^0_\lambda(q^{-m};q,t) C^0_\kappa(q^{1-m}/t;q,t)
b_\lambda(q,t)}
{C^0_\kappa(q^{-m};q,t) C^0_\lambda(q^{1-m}/t;q,t)
b_\kappa(q,t)}
P_{\lambda/\kappa}(\left[\frac{1-u^k}{1-t^k}\right];q,t).
\]
\end{lem}

\begin{proof}
The first two transformations are straightforward.
For the third transformation, we claim that in fact
\[
P_{(m^n+\lambda)/(m^n+\kappa)}(;q,t)
=
\frac{C^0_\lambda(q^{m+1} t^{n-1},t^n;q,t)C^0_\kappa(q^m t^n,q t^{n-1};q,t)}
{C^0_\kappa(q^{m+1} t^{n-1},t^n;q,t)C^0_\lambda(q^m t^n,q t^{n-1};q,t)}
P_{\lambda/\kappa}(;q,t)
\]
It suffices to compare coefficients of $P_\mu(;q,t)$ for $\ell(\mu)\le n$.
We have
\begin{align}
[P_\mu(;q,t)]
P_{(m^n+\lambda)/(m^n+\kappa)}(;q,t)
&=
[Q_{m^n+\lambda}(;q,t)]
(Q_\mu(;q,t) Q_{m^n+\kappa}(;q,t))\\
&=
\frac{
\langle P_{m^n+\lambda}(;q,t),
Q_\mu(;q,t) Q_{m^n+\kappa}(;q,t)\rangle''_n}
{
\langle P_{m^n+\lambda}(;q,t),
Q_{m^n+\lambda}(;q,t)\rangle''_n}\\
&=
\left(
\frac{C^0_\lambda(q^{m+1} t^{n-1},t^n;q,t)C^0_\kappa(q^m t^n,q t^{n-1};q,t)}
{C^0_\kappa(q^{m+1} t^{n-1},t^n;q,t)C^0_\lambda(q^m t^n,q t^{n-1};q,t)}\right)
\frac{
\langle P_\lambda(;q,t),
Q_\mu(;q,t) Q_\kappa(;q,t)\rangle''_n}
{
\langle P_\lambda(;q,t),
Q_\lambda(;q,t)\rangle''_n}\\
&=
\frac{C^0_\lambda(q^{m+1} t^{n-1},t^n;q,t)C^0_\kappa(q^m t^n,q t^{n-1};q,t)}
{C^0_\kappa(q^{m+1} t^{n-1},t^n;q,t)C^0_\lambda(q^m t^n,q t^{n-1};q,t)}
[P_\mu(;q,t)]
P_{\lambda/\kappa}(;q,t).
\end{align}

Similarly,
\[
P_{(m^n-\kappa)/(m^n-\lambda)}(;q,t)
=
(q/t)^{|\lambda/\kappa|}
\frac{C^0_\lambda(q^{-m};q,t) C^0_\kappa(q^{1-m}/t;q,t)
b_\lambda(q,t)}
{C^0_\kappa(q^{-m};q,t) C^0_\lambda(q^{1-m}/t;q,t)
b_\kappa(q,t)}
P_{\lambda/\kappa}(;q,t),
\]
and thus the fourth claim follows.
\end{proof}

\begin{cor}
For any integers $m,n\ge 0$ and partition $\lambda$ with $\ell(\lambda)\le
n$, $\lambda_n\ge m$,
\[
P_{\lambda/m^n}([(1-u^k)/(1-t^k)];q,t)
=
\lim_{Q\to q^m}
\frac{C^0_{m^n}(t^n,q t^{n-1}/Q;q,t)C^0_\lambda(q t^{n-1},u/Q;q,t)}
     {C^0_{m^n}(q t^{n-1},u/Q;q,t)C^0_\lambda(q t^{n-1}/Q;q,t)C^-_\lambda(t;q,t)}
\]
Similarly, if $\lambda\subset m^n$,
\[
\frac{
P_{m^n/\lambda}([(1-u^k)/(1-t^k)];q,t)}
{P_{m^n}([(1-u^k)/(1-t^k)];q,t)}
=
\frac{t^{n(\lambda)}
(q/u)^{|\lambda|} C^0_\lambda(t^n,q^{-m};q,t)}
{C^-_\lambda(q;q,t) C^0_\lambda(t^{n-1} q^{1-m}/u;q,t)}
\]
\end{cor}

\section{Interpolation polynomials}\label{sec:interp}

We define a ($q$-)difference operator $D^{(n)}(u_1,u_2;q,t)$ acting on
$\F[x_i^{\pm 1}]^{BC_n}$, as follows:

\begin{defn}\label{defn:diffop}
The operator $D^{(n)}(u_1,u_2;q,t)$ acts by:
\begin{align}
(D^{(n)}(u_1,u_2;q,t)f)&(x_1,x_2,\dots x_n)
=\\
&\sum_{\sigma\in \{\pm 1\}^n}
\prod_{1\le i\le n}
\frac{
(1-u_1 x_i^{\sigma_i})(1-u_2 x_i^{\sigma_i})
}{
(1-x_i^{2\sigma_i})
}
\prod_{1\le i<j\le n}
\frac{
(1-t x_i^{\sigma_i} x_j^{\sigma_j})
}
{
(1-x_i^{\sigma_i} x_j^{\sigma_j})
}
f(x_1 q^{\sigma_1/2},x_2 q^{\sigma_2/2}\dots x_n q^{\sigma_n/2})\notag
\end{align}
\end{defn}

\begin{rem}
This is one of the difference operators associated to the $BC/C$ Macdonald
polynomials; in particular, for each $u_1$, $u_2$, the eigenfunctions of
$D^{(n)}(u_1,u_2;q,t)$ are the Koornwinder polynomials
\[
K^{(n)}_\mu(;q,t;u_1,u_1\sqrt{q},u_2,u_2\sqrt{q}).
\]
See also Lemma \ref{lem:koorn:ortho1} below.
\end{rem}

Now, consider the ring $\F[s,1/s][x_i^{\pm 1}]^{BC_n}$, with basis of the
form $s^k m_\lambda(x)$ for $k\in \Z$ and $\lambda$ a partition (where
$m_\lambda(x)$ is the orbit sum of $\prod_i x_i^{\lambda_i}$); we extend
the dominance ordering to such ``monomials'' by taking
\[
(k,\lambda)\ge (l,\mu)
\]
when $\lambda\ge \mu$, $|l-k|\le |\lambda|-|\mu|$.
Now, define a difference operator $D^{(n)}_s(u;q,t)$ by:
\[
(D^{(n)}_s(u;q,t)f)(x_1,x_2,\dots x_n;s)
=
D^{(n)}(s,u/s;q,t) f(x_1,x_2,\dots x_n;s\sqrt{q}),
\]
thus acting on $s$ in addition to the $x$ variables.

\begin{lem}\label{lem:Dpoly}
The operator $D^{(n)}_s(u;q,t)$ gives a well-defined operator on
$\F[s,1/s][x_i^{\pm 1}]^{BC_n}$, and acts on monomials triangularly with
respect to the dominance ordering.  In particular,
\[
D^{(n)}_s(u;q,t) s^k m_\lambda = q^{k/2} E^{(n)}_\lambda(u;q,t) s^k m_\lambda
+ \text{dominated terms},
\]
with
\[
E^{(n)}_\lambda(u;q,t)=
q^{-|\lambda|/2} \prod_{1\le i\le n} (1-q^{\lambda_i} t^{n-i} u).
\]
\end{lem}

\begin{proof}
Since
\[
D^{(n)}_s(u;q,t) s^k f = q^{k/2} s^k D^{(n)}_s(u;q,t) f,
\]
it suffices to
consider the case $k=0$.  If we multiply $D^{(n)}_s(u;q,t) m_\lambda$ by
the $BC_n$-antisymmetric product
\[
\prod_{1\le i\le n} (x_i-1/x_i)
\prod_{1\le i<j\le n} (x_i+1/x_i-x_j-1/x_j),
\label{eq:CVandermonde}
\]
the result is manifestly a polynomial; moreover, we can use the symmetry of
$m_\lambda$ to write it as
\[
\sum_{\sigma\in \{\pm 1\}^n}
\Bigl(\prod_{1\le i\le n} \sigma_i R_{x_i}(\sigma_i)\Bigr)
F(u;x_1,x_2\dots x_n;s)
m_\lambda(\sqrt{q} x_1,\sqrt{q} x_2,\dots \sqrt{q} x_n;\sqrt{q} s),
\]
where $R_{x_i}(\pm 1)$ are the homomorphisms such that $R_{x_i}(\pm
1)x_j=x_j$, $R_{x_i}(\pm 1)x_i=x_i^{\pm 1}$, and where $F$ is a Laurent
polynomial in the ring
\[
\F[u,s,1/s][x_i^{\pm 1}]^{S_n}.
\]
Moreover, the monomials of $F$ are all dominated by the monomial
\[
x_1^n x_2^{n-1}\dots x_n + \cdots,
\]
and thus the monomials of $F m_\lambda$ are dominated by
\[
\prod_i x_i^{n+1-i+\lambda_i} + \cdots.
\]

Since the operators $\prod_i R_{x_i}(\sigma_i)$ form a normal subgroup
of $BC_n$, the resulting sum is $BC_n$-antisymmetric, and we may thus divide
back out the factor \eqref{eq:CVandermonde}, to obtain a $BC_n$-symmetric
polynomial dominated by $m_\lambda$.

For the leading coefficient, we compute:
\[
[\prod_i x_i^{\sigma_i(n+1-i)}] F(u;x_1,x_2,\dots x_n;s)
=
(-1)^n \prod_{1\le i\le n:\sigma_i=1} u t^{n-i},
\]
and thus
\begin{align}
[\prod_i x_i^{\lambda_i+n+1-i}]
\sum_{\sigma\in \{\pm 1\}^n}
\Bigl(\prod_{1\le i\le n}\sigma_i R_{x_i}(\sigma_i)\Bigr)
F m_\lambda
&=
\sum_{\sigma\in \{\pm 1\}^n}
(-1)^n
\Bigl(\prod_{1\le i\le n} \sigma_i q^{\sigma_i\lambda_i/2}\Bigr)
\prod_{1\le i\le n:\sigma_i=1} u t^{n-i}\\
&=
q^{-|\lambda|/2}
\prod_{1\le i\le n} (1-q^{\lambda_i} u t^{n-i}),
\end{align}
as required.
\end{proof}

Note that for generic $u$, the diagonal coefficients are all distinct, and
thus for each $(k,\lambda)$, there is a unique eigenfunction of
$D^{(n)}_s(u;q,t)$ of the form
$s^k m_\lambda+\text{dominated terms}$;
clearly multiplying such an eigenfunction by $s^j$ yields another
eigenfunction.  It turns out that these eigenfunctions are independent of
$u$, and are essentially just the $BC$-type interpolation polynomials of
Okounkov \cite{OkounkovA:1998}.  Given a polynomial $p\in
\F(s)[x_i^{\pm 1}]^{BC_n}$, we define
$p(\mu;s) := p(q^{\mu_i} t^{n-i} s;s)$.

\begin{thm}
The following three facts hold for all partitions $\lambda$.
\begin{itemize}
\item[a.] The operators $D^{(n)}_s(u;q,t)$ for different $u$ commute on the
monomials dominated by $m_\lambda$, and thus have a common eigenfunction
$s^k\bar{P}^{*(n)}_\mu(;s)$ for each leading monomial $s^k m_\mu$ dominated
by $m_\lambda$.
\item[b.] For any partition $\mu$, $\bar{P}^{*(n)}_\lambda(\mu;s)=0$
unless $\lambda\subset\mu$.
\item[c.] $\bar{P}^{*(n)}_\lambda(\lambda;s)\ne 0$.
\end{itemize}
\end{thm}

\begin{proof}
First, fix a partition $\lambda$, and suppose that $(a)$ holds for $\lambda$.
Then we claim that $(b)$ and $(c)$ hold as well.
Indeed, let $\mu$ be any partition different from $\lambda$; in
particular, fix $l$ such that $\mu_l\ne \lambda_l$.  By inspection of
the eigenvalues, the operator
\[
D^{(n)}_s(t^{l-n} q^{-\lambda_l};q,t)
\]
annihilates $\bar{P}^{*(n)}_\lambda$; on the other hand, we compute
\[
(D^{(n)}_s(t^{l-n} q^{-\lambda_l};q,t) \bar{P}^{*(n)}_\lambda)(\mu;s/\sqrt{q})
=
\sum_{\nu\prec\mu}
C_{\mu/\nu}
\bar{P}^{*(n)}_\lambda(\nu;s),
\]
for some coefficients $C_{\mu/\nu}$; the terms not corresponding to
partitions vanish, since then either $\sigma_n=-1$, $\nu_n=0$, and thus
$(1-s/x_n)=0$, or $\sigma_i=-1$, $\sigma_{i+1}=1$, $\nu_i=\nu_{i+1}$,
and thus $(1-t x_{i+1}/x_i)=0$.
Moreover, we compute
\[
C_{\mu/\mu}
=
\prod_{1\le i\le n}
\frac{(1-q^{\mu_i-1} t^{n-i} s^2)(1-q^{\mu_i-\lambda_l} t^{l-i})}
     {1-q^{2\mu_i-1} t^{2n-2i} s^2}
\prod_{1\le i<j\le n}
\frac{1-q^{\mu_i+\mu_j-1} t^{2n-i-j+1} s^2}
{1-q^{\mu_i+\mu_j-1} t^{2n-i-j} s^2}
\ne 0
\]

In other words, we can expand $\bar{P}^{*(n)}_\lambda(\mu;s)$ in terms of
the values $\bar{P}^{*(n)}_\lambda(\nu;s)$ for $\nu\subsetneq\mu$.
By induction, we find that $\bar{P}^{*(n)}_\lambda(\mu;s)=0$ whenever
$\mu$ does not contain $\lambda$, thus proving $(b)$.  Furthermore,
if $\bar{P}^{*(n)}_\lambda(\lambda;s)=0$, the induction would then prove
$\bar{P}^{*(n)}_\lambda(\mu;s)=0$ for all $\mu$, impossible since
$\bar{P}^{*(n)}_\lambda\ne 0$; we thus have $(c)$ as well.

Now, suppose that $(a)$ (and thus $(b)$ and $(c)$) holds for
all $\mu<\lambda$.  In particular, it follows that the operators
$D^{(n)}_s(u;q,t)$ commute on the space of polynomials {\it strictly} dominated
by $m_\lambda$, and it will thus suffice to show that they commute on
$m_\lambda$.  Let $f(x;s)$ be the unique $BC_n$-symmetric polynomial with
coefficients in $\F(s)$ of the form
\[
f(x;s)=m_\lambda(x)+\sum_{\mu<\lambda} c_{\lambda\mu}(s) \bar{P}^{*(n)}_\mu(x;s)
\]
such that $f(\mu;s)=0$ for $\mu<\lambda$; the resulting equations for
$c_{\lambda\mu}(s)$ are triangular with nonzero diagonal, by claims
$(b)$ and $(c)$.

Extend the action of $D^{(n)}_s(u;q,t)$ to polynomials with coefficients in
$k(s)$ in the obvious way.  The same inner induction as before proves
\[
(D^{(n)}_s(u;q,t) f)(\mu;s)=0
\]
whenever $\mu<\lambda$; it follows that $D^{(n)}_s(u;q,t)f$ is proportional to
$f$, and by comparing leading monomials, that:
\[
D^{(n)}_s(u;q,t) f=E^{(n)}_\lambda(u;q,t) f.
\]
We conclude that
{\allowdisplaybreaks
\begin{align}
D^{(n)}_s(u;q,t)D^{(n)}_s(v;q,t) m_\lambda
&=
D^{(n)}_s(u;q,t)D^{(n)}_s(v;q,t) (f-\sum_{\mu<\lambda} c_{\lambda\mu}(s)\bar{P}^{*(n)}_\mu)\\
&=
E^{(n)}_\lambda(u;q,t)E^{(n)}_\lambda(v;q,t) f
-
\sum_{\mu<\lambda} c_{\lambda\mu}(qs) E^{(n)}_\mu(u;q,t)E^{(n)}_\mu(v;q,t)
\bar{P}^{*(n)}_\mu\\
&=
D^{(n)}_s(v;q,t)D^{(n)}_s(u;q,t) m_\lambda,
\end{align}
}
as required.
\end{proof}

In the sequel, we will write the common eigenfunctions of $D^{(n)}_s(u;q,t)$ as
$\bar{P}^{*(n)}_\lambda(;q,t,s)$, and refer to them as interpolation
polynomials.

\begin{rem}
In particular, it follows that these polynomials agree up to a factor in
$\F(s)$ and a shifting of the arguments with the interpolation polynomials
of \cite{OkounkovA:1998}.  It was shown there that the interpolation
polynomials are not common eigenfunctions of any rational difference
operators in $x_1\dots x_n$; the loophole, of course, is that our
operators also act on $s$.
\end{rem}

\begin{rem}
The interpolation polynomials are, naturally, independent of the choice of
square root of $q$ used to define $D^{(n)}_s(u;q,t)$; the proof shows them
to be characterized by the triangularity and vanishing properties, neither
of which depends on that choice.
\end{rem}

\begin{cor}
The operators $D^{(n)}(u_1,u_2;q,t)$ satisfy the quasi-commutation relation
\[
D^{(n)}(u_1,\sqrt{q} u_2;q,t)D^{(n)}(\sqrt{q} u_1,u_3;q,t)
=
D^{(n)}(u_1,\sqrt{q} u_3;q,t)D^{(n)}(\sqrt{q} u_1,u_2;q,t)
\]
\end{cor}

\begin{proof}
For any function $f\in \F(s)[x_i^{\pm 1}]^{BC_n}$, we have:
\[
(D^{(n)}_s(u;q,t)D^{(n)}_s(v;q,t) f)(x_1,x_2,\dots x_n;s)
=
(D^{(n)}(s,u/s;q,t)D^{(n)}(s\sqrt{q},\frac{v}{s\sqrt{q}};q,t) f)(x_1,x_2,\dots x_n;s q);
\]
as this must be symmetric in $u$ and $v$, the corrolary follows.
\end{proof}

\begin{cor}
The interpolation polynomials satisfy the difference equation
\begin{align}
q^{-|\lambda|/2} &\Bigl(\!\prod_{1\le i\le n} (1-q^{\lambda_i} t^{n-i} u)\Bigr)
\bar{P}^{*(n)}_\lambda(x_1,x_2,\dots x_n;q,t,s)
=\\
&\sum_{\sigma\in \{\pm 1\}^n}
\prod_{1\le i\le n}
\frac{
(1-s x_i^{\sigma_i})(1-u x_i^{\sigma_i}/s)
}{
(1-x_i^{2\sigma_i})
}
\prod_{1\le i<j\le n}
\frac{
(1-t x_i^{\sigma_i} x_j^{\sigma_j})
}
{
(1-x_i^{\sigma_i} x_j^{\sigma_j})
}
\bar{P}^{*(n)}_\lambda(x_1 q^{\sigma_1/2},x_2 q^{\sigma_2/2}\dots x_n q^{\sigma_n/2};q,t,s\sqrt{q})\notag
\end{align}
\end{cor}

Two limiting cases of the theorem are of special interest.  First, we find
that the limit
\[
\bar{P}^{*(n)}_\lambda(x_i;q,t)
=
\lim_{s\to\infty} s^{-|\lambda|} \bar{P}^{*(n)}_\lambda(x_i s;q,t,s)
\]
is well-defined, and produces an ordinary symmetric polynomial vanishing
when $x_i=q^{\mu_i} t^{n-i}$ for $\mu\not\subset\lambda$.  We thus recover
the symmetric version of the shifted Macdonald polynomials.  In this limit,
our difference operators converge to those of Knop and Sahi
\cite{KnopF:1997,SahiS:1996}.  In particular, we see immediately that the
shifted Macdonald polynomials are limits of $BC/C$-type Macdonald
polynomials in multiple ways \cite{BakerTH/ForresterPJ:2000}, as they are
eigenfunctions of limiting versions of $D^{(n)}(u_1,u_2;q,t)$; the action
on $s$ becomes trivial in the limit.  Similarly, the ``leading term'' limit
\[
\lim_{a\to\infty} a^{-|\lambda|} \bar{P}^{*(n)}_\lambda(a x_i;q,t,s)
\]
satisfies the difference equation of the ordinary Macdonald polynomials,
and thus
\[
\lim_{a\to\infty} a^{-|\lambda|} \bar{P}^{*(n)}_\lambda(a x_i;q,t,s)
=
P_\lambda(x_i;q,t).
\]
We will prove a refinement of this fact in Theorem \ref{thm:PsPtriangular}
below.

We will now recall some basic properties of the interpolation polynomials.
First, some symmetries:

\begin{lem}\cite{OkounkovA:1998}
The interpolation polynomials satisfy the identities
\begin{align}
\bar{P}^{*(n)}_\lambda(x_1,\dots x_n;q,t,s)
&=
\bar{P}^{*(n)}_\lambda(x_1,\dots x_n;1/q,1/t,1/s)\\
&=
(-1)^{|\lambda|}
\bar{P}^{*(n)}_\lambda(-x_1,\dots -x_n;q,t,-s)\label{eq:eveninterp}
\end{align}
for all partitions $\lambda$.
\end{lem}

Next, identities for shrinking the dimension or the indexing partition:

\begin{lem}\label{lem:interp:dec_mn}\cite{OkounkovA:1998}
For any partition $\lambda$,
\[
\bar{P}^{*(n+m)}_\lambda(x_1,\dots x_n,s,st,\dots st^{m-1};q,t,s)
=
\begin{cases}
\bar{P}^{*(n)}_\lambda(x_1,\dots x_n;q,t,st^m)&\lambda_{n+1}=0\\
0 & \lambda_{n+1}>0
\end{cases}
\]
and
\[
\bar{P}^{*(n)}_{m^n+\lambda}(x_1,x_2\dots x_n;q,t,s)
=
\prod_{(i,j)\in m^n} (x_i+1/x_i-q^{j-1} s-q^{1-j}/s)
\bar{P}^{*(n)}_{\lambda}(x_1,x_2,\dots x_n;q,t,sq^m).
\]
\end{lem}

\begin{proof}
In each case, both sides are monic, triangular, and vanish on the
appropriate partitions, so must be the same.
\end{proof}

The next corollary then follows by induction.

\begin{cor}
We have the normalization
\[
\bar{P}^{*(n)}_\lambda(\lambda;q,t,s)
=
(q t^{n-1} s)^{-|\lambda|} t^{n(\lambda)} q^{-2n(\lambda')}
C^-_\lambda(q;q,t)
C^+_\lambda(t^{2n-2} s^2;q,t)
\]
\end{cor}

\begin{rems}
In particular, we have:
\[
P^{*(n)}_\lambda(x;q,t,s)
=
(t^{n-1} s)^{-|\lambda|} \bar{P}^{*(n)}_\lambda(x_i t^{n-i} s;q,t,s),
\]
where $P^{*(n)}_\lambda$ is Okounkov's interpolation polynomial in $n$
variables; this also follows by comparing leading terms.
\end{rems}

\begin{rems}
If we compute $\bar{P}^{*(n)}_\lambda$ by solving the appropriate
triangular system of equations, we find that the denominators of the
coefficients of $\bar{P}^{*(n)}_\lambda$ must divide the determinant of the
system, i.e.
\[
\prod_{\mu<\lambda} \bar{P}^{*(n)}_\mu(\mu;q,t,s).
\]
Since the coefficients of $\bar{P}^{*(n)}_\lambda$ are in $\F[s,1/s]$, we
conclude that the only possible denominator factors are $q$, $t$, $s$, and
$(1-q^i t^j)$ for $i,j\ge 0$.
\end{rems}

We also note the following special case of the difference equation;
this is of interest because the difference operator involved is independent
of $s$.

\begin{cor}
The interpolation polynomials satisfy the difference equation
\begin{align}
\sum_{\sigma\in \{\pm 1\}^n}
\prod_{1\le i\le n}
\frac{
x_i^{\sigma_i}
}{
(1-x_i^{2\sigma_i})
}
&\prod_{1\le i<j\le n}
\frac{
(1-t x_i^{\sigma_i} x_j^{\sigma_j})
}
{
(1-x_i^{\sigma_i} x_j^{\sigma_j})
}
\bar{P}^{*(n)}_\lambda(x_1 q^{\sigma_1/2},x_2 q^{\sigma_2/2}\dots x_n
q^{\sigma_n/2};q,t,s)
\notag\\*
&=
\begin{cases}
q^{-|\lambda|/2}
\Bigl(\!\prod_{1\le i\le n} (1-q^{\lambda_i} t^{n-i})\Bigr)
\bar{P}^{*(n)}_{\mu}(x_1,x_2,\dots x_n;q,t,s\sqrt{q})&\lambda=1^n+\mu,\\
0&\text{otherwise.}
\end{cases}
\end{align}
\end{cor}

\begin{proof}
Take $u=1$ in the full difference equation, then divide both sides by
$\prod_{1\le i\le n} (x_i+1/x_i-s-1/s)$.
\end{proof}

The branching rule for interpolation polynomials extends to the following
``bulk'' branching rule.

\begin{thm}\label{thm:interp:bulk_branch}
We have
\[
\bar{P}^{*(n+m)}_\lambda(x_1,x_2,\dots x_n,t^{m-1} v,t^{m-2} v,\dots v;q,t,s)
=
\sum_{\substack{\mu\subset\lambda\\\ell(\mu)\le n}}
\psi^{(B)}_{\lambda/\mu}(v,vt^m;q,t,s t^n)
\bar{P}^{*(n)}_\mu(x_1,x_2,\dots x_n;q,t,s),
\]
where
\[
\psi^{(B)}_{\lambda/\mu}(v,v';q,t,s)
=
\frac{
C^0_\lambda(s/v;q,t)
C^0_\lambda(t/sv';1/q,1/t)}
{
C^0_\mu(s/v;q,t)
C^0_\mu(t/sv';1/q,1/t)}
P_{\lambda/\mu}(\left[\frac{v^k-v^{\prime k}}{1-t^k}\right];q,t)
\]
\end{thm}

\begin{proof}
Apply the difference equation with $u=vs$.  The only terms that contribute
are those with $\sigma_{n+i}=1$ for $1\le i\le m$, in which case the
difference operator simplifies to an $n$-dimensional operator.  We thus
find
\begin{align}
q^{-|\lambda|/2}
\prod_{1\le i\le n+m} &(1-q^{\lambda_i} t^{n+m-i} v s)
\bar{P}^{*(n+m)}_\lambda((x_i),(t^{m-j} v);q,t,s)
=\label{eq:interp:branch_deq}\\*
&\prod_{n<i\le n+m} (1-t^{n+m-i} v s)
D^{(n)}(t^m v,s;q,t)
\bar{P}^{*(n+m)}_\lambda((x_i),(t^{m-j} v\sqrt{q});q,t,s\sqrt{q}).
\notag
\end{align}
Now, we have an expansion of the form
\[
\bar{P}^{*(n+m)}_\lambda(x_i,v;q,t,s)
=
\sum_{\substack{\mu\le\lambda\\\ell(\mu)\le n}}
c_{\lambda\mu}(v;q,t,s)
\bar{P}^{*(n)}_\mu(x_i;q,t,s),
\]
with unknown coefficients $c_{\lambda\mu}$, since the interpolation
polynomials are monic and triangular.  As
\[
D^{(n)}(t^m v,s)
\bar{P}^{*(n)}_\mu(x_i;q,t,s\sqrt{q})
=
q^{-|\mu|/2} \prod_{1\le i\le n} (1-q^{\mu_i} t^{n+m-i} v s)
\bar{P}^{*(n)}_\mu(x_i;q,t,s),
\]
we can substitute this expansion into \eqref{eq:interp:branch_deq}, and
compare coefficients of $\bar{P}^{*(n)}_\mu(x_i;q,t,s)$.  We thus find
\[
q^{-|\lambda|/2}
\prod_{1\le i\le n+m} (1-q^{\lambda_i} t^{n+m-i} v s)
c_{\lambda\mu}(v;q,t,s)
=
q^{-|\mu|/2}
\prod_{1\le i\le n+m} (1-q^{\mu_i} t^{n+m-i} v s)
c_{\lambda\mu}(v\sqrt{q};q,t,s\sqrt{q}),
\]
and by symmetry ($v\mapsto 1/(t^{m-1} v)$),
\[
q^{-|\lambda|/2}
\prod_{1\le i\le n+m} (1-q^{\lambda_i} t^{n+1-i}s/v)
c_{\lambda\mu}(v;q,t,s)
=
q^{-|\mu|/2}
\prod_{1\le i\le n+m} (1-q^{\mu_i} t^{n+1-i}s/v)
c_{\lambda\mu}(v/\sqrt{q};q,t,s\sqrt{q}).
\]
Solving these difference equations gives
\[
c_{\lambda\mu}(v;q,t,s)
=
c_{\lambda\mu}(q,t)
\frac{\prod_{(i,j)\in \lambda}
(-s)^{-1} q^{1-j} t^{i-n-m}
(1-q^{j-1} t^{n+1-i} s/v)
(1-q^{j-1} t^{n+m-i} s v)}
{\prod_{(i,j)\in \mu}
(-s)^{-1} q^{1-j} t^{i-n-m}
(1-q^{j-1} t^{n+1-i} s/v)
(1-q^{j-1} t^{n+m-i} s v)},
\]
where $c_{\lambda\mu}(q,t)$ is independent of $s$ and $v$.
This remaining factor is then determined by the Macdonald polynomial
limit; we find that
\[
c_{\lambda\mu}(q,t)
=
\lim_{v\to\infty} v^{-|\lambda|+|\mu|}
c_{\lambda\mu}(v;q,t,s)
=
P_{\lambda/\mu}(1,t,\dots t^{m-1};q,t),
\]
as required.
\end{proof}

\begin{cor}\cite{OkounkovA:1998}
We have
\[
\bar{P}^{*(n+1)}_\lambda(x_1,x_2,\dots x_n,v;q,t,s)
=
\sum_{\substack{\mu'\prec\lambda'\\\mu_{n+1}=0}}
\psi^{(b)}_{\lambda/\mu}(v;q,t,s t^n)
\bar{P}^{*(n)}_\mu(x_1,x_2,\dots x_n;q,t,s),
\]
where
\[
\psi^{(b)}_{\lambda/\mu}(v;q,t,s)
=
\psi_{\lambda/\mu}(q,t)
\prod_{(i,j)\in \lambda/\mu}
(v+1/v-q^{j-1} t^{1-i}s-q^{1-j}t^{i-1}/s)
\]
\end{cor}

\begin{cor}\label{cor:interp:rect1}
For any partition $\lambda$ and integer $n\ge\ell(\lambda)$,
\[
\bar{P}^{*(n)}_\lambda(xt^{n-1}s,x t^{n-2}s,\dots xs;q,t,s)
=
(-s t^{n-1})^{-|\lambda|}
t^{2n(\lambda)} q^{-n(\lambda')}
\frac{C^0_\mu(t^n,1/x,x s^2 t^{n-1};q,t)}{C^-_\mu(t;q,t)}
\]
\end{cor}

\begin{rem}
Equivalently, this gives a formula for evaluating Okounkov's original
version of the interpolation polynomials at a constant.
\end{rem}

The bulk branching rule also implies the following connection coefficient
result.

\begin{thm}
For any partitions $\lambda$, $\mu$,
\[
[\bar{P}^{*(n)}_\mu(;q,t,s)]
\bar{P}^{*(n)}_\lambda(;q,t,s')
=
\frac{
C^0_\lambda(t^n;q,t)
C^0_\lambda(t^{1-n}/s s';1/q,1/t)}
{
C^0_\mu(t^n;q,t)
C^0_\mu(t^{1-n}/s s';1/q,1/t)}
P_{\lambda/\mu}([\frac{s^k-s^{\prime k}}{1-t^k}];q,t).
\]
\end{thm}

\begin{proof}
Both sides are rational in $q$, $t$, $s$, $s'$, so it suffices to prove
this under the assumption $s'=s t^m$ for some integer $m\ge 0$.
Then by Lemma \ref{lem:interp:dec_mn} and the bulk branching rule,
\begin{align}
\bar{P}^{*(n)}_\lambda(x_1,\dots x_n;q,t,s t^m)
&=
\bar{P}^{*(n+m)}_\lambda(x_1,\dots x_n,s,st,st^2,\dots st^{m-1};q,t,s)\\
&=
\sum_{\substack{\mu\subset\lambda\\\ell(\mu)\le n}}
\psi^{(B)}_{\lambda/\mu}(s,s';q,t,s t^n)
\bar{P}^{*(n)}_\mu(x_1,x_2,\dots x_n;q,t,s).
\end{align}
The theorem follows immediately.
\end{proof}

\begin{rem}
If we expand an interpolation polynomial using the connection coefficient
identity, we cannot in general insist that the polynomials on both sides
are evaluated at a partition.  A notable exception is when $s s' = t^{n-1}
q^{-m}$, since then $\bar{P}^{*(n)}(\mu;q,t,s)$ and
$\bar{P}^{*(n)}(m^n-\mu;q,t,s')$ are evaluated at the same point.
\end{rem}

There is also a bulk version of the Pieri identity.

\begin{thm}\label{thm:interp:bulk_Pieri}
For any integer $n\ge 0$ and partition $\mu$ of length $\le n$, the
following identity holds in the power series ring
$\F(s)[x_1^{\pm 1},\dots x_n^{\pm 1}][[u,v]]$.
\begin{align}
\prod_{1\le i\le n} \frac{(v x_i,v/x_i;q)}{(u x_i,u/x_i;q)}
\bar{P}^{*(n)}_\mu(x_1,\dots x_n;q,t,s)
&=
\prod_{1\le i\le n}
\frac{(v q^{\mu_i} t^{n-i} s,v q^{-\mu_i} t^{i-n}/s;q)}
{(u q^{\mu_i} t^{n-i} s,u q^{-\mu_i} t^{i-n}/s;q)}\notag\\*
&\phantom{{}={}}
\sum_{\lambda\supset\mu}
\psi^{(P)}_{\lambda/\mu}(u,v;q,t;s t^n)
\bar{P}^{*(n)}_\lambda(x_1,\dots x_n;q,t,s),
\end{align}
where
\[
\psi^{(P)}_{\lambda/\mu}(u,v;q,t;s)=
\frac{C^0_\mu(sv/t;q,t)C^0_\mu(tu/qs;1/q,1/t)}
{C^0_\lambda(sv/t;q,t)C^0_\lambda(tu/qs;1/q,1/t)}
Q_{\lambda/\mu}([(u^k-v^k)/(1-t^k)];q,t).
\]
\end{thm}

\begin{proof}
It suffices to consider the case $u=q^m v$ for an integer $m\ge 0$.
We certainly have an expansion of the form
\[
\prod_{1\le i\le n} (v x_i,v/x_i;q)_m
\bar{P}^{*(n)}_\mu(x_1,\dots x_n;q,t,s)
=
\sum_\lambda c^{(n,m)}_{\lambda\mu}(v;q,t,s)
\bar{P}^{*(n)}_\lambda(x_1,\dots x_n;q,t,s),
\label{eq:bPieri:teq1}
\]
for $c^{(n,m)}_{\lambda\mu}(v;q,t,s)\in \F[s,1/s,v]$; the content of the
theorem is that
\[
c^{(n,m)}_{\lambda\mu}(v;q,t,s)
=
Q_{\lambda/\mu}([\frac{q^{mk}-1}{1-t^k}];q,t)
v^{|\lambda|-|\mu|}
\prod_{(i,j)\in (m^n+\mu)/\lambda}
(1-v q^{j-1} t^{n-i} s)
(1-v q^{m-j} t^{i-n}/s)
.
\]
By the Macdonald limit,
\begin{align}
\lim_{v\to 0}
v^{|\mu|-|\lambda|} c^{(n,m)}_{\lambda\mu}(v;q,t,s)
&=
Q_{\lambda/\mu}([\frac{q^{mk}-1}{1-t^k}];q,t)\\
\lim_{v\to \infty}
v^{|\lambda|-|\mu|-2mn}
c^{(n,m)}_{\lambda\mu}(v;q,t,s)
&=
q^{(m-1)(mn+|\mu|-|\lambda|)}
Q_{\lambda/\mu}([\frac{q^{mk}-1}{1-t^k}];q,t),
\end{align}
and thus $v^{|\mu|-|\lambda|} c^{(n,m)}_{\lambda\mu}(v;q,t,s)$
is a polynomial of degree at most $2mn+2|\mu|-2|\lambda|$, with
constant term $Q_{\lambda\mu}([(q^{mk}-1)/(1-t^k)];q,t)$.

Now, if we evaluate both sides of \eqref{eq:bPieri:teq1} at $\kappa$, the
left-hand side vanishes if either $\mu\not\subset\kappa$ or
\[
v\in \{t^{n-i} q^{\kappa_i+j-m} s,t^{i-n} q^{1-j-\kappa_i}/s:
(i,j)\in (m^n+\kappa)/\kappa
\}.
\]
Thus by induction in $\lambda$, we find $v^{|\mu|-|\lambda|}
c^{(n,m)}_{\lambda\mu}(v;q,t,s)$ vanishes if $\mu\not\subset\lambda$, and
is otherwise a multiple of
\[
\prod_{\substack{(i,j)\in m^n+\mu\\(i,j)\notin\lambda}}
(1-v q^{j-1} t^{n-i} s)
(1-v q^{m-j} t^{i-n}/s).
\]
This has degree $\ge 2mn+2|\mu|-2|\lambda|$, and thus it follows
immediately that
\[
v^{|\mu|-|\lambda|}
c^{(n,m)}_{\lambda\mu}(v;q,t,s)
=
Q_{\lambda/\mu}([\frac{q^{mk}-1}{1-t^k}];q,t)
\prod_{\substack{(i,j)\in m^n+\mu\\(i,j)\notin\lambda}}
(1-v q^{j-1} t^{n-i} s)
(1-v q^{m-j} t^{i-n}/s)
\]
as required; note that the skew Macdonald polynomial vanishes unless 
$\lambda\subset m^n+\mu$.
\end{proof}

\begin{rem} When $m=1$, this is essentially the proof of
\cite{OkounkovA:1998} for the ordinary Pieri identity.
\end{rem}

The case $u=qv$ gives the ordinary $e$-type Pieri identity.

\begin{cor}\cite{OkounkovA:1998}
For any integer $n\ge 0$ and partition $\mu$ of length $\le n$,
\begin{align}
\prod_{1\le i\le n} (v+1/v+{}&x_i+1/x_i)
\bar{P}^{*(n)}_\mu(x_1,\dots x_n;q,t,s)\\*
&=
\prod_{1\le i\le n} (v+1/v+q^{\mu_i} t^{n-i} s+q^{-\mu_i} t^{i-n}/s)
\sum_{\lambda\succ\mu}
\psi^{(e)}_{\lambda/\mu}(v;q,t,s t^n)
\bar{P}^{*(n)}_\lambda(x_1,\dots x_n;q,t,s),
\notag\end{align}
where
\[
\psi^{(e)}_{\lambda/\mu}(v;q,t,s)
=
\psi'_{\lambda/\mu}(q,t)
\prod_{(i,j)\in \lambda/\mu}
(v+1/v+q^{j-1}t^{-i} s+q^{1-j} t^{i}/s)^{-1}.
\]
\end{cor}

We also obtain a $g$-type Pieri identity, by taking $v=tu$.

\begin{cor}
For any integer $n\ge 0$ and partition $\mu$ of length $\le n$, the
following identity holds in the power series ring
$\F(s)[x_1^{\pm 1},\dots x_n^{\pm 1}][[u]]$.
\begin{align}
\prod_{1\le i\le n} \frac{(tu x_i,tu/x_i;q)}{(u x_i,u/x_i;q)}
\bar{P}^{*(n)}_\mu(x_1,\dots x_n;q,t,s)
&=
\prod_{1\le i\le n}
\frac{(tu q^{\mu_i} t^{n-i} s,tu q^{-\mu_i} t^{i-n}/s;q)}
{(u q^{\mu_i} t^{n-i} s,u q^{-\mu_i} t^{i-n}/s;q)}\notag\\*
&\phantom{{}={}}
\sum_{\lambda\supset\mu}
\varphi^{(g)}_{\lambda/\mu}(u;q,t;s t^n)
\bar{P}^{*(n)}_\lambda(x_1,\dots x_n;q,t,s),
\end{align}
where
\[
\varphi^{(g)}_{\lambda/\mu}(u;q,t;s)=
\varphi_{\lambda/\mu}(q,t)
\prod_{(i,j)\in \lambda/\mu}
\frac{u}{(1-u q^{j-1} t^{1-i} s)(1-u q^{-j} t^i/s)}
\]
\end{cor}

As observed in \cite{OkounkovA:1998}, the branching rule and Pieri 
identity are connected via the Cauchy identity.

\begin{thm}\cite{OkounkovA:1998}
For any integers $m$, $n$,
\[
\sum_{\lambda\subset m^n}
(-1)^{mn-|\lambda|}
\bar{P}^{*(n)}_\lambda(x;q,t,s)
\bar{P}^{*(m)}_{n^m-\lambda'}(y;t,q,s)
=
\prod_{1\le i\le n}
\prod_{1\le j\le m}
(x_i+1/x_i-y_j-1/y_j)
\]
\end{thm}

\begin{proof}
By induction in $n$; if we multiply both sides by
\[
\prod_{1\le j\le m}
(x_{n+1}+1/x_{n+1}-y_j-1/y_j),
\]
then expanding via the ($e$-type) Pieri identity and simplifying via the
branching rule turns the resulting left-hand side into the left-hand
side of the Cauchy identity for $n+1$.
\end{proof}

\section{Hypergeometric transformations}\label{sec:hyperg}

Define the ``binomial coefficients''
\begin{align}
\binomQ{\lambda}{\mu}_{q,t,s}
&:=
\frac{\bar{P}^{*(n)}_\mu(\lambda;q,t,t^{1-n} s)}
     {\bar{P}^{*(n)}_\mu(\mu;q,t,t^{1-n} s)}\\
\binomI{\lambda}{\mu}_{q,t,s}
&:=
\frac{
\langle P_\mu,P_\mu\rangle''_n
\bar{P}^{*(n)}_{\lambda}(\lambda;q,t,t^{1-n} s)
\bar{P}^{*(n)}_{m^n-\lambda}(m^n-\mu;q,t,q^{-m}/s)
}{
\langle P_\lambda,P_\lambda\rangle''_n
\bar{P}^{*(n)}_{\mu}(\mu;q,t,t^{1-n} s)
\bar{P}^{*(n)}_{m^n-\mu}(m^n-\mu;q,t,q^{-m}/s)
},
\end{align}
where $m$ and $n$ are chosen so that $\lambda,\mu\subset m^n$.  Note that
each binomial coefficient vanishes unless $\mu\subset \lambda$, and is
equal to 1 when $\mu=\lambda$; furthermore, each is preserved by the
substitutions $s\mapsto -s$ and $(q,t,s)\mapsto (1/q,1/t,1/s)$.  That the
first binomial coefficient is independent of $n$ and the second is
independent of $m$ follows from Lemma \ref{lem:interp:dec_mn}; that the
second is also independent of $n$ follows from Theorem
\ref{thm:interp:inversion} below.  First, though, some transformations and
special values.

\begin{prop}
For any partitions $\mu\subset\lambda$,
\begin{align}
\binomQ{\lambda}{\mu}_{q,t,s} &= \binomQ{\lambda}{\mu}_{1/q,1/t,1/s}
&\binomI{\lambda}{\mu}_{q,t,s} &= \binomI{\lambda}{\mu}_{1/q,1/t,1/s}\\
\binomQ{\lambda}{\mu}_{q,t,s} &= \binomQ{\lambda}{\mu}_{q,t,-s}
&\binomI{\lambda}{\mu}_{q,t,s} &= \binomI{\lambda}{\mu}_{q,t,-s}
\end{align}
For integers $m,n\ge 0$ with $\ell(\lambda)\le n$,
\begin{align}
\binomQ{m^n+\lambda}{m^n+\mu}_{q,t,s} &=
q^{-m|\lambda/\mu|}
\frac{C^0_\lambda(t^{1-n} q^{2m} s^2,t^{n-1} q^{m+1};q,t)
C^0_\mu(t^{1-n} q^m s^2,t^{n-1} q;q,t)}
{C^0_\mu(t^{1-n} q^{2m} s^2,t^{n-1} q^{m+1};q,t)
C^0_\lambda(t^{1-n} q^m s^2,t^{n-1} q;q,t)}
\binomQ{\lambda}{\mu}_{q,t,q^m s}\\
\binomI{m^n+\lambda}{m^n+\mu}_{q,t,s} &=
q^{-m|\lambda/\mu|}
\frac{C^0_\lambda(t^{1-n} q^{2m} s^2,t^{n-1} q^{m+1};q,t)
C^0_\mu(t^{1-n} q^m s^2,t^{n-1} q;q,t)}
{C^0_\mu(t^{1-n} q^{2m} s^2,t^{n-1} q^{m+1};q,t)
C^0_\lambda(t^{1-n} q^m s^2,t^{n-1} q;q,t)}
\binomI{\lambda}{\mu}_{q,t,q^m s}.
\end{align}
If in fact $\lambda\subset m^n$,
\begin{align}
\binomQ{m^n-\mu}{m^n-\lambda}_{q,t,s}
&=
\frac{
\langle P_\lambda,P_\lambda\rangle''_n
\bar{P}^{*(n)}_{\mu}(\mu;q,t,q^{-m}/s)
\bar{P}^{*(n)}_{m^n-\mu}(m^n-\mu;q,t,t^{1-n} s)
}{
\langle P_\mu,P_\mu\rangle''_n
\bar{P}^{*(n)}_{\lambda}(\lambda;q,t,q^{-m}/s)
\bar{P}^{*(n)}_{m^n-\lambda}(m^n-\lambda;q,t,t^{1-n} s)
}
\binomI{\lambda}{\mu}_{q,t,t^{n-1}/q^m s}\\
\binomI{m^n-\mu}{m^n-\lambda}_{q,t,s}
&=
\frac{
\langle P_\lambda,P_\lambda\rangle''_n
\bar{P}^{*(n)}_{\mu}(\mu;q,t,q^{-m}/s)
\bar{P}^{*(n)}_{m^n-\mu}(m^n-\mu;q,t,t^{1-n} s)
}{
\langle P_\mu,P_\mu\rangle''_n
\bar{P}^{*(n)}_{\lambda}(\lambda;q,t,q^{-m}/s)
\bar{P}^{*(n)}_{m^n-\lambda}(m^n-\lambda;q,t,t^{1-n} s)
}
\binomQ{\lambda}{\mu}_{q,t,t^n/q^m s}
\end{align}
For any partition $\lambda\subset m^n$,
\begin{align}
\binomQ{\lambda}{0}_{q,t,s} &= 1\\
\binomI{\lambda}{0}_{q,t,s} &=
(-1)^{|\lambda|} t^{n(\lambda)} q^{-n(\lambda')}
\frac{C^+_\lambda(s^2;q,t)}
{C^0_\lambda(q s^2;q,t)}\\
\binomQ{m^n}{\lambda}_{q,t,s} &=
(-q)^{|\lambda|} t^{n(\lambda)} q^{n(\lambda')}
\frac{C^0_\lambda(t^n,q^{-m},q^m s^2/t^{n-1};q,t)}
{C^-_\lambda(q,t;q,t)C^+_\lambda(s^2;q,t)}
\\
\binomI{m^n}{\lambda}_{q,t,s} &=
\frac{
(-1)^{mn} t^{n(m^n)} C^0_{m^n}(q^m s^2/t^{n-1};q,t)
(q^m/t^{n-1})^{|\lambda|}
t^{2n(\lambda)}
C^0_\lambda(t^n,q^{-m};q,t) C^0_{2\lambda^2}(s^2q;q,t)}
{q^{n(n^m)} C^0_{m^n}(q s^2;q,t)
C^-_\lambda(q,t;q,t) C^+_\lambda(s^2,s^2q/t;q,t)
C^0_\lambda(q^{m+1}s^2,s^2 q/t^n;q,t)}
\end{align}
When $n=1$, $\lambda=l$,
\begin{align}
\binomQ{m}{l}_{q,t,s} &=
(-1)^l q^{l(l+1)/2}
\frac{(q^{-m},q^m s^2;q)_l}{(q^l s^2,q;q)_l}
\\
\binomI{m}{l}_{q,t,s}
&=
(-1)^m q^{lm-m(m-1)/2}
\frac{(q^{-m},s^2;q)_l (1-q^{2l} s^2) (q^{m+1} s^2;q)_m}
{(q^{m+1} s^2,q;q)_l (1-q^{2m} s^2) (s^2;q)_m}
\end{align}
\end{prop}

If we state the bulk Pieri identity in terms of binomial coefficients,
we obtain the following generalized $q$-Saalsch\"utz formula.

\begin{thm}\label{thm:qSaal}
For any partitions $\kappa\subset\lambda$,
\begin{align}
\sum_{\kappa\subset \mu\subset\lambda}
q^{n(\kappa')-n(\mu')}
\frac{(-1)^{|\mu/\kappa|}
C^0_\kappa(b,c;q,t) C^-_\mu(t;q,t) C^+_\mu(a;q,t)}
{C^0_\mu(qa/b,qa/c;q,t) C^-_\kappa(t;q,t) C^+_\kappa(a;q,t)}
&P_{\mu/\kappa}([\frac{1-(qa/bc)^k}{1-t^k}];q,t)
\binomQ{\lambda}{\mu}_{q,t,\sqrt{a}}\\*
&=
(qa/bc)^{|\lambda/\kappa|}
\frac{C^0_\lambda(b,c;q,t)}
{C^0_\lambda(qa/b,qa/c;q,t)}
\binomQ{\lambda}{\kappa}_{q,t,\sqrt{a}}\notag
\end{align}
\end{thm}

\begin{proof}
For fixed $\lambda$, $\kappa$, both sides are rational functions of $b$ and
$c$; moreover, if we multiply both sides by $c^{-|\kappa|}$, the results
are well defined in the limit $(b,c)\to (0,\infty)$.  We may thus work in
the power series ring $\F(\sqrt{a})[[b,1/c]]$.  If we evaluate both sides
of the bulk Pieri identity (Theorem \ref{thm:interp:bulk_Pieri}) at a partition,
substitute
\[(s,u,v)\mapsto (t^{1-n}\sqrt{a},b/\sqrt{a},q\sqrt{a}/c),\]
and simplify, the result follows.
\end{proof}

One consequence is the following symmetry (``duality'').

\begin{cor}\label{cor:bc_duality}
For any partitions $\mu$ and $\lambda$,
\begin{align}
\binomQ{\lambda}{\mu}_{q,t,s}
&=
\binomQ{\lambda'}{\mu'}_{t,q,1/\sqrt{qt}s}\\
\binomI{\lambda}{\mu}_{q,t,s}
&=
\binomI{\lambda'}{\mu'}_{t,q,1/\sqrt{qt}s}.
\end{align}
\end{cor}

\begin{proof}
For the first
equation, we observe that if we conjugate $\kappa$, $\mu$, $\lambda$ in the
generalized $q$-Saalsch\"utz formula and substitute
\[
(q,t,a,b,c)\mapsto (t,q,1/qta,1/b,1/c),
\]
we obtain the generalized $q$-Saalsch\"utz formula again, except with
\[
\binomQ{\lambda}{\mu}_{q,t,\sqrt{a}}
\quad\text{replaced by}\quad
\binomQ{\lambda'}{\mu'}_{t,q,1/\sqrt{qta}}.
\]
Thus both binomial coefficients satisfy the same set of recurrences, with
the same initial conditions, and must therefore be the same.  This proves
the first equation; the second equation follows immediately.
\end{proof}

\begin{rem}
An alternate proof is given in Corollary \ref{cor:bc_duality2} below.
\end{rem}

The $q$-Saalsch\"utz formula can also be written in the following form,
obtained by ``reversing the order of summation''; that is, replacing
$\lambda$, $\mu$, $\kappa$ by their complements $m^n-\lambda$, $m^n-\mu$,
$m^n-\kappa$ for sufficiently large $m$ and $n$.

\begin{cor}\label{cor:qSaal:rev}
For any partitions $\kappa\subset\lambda$,
\begin{align}
\sum_{\kappa\subset\mu\subset\lambda}
q^{n(\mu')-n(\lambda')}
\frac{(-1)^{|\lambda/\mu|}
C^0_\mu(b,c;q,t) C^-_\lambda(t;q,t) C^+_\lambda(a;q,t)}
{C^0_\lambda(qa/b,qa/c;q,t) C^-_\mu(t;q,t) C^+_\mu(a;q,t)}
&P_{\lambda/\mu}([\frac{(bc/qa)^k-1}{1-t^k}];q,t)
\binomI{\mu}{\kappa}_{q,t,\sqrt{a}}\\*
&=
(bc/aq)^{|\lambda/\kappa|}
\frac{C^0_\kappa(b,c;q,t)}
{C^0_\kappa(qa/b,qa/c;q,t)}
\binomI{\lambda}{\kappa}_{q,t,\sqrt{a}}\notag
\end{align}
\end{cor}

Our definition of the second kind of binomial coefficients is justified
by the following result.

\begin{thm}\label{thm:interp:inversion}
The binomial coefficients satisfy the inversion identities
\[
\sum_{\kappa\subset\mu\subset\lambda}
\binomQ{\mu}{\kappa}_{q,t,s}
\binomI{\lambda}{\mu}_{q,t,s}
=
\delta_{\lambda\kappa}
\]
and
\[
\sum_{\kappa\subset\mu\subset\lambda}
\binomI{\mu}{\kappa}_{q,t,s}
\binomQ{\lambda}{\mu}_{q,t,s}
=
\delta_{\lambda\kappa}.
\]
\end{thm}

\begin{proof}
Fix integers $m,n\ge 0$, and define the matrix
$
\binomI{\lambda}{\mu}'_{q,t,s}
$
indexed by partitions $\lambda,\mu\subset m^n$ to be the inverse of the matrix
$
\binomQ{\mu}{\kappa}_{q,t,s}$.  The theorem is then equivalent to the equation
\[
\binomI{\lambda}{\mu}'_{q,t,s} = \binomI{\lambda}{\mu}_{q,t,s}.
\]

If we multiply both sides of the $q$-Saalsch\"utz formua by
\[
(bc/qa)^{|\rho/\kappa|}
\binomI{\rho}{\lambda}'_{q,t,\sqrt{a}}
\binomI{\kappa}{\nu}'_{q,t,\sqrt{a}}
\]
and sum over $\lambda$ and $\kappa$, we find that Corollary
\ref{cor:qSaal:rev} is satisfied by the alternate binomial coefficients
as well; as in Corollary \ref{cor:bc_duality}, this implies that the
two sets of binomial coefficients are the same.
\end{proof}

\begin{rem}
When $\lambda\ne \kappa$ both have length 1, the first sum becomes a
${}_4\phi_3$, summed to 0 by \cite[Eq. 2.3.4]{GasperG/RahmanM:1990}.  The
second sum is a ${}_2\phi_1$ in the univariate case.
\end{rem}

Together with inversion, the generalized $q$-Saalsch\"utz formula implies
the following identity, which generalizes the sum of a terminating
very-well-poised ${}_6\phi_5$.

\begin{cor}\label{cor:6vwp5}
For any partitions $\kappa\subset\lambda$,
\begin{align}
\sum_{\kappa\subset \mu\subset\lambda}
(qa&/bc)^{|\mu/\kappa|}\frac{C^0_\mu(b,c;q,t)}{C^0_\mu(qa/b,qa/c;q,t)}
\binomI{\lambda}{\mu}_{q,t,\sqrt{a}}
\binomQ{\mu}{\kappa}_{q,t,\sqrt{a}}\notag\\*
&=
(-1)^{|\lambda/\kappa|}
q^{n(\kappa')-n(\lambda')}
\frac{C^-_\lambda(t;q,t) C^+_\lambda(a;q,t)C^0_\kappa(b,c;q,t)}
{C^-_\kappa(t;q,t) C^+_\kappa(a;q,t)C^0_\lambda(aq/b,aq/c;q,t)}
P_{\lambda/\kappa}([\frac{1-(qa/bc)^k}{1-t^k}];q,t).
\end{align}
\end{cor}

\begin{proof}
In the left-hand side, expand
\[
(qa/bc)^{|\mu/\kappa|}
\frac{C^0_\mu(b,c;q,t)}{C^0_\mu(qa/b,qa/c;q,t)}
\binomQ{\mu}{\kappa}_{q,t,\sqrt{a}}
\]
by applying the generalized $q$-Saalsch\"utz formula in reverse.
We can then sum over $\mu$ using inversion, obtaining the desired sum.
\end{proof}

\begin{rem}
This generalizes equation $(2.4.2)$ of \cite{GasperG/RahmanM:1990}; the above
proof is a direct adaptation of the proof in the univariate case.
The special case $\kappa=0$, $\lambda=m^n$ was shown in \cite{vanDiejen:1997}.
\end{rem}

Iterating the above argument gives a generalization of Watson's
tranformation between a terminating very-well-poised ${}_8\phi_7$ and
a balanced terminating ${}_4\phi_3$.

\begin{thm}\label{thm:8vwp7}
For any partitions $\kappa\subset\lambda$
\begin{align}
\sum_{\kappa\subset \mu\subset\lambda}&
\left(\frac{a^2q^2}{bcde}\right)^{|\mu|-|\kappa|}
\frac{C^0_\mu(b,c,d,e;q,t)}{C^0_\mu(aq/b,aq/c,aq/d,aq/e;q,t)}
\binomI{\lambda}{\mu}_{q,t,\sqrt{a}}
\binomQ{\mu}{\kappa}_{q,t,\sqrt{a}}\notag\\*
&=
(-1)^{|\lambda|-|\kappa|}
q^{n(\kappa')-n(\lambda')}
\frac
{C^0_\kappa(b,c;q,t)C^-_\lambda(t;q,t) C^+_\lambda(a;q,t)}
{C^0_\lambda(aq/d,aq/e;q,t)C^-_\kappa(t;q,t) C^+_\kappa(a;q,t)}\notag\\*
&\phantom{{}={}}
\sum_{\kappa\subset \mu\subset\lambda}
\frac
{C^0_\mu(d,e;q,t)}
{C^0_\mu(aq/b,aq/c;q,t)}
P_{\lambda/\mu}([\frac{1-(aq/de)^k}{1-t^k}];q,t)
P_{\mu/\kappa}([\frac{(aq/de)^k-(a^2q^2/bcde)^k}{1-t^k}];q,t)
\end{align}
\end{thm}

\begin{rem}
Taking $b=aq/c$ or $d=aq/e$ recovers Corollary \ref{cor:6vwp5}.
\end{rem}

If we exchange $c$ and $d$, the left-hand side is unchanged, thus leading
to a transformation of the right-hand side, a multivariate analogue of
Sears' transformation of a balanced terminating ${}_4\phi_3$.

\begin{cor}\label{cor:Sears}
For any partitions $\kappa\subset\lambda$,
\begin{align}
\sum_{\kappa\subset \mu\subset\lambda}
&
\frac{
C^0_\lambda(aq/b,aq/c;q,t)
C^0_\mu(d,e;q,t)
}{
C^0_\mu(aq/b,aq/c;q,t)
C^0_\kappa(d,e;q,t)
}
P_{\lambda/\mu}([\frac{1-(aq/de)^k}{1-t^k}];q,t)
P_{\mu/\kappa}([\frac{(aq/de)^k-(a^2q^2/bcde)^k}{1-t^k}];q,t)\\*
=&
\sum_{\kappa\subset \mu\subset\lambda}
\frac{
C^0_\lambda(aq/b,aq/d;q,t)
C^0_\mu(c,e;q,t)
}{
C^0_\mu(aq/b,aq/d;q,t)
C^0_\kappa(c,e;q,t)
}
P_{\lambda/\mu}([\frac{1-(aq/ce)^k}{1-t^k}];q,t)
P_{\mu/\kappa}([\frac{(aq/ce)^k-(a^2q^2/bcde)^k}{1-t^k}];q,t)\notag
\end{align}
\end{cor}

This implies two more transformations, the first of which is another
generalized $q$-Saalsch\"utz formula (not used below).

\begin{cor}
For any partitions $\kappa\subset\lambda$,
\[
\sum_{\kappa\subset\mu\subset\lambda}
\frac{C^0_\mu(a;q,t)}
{C^0_\mu(c;q,t)}
P_{\lambda/\mu}([\frac{a^k-b^k}{1-t^k}];q,t)
P_{\mu/\kappa}([\frac{b^k-c^k}{1-t^k}];q,t)
=
\frac{C^0_\kappa(a;q,t)C^0_\lambda(b;q,t)}
{C^0_\kappa(b;q,t)C^0_\lambda(c;q,t)}
P_{\lambda/\kappa}([\frac{a^k-c^k}{1-t^k}];q,t)
\label{eq:weakqSaal}
\]
Similarly,
\begin{align}
\sum_{\kappa\subset \mu\subset\lambda}
&
\frac{
C^0_\lambda(aq/b,aq/c;q,t)
C^0_\mu(d,e;q,t)
}{
C^0_\mu(aq/b,aq/c;q,t)
C^0_\kappa(d,e;q,t)
}
P_{\lambda/\mu}([\frac{1-(aq/de)^k}{1-t^k}];q,t)
P_{\mu/\kappa}([\frac{(aq/de)^k-(a^2q^2/bcde)^k}{1-t^k}];q,t)\\*
=&\sum_{\kappa\subset \mu\subset\lambda}
\frac{
C^0_\lambda(aq/d,aq/e;q,t)
C^0_\mu(b,c;q,t)
}{
C^0_\mu(aq/d,aq/e;q,t)
C^0_\kappa(b,c;q,t)
}
P_{\lambda/\mu}([\frac{1-(aq/bc)^k}{1-t^k}];q,t)
P_{\mu/\kappa}([\frac{(aq/bc)^k-(a^2q^2/bcde)^k}{1-t^k}];q,t)
\notag
\end{align}
\end{cor}

\begin{proof}
The first identity follows from the special case $d=aq/e$ of Corollary
\ref{cor:Sears}; the second identity follows from two applications of that
corollary.
\end{proof}

\begin{rems}
Setting $\kappa=0$ in \eqref{eq:weakqSaal}, multiplying both sides by
$Q_\lambda(x;q,t)$ and and summing over $\lambda$ gives a multivariate
analogue of Euler's transformation (see \cite{BakerTH/ForresterPJ:1999} for
an alternate proof).
\end{rems}

\begin{rems}
These two identities are precisely the conditions required for the bulk
branching rule (or the bulk Pieri identity) to be self-consistent: the
first allows one to combine two adjoining applications into one, while the
second allows one to exchange consecutive applications.
\end{rems}

The case $\lambda=m^n$, $\kappa=0$ of Theorem \ref{thm:8vwp7} is of special
interest:

\begin{cor}
For any integers $m,n\ge 0$,
\begin{align}
\sum_{\mu\subset m^n}
&\frac{
C^0_{2\mu^2}(aq;q,t) C^0_\mu(t^n,q^{-m},b,c,d,e;q,t)t^{2n(\mu)}
(a^2q^{m+2}/t^{n-1}bcde)^{|\mu|}}
{
C^+_\mu(a,qa/t;q,t)
C^0_\mu(aq/t^n,q^{m+1}a,aq/b,aq/c,aq/d,aq/e;q,t) C^-_\mu(q,t;q,t)}\notag\\*
&\phantom{C^+_\mu(a,t a;q,t)}=
\frac{C^0_{m^n}(aq,aq/de;q,t)}{C^0_{m^n}(aq/d,aq/e;q,t)}
\sum_{\mu\subset m^n}
\frac{C^0_\mu(t^n,q^{-m},d,e,aq/bc;q,t) t^{2n(\mu)} q^{|\mu|}}
{C^0_\mu(aq/b,aq/c,t^{n-1} q^{-m}de/a;q,t) C^-_\mu(q,t;q,t)}.
\end{align}
\end{cor}

For future use, we define (horizontally) nonterminating versions of these sums:
\begin{align}
{}_8W^{(n)}_7(a;b,c,d,e,f;q,t;z)
&:=
\sum_{\ell(\mu)\le n}
\frac{
C^0_{2\mu^2}(aq;q,t) C^0_\mu(t^n,b,c,d,e,f;q,t)t^{2n(\mu)}
z^{|\mu|}}
{
C^+_\mu(a,qa/t;q,t)
C^0_\mu(aq/t^n,aq/b,aq/c,aq/d,aq/e,aq/f;q,t) C^-_\mu(q,t;q,t)}\\
{}_4\Phi^{(n)}_3\left(\genfrac{}{}{0pt}{}{a,b,c,d}{e,f,g};q,t;z\right)
&:=
\sum_{\ell(\mu)\le n}
\frac{C^0_\mu(t^n,a,b,c,d;q,t) t^{2n(\mu)} z^{|\mu|}}
{C^0_\mu(e,f,g;q,t) C^-_\mu(q,t;q,t)}.
\end{align}
Both sums converge if $|q|,|z|<1$; in contrast, convergence of a similar
vertically nonterminating sum would require $|t|>1$.  When $n=1$, we have
\begin{align}
{}_8W^{(1)}_7(a;b,c,d,e,f;q,t;z)&={}_8W_7(a;b,c,d,e,f;q;z)\\
{}_4\Phi^{(1)}_3\left(\genfrac{}{}{0pt}{}{a,b,c,d}{e,f,g};q,t;z\right)
&=
{}_4\phi_3\left(\genfrac{}{}{0pt}{}{a,b,c,d}{e,f,g};q;z\right).
\end{align}


\begin{cor}\cite{vanDiejenJF/SpiridonovVP:2000}
If $a^2 q^{m+1}=t^{n-1} bcde$, then
\[
{}_8W^{(n)}_7(a;b,c,d,e,q^{-m};q,t;q)
=
\frac{C^0_{m^n}(aq,aq/cd,aq/ce,aq/de;q,t)}{C^0_{m^n}(aq/c,aq/d,aq/e,aq/cde;q,t)}
\]
\end{cor}

\begin{rem}
This generalizes Jackson's sum of a balanced, very-well-poised, terminating
${}_8\phi_7$.  In \cite{bctheta}, we will derive a version of this indexed
by skew diagrams, as well as an associated analogue of Bailey's
transformation; in particular, we obtain Warnaar's conjectured multivariate
elliptic Bailey transform, \cite[Conjecture 6.1]{WarnaarSO:2002}.  Warnaar
also conjectured the elliptic analogue of the above sum, since proved in
\cite{RosengrenH:2001}.
\end{rem}

In the sequel, it will be useful to know how the difference equation is
expressed in terms of the binomial coefficients.

\begin{thm}\label{thm:interp:bc_deq}
For any partitions $\mu\subset\lambda$,
\begin{align}
\psi^{(d)}_{\mu/\mu}(u;q,t,s)
\binomQ{\lambda}{\mu}_{q,t,s}
&=
\sum_{\kappa\prec\lambda}
\psi^{(d)}_{\lambda/\kappa}(u;q,t,s)
\binomQ{\kappa}{\mu}_{q,t,s\sqrt{q}}\\
\psi^{(d)}_{\lambda/\lambda}(u;q,t,s)
\binomI{\lambda}{\mu}_{q,t,s\sqrt{q}}
&=
\sum_{\kappa\succ\mu}
\psi^{(d)}_{\kappa/\mu}(u;q,t,s)
\binomI{\lambda}{\kappa}_{q,t,s},
\end{align}
where
{\allowdisplaybreaks
\begin{align}
\psi^{(d)}_{\lambda/\kappa}(u;q,t,s)
&=
(-u/t)^{|\lambda/\kappa|} t^{n(\kappa)-n(\lambda)}
\frac{C^0_\lambda(qt s^2/u;q,t) C^0_\kappa(qu/t;q,t)}
     {C^0_\lambda(u/t;q,t) C^0_\kappa(qt s^2/u;q,t)}\notag\\*
&\phantom{{}={}}
\prod_{\substack{(i,j)\in \lambda\\\lambda_i=\kappa_i}}
\frac{1-q^{\lambda_i+j-1} t^{2-\lambda'_j-i} s^2}
     {1-q^{\kappa_i-j} t^{\kappa'_j-i+1}}
\prod_{\substack{(i,j)\in \lambda\\\lambda_i\ne \kappa_i}}
\frac{1-q^{\lambda_i-j+1} t^{\lambda'_j-i}}
     {1-q^{\kappa_i+j+1} t^{1-\kappa'_j-i} s^2}\notag\\*
&\phantom{{}={}}
\prod_{\substack{(i,j)\in \kappa\\\lambda_i=\kappa_i}}
\frac{1-q^{\lambda_i-j}t^{\lambda'_j-i+1}}
     {1-q^{\kappa_i+j} t^{2-\kappa'_j-i} s^2}
\prod_{\substack{(i,j)\in \kappa\\\lambda_i\ne \kappa_i}}
\frac{1-q^{\lambda_i+j} t^{1-\lambda'_j-i} s^2}
     {1-q^{\kappa_i-j+1} t^{\kappa'_j-i}}
\label{eq:interp:psid1}\\
\psi^{(d)}_{\lambda/\lambda}(u;q,t,s)
&=
\frac{
C^+_\lambda(s^2;q,t)
}{
C^+_\lambda(s^2q;q,t)}
\prod_{1\le i\le \ell(\lambda)}
\frac{1-q^{\lambda_i} t^{-i} u}
{1-t^{-i} u}.
\end{align}
}
\end{thm}

\begin{proof}
For the first claim, multiply $s$ and $u$ by $t^{-n}$ in the difference
equation, divide both sides by
\[
\bar{P}^{*(n)}_\mu(\mu;q,t,t^{1-n} s\sqrt{q})
\prod_{1\le i\le n} (1-t^{-i} u),
\]
and evaluate at $\lambda$.  We thus obtain an equation of the desired
form; it remains to simplify the coefficients on the right hand side,
namely
\begin{align}
&\prod_{i\in R'}
\frac{
(1-q^{\lambda_i} t^{2-n-i} s^2)(1-u q^{\lambda_i} t^{-i})
}{
(1-q^{2\lambda_i} t^{2-2i} s^2)(1-u t^{-i})
}
\prod_{i\in R}
\frac{
(1-q^{-\lambda_i} t^{i-n})(1-u q^{-\lambda_i} t^{i-2}/s^2)
}{
(1-q^{-2\lambda_i} t^{2i-2} s^{-2})(1-u t^{-i})
}\notag\\\times
&\prod_{i<j\in R'}
\frac{
(1-q^{\lambda_i+\lambda_j} t^{3-i-j} s^2)
}
{
(1-q^{\lambda_i+\lambda_j} t^{2-i-j} s^2)
}
\prod_{i\in R',j\in R}
\frac{
(1-q^{\lambda_i-\lambda_j} t^{j+1-i})
}
{
(1-q^{\lambda_i-\lambda_j} t^{j-i})
}
\prod_{i<j\in R}
\frac{
(1-q^{-\lambda_i-\lambda_j} t^{i+j-1} s^{-2})
}
{
(1-q^{-\lambda_i-\lambda_j} t^{i+j-2} s^{-2})
},
\end{align}
where
\begin{align}
R&=\{i:i\in \{1,2,\dots n\}|\lambda_i=\kappa_i+1\}\\
R'&=\{i:i\in \{1,2,\dots n\}|\lambda_i=\kappa_i\}.
\end{align}

We have, for instance,
{\allowdisplaybreaks
\begin{align}
\frac{\prod_{i<j\in R} (1-q^{-\lambda_i-\lambda_j} t^{i+j-1} s^{-2})}
{\prod_{i\le j\in R} (1-q^{-\lambda_i-\lambda_j} t^{i+j-2} s^{-2})}
&\propto
\frac{\prod_{i<j\in R} (1-q^{\lambda_i+\lambda_j} t^{1-i-j} s^2)}
{\prod_{i\le j\in R} (1-q^{\lambda_i+\lambda_j} t^{2-i-j} s^2)}\\
&=
\phantom{{}\times{}}
\prod_{i\in R}
\prod_{1\le k\le \kappa_i}
\prod_{\substack{j\in R\\\lambda_j=k}}
\frac{1-q^{\lambda_i+k} t^{1-i-j} s^2}
{1-q^{\lambda_i+k} t^{2-i-j} s^2}\notag\\*
&\phantom{{}={}}\times
\prod_{1\le k}
\frac{
\prod_{\kappa'_k<i<j\le \lambda'_k}
(1-q^{2k} t^{1-i-j} s^2)}
{
\prod_{\kappa'_k<i\le j\le \lambda'_k}
(1-q^{2k} t^{2-i-j} s^2)}
\\
&=
\prod_{i\in R}
\prod_{1\le k\le \kappa_i}
\frac{1-q^{\lambda_i+k} t^{1-\lambda'_k-i} s^2}
     {1-q^{\lambda_i+k} t^{1-\kappa'_k-i} s^2}
\prod_{1\le k}
\prod_{\substack{j\in R\\\lambda_j=k}}
\frac{1}{1-q^{\lambda_j+k} t^{1-\kappa'_k-j} s^2}\\
&=
\prod_{i\in R}
\frac{\prod_{1\le j\le \kappa_i}(1-q^{\lambda_i+j} t^{1-\lambda'_j-i} s^2)}
     {\prod_{1\le j\le \lambda_i} (1-q^{\kappa_i+1+j} t^{1-\kappa'_j-i} s^2)},
\end{align}
}
where the constant of proportionality is
\[
\frac{\prod_{i<j\in R} (-q^{-\lambda_i-\lambda_j} t^{i+j-1} s^{-2})}
{\prod_{i\le j\in R} (-q^{-\lambda_i-\lambda_j} t^{i+j-2} s^{-2})}
=
t^{|R|(|R|-1)/2} \prod_{i\in R} (-q^{2\lambda_i} t^{2-2i} s^2).
\]
We thus obtain two of the factors of \eqref{eq:interp:psid1}.  Here we used
the fact that $i\in R$ and $\lambda_i=k$ if and only if $\kappa'_k<i\le
\lambda'_k$; similarly, $i\in R'$ and $\lambda_i=k$ if and only if
$\lambda_{k+1}<i\le \kappa'_k$.  The remaining simplifications are analogous.

The second equation follows in a similar way; here we substitute
$(s,u)\mapsto (q^{-m}/ts,t q^{-m}/u)$.  Alternatively, it follows
immediately from the first via inversion.
\end{proof}

Dualizing the difference equations (via Corollary \ref{cor:bc_duality})
gives ``integral equations''.

\begin{cor}
For any partition $\mu\subset\lambda$,
\begin{align}
\psi^{(i)}_{\mu/\mu}(u;q,t,s)
\binomQ{\lambda}{\mu}_{q,t,s\sqrt{t}}
&=
\sum_{\kappa'\prec\lambda'}
\psi^{(i)}_{\lambda/\kappa}(u;q,t,s)
\binomQ{\kappa}{\mu}_{q,t,s}\\
\psi^{(i)}_{\lambda/\lambda}(u;q,t,s)
\binomI{\lambda}{\mu}_{q,t,s}
&=
\sum_{\nu'\succ\mu'}
\psi^{(i)}_{\nu/\mu}(u;q,t,s)
\binomI{\lambda}{\nu}_{q,t,s\sqrt{t}},
\end{align}
where
\begin{align}
\psi^{(i)}_{\lambda/\kappa}(u;q,t,s)
&=
(u/t)^{|\lambda|-|\kappa|}
t^{n(\kappa)-n(\lambda)}
\frac{C^0_\lambda(s^2qt/u;q,t) C^0_\kappa(u/t;q,t)}
     {C^0_\lambda(u;q,t) C^0_\kappa(s^2qt/u;q,t)}
\notag\\
&\phantom{{}={}}
\prod_{\substack{(i,j)\in \lambda\\\lambda'_j=\kappa'_j}}
\frac{1-q^{\lambda_i+j-1} t^{-\lambda'_j-i+3} s^2}
     {1-q^{\kappa_i-j+1} t^{\kappa'_j-i}}
\prod_{\substack{(i,j)\in \lambda\\\lambda'_j\ne \kappa'_j}}
\frac{1-q^{\lambda_i-j} t^{\lambda'_j-i+1}}
     {1-q^{\kappa_i+j} t^{-\kappa'_j-i+1} s^2}\notag\\
&\phantom{{}={}}
\prod_{\substack{(i,j)\in \kappa\\\lambda'_j=\kappa'_j}}
\frac{1-q^{\lambda_i-j+1}t^{\lambda'_j-i}}
     {1-q^{\kappa_i+j-1}t^{2-\kappa'_j-i} s^2}
\prod_{\substack{(i,j)\in \kappa\\\lambda'_j\ne \kappa'_j}}
\frac{1-q^{\lambda_i+j} t^{2-\lambda'_j-i} s^2}
     {1-q^{\kappa_i-j} t^{\kappa'_j-i+1}}\\
\psi^{(i)}_{\lambda/\lambda}(u;q,t,s)
&=
\frac{C^0_\lambda(u/t;q,t)}
     {C^0_\lambda(u;q,t)}
\frac{C^+_\lambda(s^2 t)}
     {C^+_\lambda(s^2)}
\end{align}
\end{cor}

\begin{rem}
These identities can be analytically continued to give a one-parameter
family of integral equations for interpolation polynomials, having the
integral representation of \cite{OkounkovA:1998} as the case $u=t^n$,
$\ell(\mu)<n$.  In fact, we discovered these integral equations first (in
order to prove Theorem \ref{thm:koorn:connt}), then deduced the likely
existence of difference equations via duality.  As the integral operators
are rather complicated, and unnecessary for our purposes, we omit the
details, and note simply that they correspond to the operators defined in
\cite{xforms} via contour integrals.
\end{rem}

\section{Koornwinder polynomials}\label{sec:koorn}

\begin{defn}
The {\it Koornwinder polynomials} are the unique
family of $BC_n$-symmetric polynomials
\[
K^{(n)}_\lambda(;q,t;t_0,t_1,t_2,t_3)
\]
such that
\begin{itemize}
\item (Triangularity) $K^{(n)}_\lambda(;q,t;t_0,t_1,t_2,t_3)=m_\lambda+\text{dominated terms}$.
\item (Evaluation symmetry)
For any pair of partitions $\mu<\lambda$,
\[
\frac{
K^{(n)}_\lambda(q^{\mu_i} t^{n-i} t_0;q,t;t_0,t_1,t_2,t_3)
}{
k^0_\lambda(q,t,t^n;t_0{:}t_1,t_2,t_3)
}
=
\frac{
K^{(n)}_\mu(q^{\lambda_i} t^{n-i} \hat{t}_0;q,t;\hat{t}_0,\hat{t}_1,\hat{t}_2,\hat{t}_3)
}{
k^0_\mu(q,t,t^n;\hat{t}_0{:}\hat{t}_1,\hat{t}_2,\hat{t}_3)
},
\]
where
\begin{align}
\hat{t}_0=\sqrt{t_0t_1t_2t_3/q}&\text{;\qquad}
\hat{t}_i=t_0t_i/\hat{t}_0,\text{ $i\in \{1,2,3\}$}\\
k^0_\lambda(q,t,T;t_0{:}t_1,t_2,t_3)
&=
(t_0 T/t)^{-|\lambda|} t^{n(\lambda)}
\frac{C^0_\lambda(T,T t_0 t_1/t,T t_0 t_2/t,T t_0 t_3/t;q,t)}
{C^-_\lambda(t;q,t) C^+_\lambda(T^2 \hat{t}_0^2/t^2;q,t)}.
\end{align}
\end{itemize}
\end{defn}

This differs from the definition in the literature (in which evaluation
symmetry is replaced by orthogonality with respect to the Koornwinder inner
product); that our definition is equivalent to the usual definition will
be shown below.

\begin{thm}
The Koornwinder polynomials are well-defined, and are given
by the expansion
\[
K^{(n)}_\lambda(;q,t;t_0,t_1,t_2,t_3)
=
\sum_{\mu\subset\lambda}
\binomQ{\lambda}{\mu}_{q,t,t^{n-1} \hat{t}_0}
\frac{k^0_\lambda(q,t,t^n;t_0{:}t_1,t_2,t_3)}
     {k^0_\mu(q,t,t^n;t_0{:}t_1,t_2,t_3)}
\bar{P}^{*(n)}_\mu(;q,t;t_0).\label{eq:koorn:binomial}
\]
\end{thm}

\begin{proof}
A monic triangular $BC_n$-symmetric polynomial with leading term
$m_\lambda$ is uniquely determined by its values at $q^{\mu_i} t^{n-i} t_0$
for $\mu<\lambda$.  Indeed, we can write it as
\[
m_\lambda + \sum_{\mu<\lambda} c_{\lambda\mu} \bar{P}^{*(n)}_\mu,
\]
where the coefficients $c_{\lambda\mu}$ are determined by a triangular
system of linear equations with nonzero diagonal.

We find that, if the expansion holds for all $\mu<\lambda$, then
\begin{align}
\frac{
K^{(n)}_\lambda(q^{\mu_i}t^{n-i} t_0;q,t;t_0,t_1,t_2,t_3)}
{k^0_\lambda(q,t,t^n;t_0{:}t_1,t_2,t_3)}
&=
\frac{
K^{(n)}_\mu(q^{\lambda_i} t^{n-i} \hat{t}_0;q,t
                               ;\hat{t}_0,\hat{t}_1,\hat{t}_2,\hat{t}_3)}
{k^0_\mu(q,t,t^n;\hat{t}_0{:}\hat{t}_1,\hat{t}_2,\hat{t}_3)}\\
&=
\sum_{\nu\subset\mu}
\binomQ{\mu}{\nu}_{q,t,t^{n-1} t_0}
\binomQ{\lambda}{\nu}_{q,t,t^{n-1} \hat{t}_0}
\frac{\bar{P}^{*(n)}_\nu(\nu;q,t,\hat{t}_0)}
     {k^0_\nu(q,t,t^n;\hat{t}_0{:}\hat{t}_1,\hat{t}_2,\hat{t}_3)}\\
&=
\sum_{\nu\subset\lambda}
\binomQ{\mu}{\nu}_{q,t,t^{n-1} t_0}
\binomQ{\lambda}{\nu}_{q,t,t^{n-1} \hat{t}_0}
\frac{\bar{P}^{*(n)}_\nu(\nu;q,t,t_0)}
     {k^0_\nu(q,t,t^n;t_0{:}t_1,t_2,t_3)}\\
&=
\sum_{\nu\subset\lambda}
\binomQ{\lambda}{\nu}_{q,t,t^{n-1} t_0}
\frac{\bar{P}^{*(n)}_\nu(\mu;q,t,t_0)}
     {k^0_\nu(q,t,t^n;t_0{:}t_1,t_2,t_3)},
\end{align}
where the second-to-last step follows from the fact that
\[
k^0_\nu(q,t,t^n;t_0{:}t_1,t_2,t_3)
\bar{P}^{*(n)}_\nu(\nu;q,t,\hat{t}_0)
\]
depends on $t_0,\dots t_3$ only through their pairwise products, and
is thus preserved by the ``hat'' involution.  In particular, both sides of
\eqref{eq:koorn:binomial} are monic triangular with leading term $\lambda$,
and agree at $q^{\mu_i} t^{n-i} t_0$ whenever $\mu<\lambda$.  The identity
follows.
\end{proof}

\begin{rems}
This, of course, is essentially Okounkov's binomial formula
\cite{OkounkovA:1998}; the difference is that the principal
specializations in Okounkov's formula have been replaced with the
appropriate product.  Note that in the univariate case, this is precisely
the expansion of an Askey-Wilson polynomial as a ${}_4\phi_3$
\cite{AskeyR/WilsonJ:1985}.
\end{rems}

\begin{rems}
We will tend to avoid the $\hat{t}_i$ notation in the sequel, as it is
really only suited to contexts in which the parameters are fixed;
in other contexts, it can lead to serious ambiguities.
\end{rems}

\begin{cor}
The evaluation symmetry property holds without restriction on $\mu$ and
$\lambda$.
\end{cor}

\begin{cor}
For any partition $\lambda$, the only possible factors of the denominators
of the coefficients of the polynomial
\[
C^+_\lambda(t^{2n-2} t_0t_1t_2t_3/q;q,t)
\bar{K}_\lambda(;q,t;t_0,t_1,t_2,t_3)
\]
are binomials of the form $1-q^i t^j$ for $i,j\ge 0$.
\end{cor}

\begin{proof}
Taking into account the denominator factors introduced by the binomial
coefficient and the interpolation polynomial in the binomial formula, we
conclude that the only possible other denominator factors are $q$, $t$,
$t_0$.  Now, the Koornwinder inner product is well-defined whenever any of
$q=0$, $t=0$, or $t_0=0$; thus Theorem \ref{thm:koorn:ortho1} below implies
that $q$ and $t$ and $t_0$ are not denominator factors of $K^{(n)}$.
\end{proof}

The symmetries of interpolation polynomials induce corresponding
symmetries of Koornwinder polynomials.

\begin{cor}\label{cor:koorn:trivsim}
For any partition $\lambda$,
\begin{align}
K^{(n)}_\lambda(;q,t;t_0,t_1,t_2,t_3)
&=
K^{(n)}_\lambda(;1/q,1/t;1/t_0,1/t_1,1/t_2,1/t_3)\\
K^{(n)}_\lambda(;q,t;t_0,t_1,t_2,t_3)
&=
(-1)^{|\lambda|} K^{(n)}_\lambda(-x;q,t;-t_0,-t_1,-t_2,-t_3)\\
\end{align}
\end{cor}

A further symmetry follows from the $q$-Saalsch\"utz formula.

\begin{thm}
For any partition $\lambda$,
\[
K^{(n)}_\lambda(;q,t;t_0,t_1,t_2,t_3)
=
K^{(n)}_\lambda(;q,t;t_1,t_0,t_2,t_3).
\]
Thus $K^{(n)}_\lambda$ is invariant under permutations of $t_0$, $t_1$,
$t_2$, $t_3$.
\end{thm}

\begin{proof}
For all $\kappa$, the coefficient of $P^{(n)}_\kappa(;q,t;t_1)$
in both sides is the same.
\end{proof}

\begin{rem}
This is another multivariate analogue of Sears' ${}_4\phi_3$
transformation.
\end{rem}

We will refer to this fact as ``parameter symmetry''.

\bigskip
Before delving further into the properties of the Koornwinder polynomials
as we have defined them, we first must justify the name.  We could of
course simply refer to the proofs in the literature
\cite{vanDiejen:1996,SahiS:1999,StokmanJV:2000b} that the usual Koornwinder
polynomials satisfy evaluation symmetry; instead, we will give a direct
proof (in particular, avoiding any use of double affine Hecke algebra
machinery).  We will, in fact, give two proofs.  The first uses difference
operators to show orthogonality with respect to the Koornwinder weight
function, while the second uses our generalized hypergeometric
transformations to show orthogonality with respect to the $q$-Racah weight
function.

First, recall the difference operators $D^{(n)}(u_1,u_2;q,t)$ of Definition
\ref{defn:diffop}.  These act on our polynomials $K^{(n)}$ as follows.

\begin{lem}\label{lem:koorn:ortho1}
For any integer $n$ and partition $\lambda$ with $\ell(\lambda)\le n$,
\[
D^{(n)}(t_0,t_1;q,t)K^{(n)}_\lambda(;q,t;t_0\sqrt{q},t_1\sqrt{q},t_2,t_3)
=
E^{(n)}_\lambda(t_0t_1;q,t)
K^{(n)}_\lambda(;q,t;t_0,t_1,t_2\sqrt{q},t_3\sqrt{q}),
\]
where
\[
E^{(n)}_\lambda(u;q,t) = q^{-|\lambda|/2} \prod_{1\le i\le n}
(1-q^{\lambda_i} t^{n-i} u).
\]
\end{lem}

\begin{proof}
By the binomial formula, we have:
\begin{align}
D^{(n)}(t_0,t_1;q,t)&K^{(n)}_\lambda(;q,t;t_0\sqrt{q},t_1\sqrt{q},t_2,t_3)\notag\\* 
&=
\sum_{\mu\subset\lambda}
\binomQ{\lambda}{\mu}_{q,t,t^{n-1}\sqrt{t_0t_1t_2t_3}}
\frac{k^0_\lambda(q,t,t^n;t_0\sqrt{q}{:}t_1\sqrt{q},t_2,t_3)}
     {k^0_\mu(q,t,t^n;t_0\sqrt{q}{:}t_1\sqrt{q},t_2,t_3)}
D^{(n)}(t_0,t_1;q,t)
\bar{P}^*_\mu(;q,t;t_0\sqrt{q})\\
&=
\sum_{\mu\subset\lambda}
\binomQ{\lambda}{\mu}_{q,t,t^{n-1}\sqrt{t_0t_1t_2t_3}}
\frac{k^0_\lambda(q,t,t^n;t_0\sqrt{q}{:}t_1\sqrt{q},t_2,t_3)}
     {k^0_\mu(q,t,t^n;t_0\sqrt{q}{:}t_1\sqrt{q},t_2,t_3)}
E^{(n)}_\mu(t_0t_1;q,t)
\bar{P}^*_\mu(;q,t;t_0).
\end{align}
Since
\[
k^0_\mu(q,t,t^n;t_0\sqrt{q}{:}t_1\sqrt{q},t_2,t_3)
=
\frac{E^{(n)}_\mu(t_0t_1;q,t)}{\prod_{1\le i\le n} (1-t^{n-i} t_0 t_1)}
k^0_\mu(q,t,t^n;t_0{:}t_1,t_2\sqrt{q},t_3\sqrt{q}),
\]
the result follows.
\end{proof}

\begin{rem}
In particular, we find that the polynomials
$K^{(n)}_\mu(;q,t;t_0,t_1,t_0\sqrt{q},t_1\sqrt{q})$ are eigenfunctions of
the operator $D^{(n)}(t_0,t_1;q,t)$, and thus are $BC_n/C_n$-Macdonald
polynomials.  More generally, the polynomials
$K^{(n)}_\mu(;q,t;t_0,t_1,t_2,t_3)$ are eigenfunctions of
$D^{(n)}(t_0,t_1;q,t)D^{(n)}(q^{-1/2} t_2,q^{-1/2} t_3;q,t)$.
\end{rem}

Let $w^{(n)}_K(;q,t;t_0,t_1,t_2,t_3)$ denote the Koornwinder weight function
\cite{KoornwinderTH:1992}
\[
w^{(n)}_K(x_1,x_2,\dots x_n;q,t;t_0,t_1,t_2,t_3)
=
\prod_{1\le i\le n} \frac{(x_i^{\pm 2};q)}
{(t_0 x_i^{\pm 1},t_1 x_i^{\pm 1},t_2 x_i^{\pm 1},t_3 x_i^{\pm 1};q)}
\prod_{1\le i<j\le n} \frac{(x_i^{\pm 1} x_j^{\pm 1};q)}{(t x_i^{\pm
1} x_j^{\pm 1};q)}.
\]
Also, for a multivariate function $f$ analytic in a neighborhood of the
unit torus (the locus in which all variables have magnitude 1), $\int f d\T$
denotes the integral of $f$ with respect to the uniform density on the
torus; i.e.,
\begin{align}
\int f d\T
&=
\int_{[-\pi,\pi]^n} f(e^{i\theta_1},e^{i\theta_2},\dots e^{i\theta_n})
\prod_j \frac{d\theta_j}{2\pi}\\
&=
\int_{|x_j|=1} f(x_1,x_2,\dots x_n)
\prod_j \frac{dx_j}{2\pi i x_j},
\end{align}
or equivalently the constant coefficient of the Laurent expansion of $f$.

A straightforward adaptation of the standard adjointness argument for
$BC_n/C_n$-Macdonald polynomials \cite{MacdonaldIG:2000} proves the following.

\begin{lem}
Let $q$, $t$, $t_0$, $t_1$, $t_2$, $t_3$ be arbitrary complex numbers of
magnitude $<1$.  Then for any integer $n$ and $BC_n$-symmetric polynomials
$f$, $g$,
\begin{align}
\int
(D^{(n)}(t_0,t_1;q,t)g)f\,
w^{(n)}(;q,t;t_0,t_1,t_2\sqrt{q},t_3\sqrt{q})
d\T
=
\int
(D^{(n)}(t_2,t_3;q,t)f)g\,
w^{(n)}(;q,t;t_2,t_3,t_0\sqrt{q},t_1\sqrt{q})
d\T.
\notag
\end{align}
\end{lem}

\begin{proof}
Factor the two weight functions as:
\begin{align}
w^{(n)}(x_1,x_2,\dots x_n;q,t;t_0,t_1,t_2\sqrt{q},t_3\sqrt{q})
&=
\Delta^{(n)}_1(x_1,x_2,\dots x_n) \Delta^{(n)}_1(x_1^{-1},x_2^{-1},\dots x_n^{-1})\\
w^{(n)}(x_1,x_2,\dots x_n;q,t;t_2,t_3,t_0\sqrt{q},t_1\sqrt{q})
&=
\Delta^{(n)}_2(x_1,x_2,\dots x_n) \Delta^{(n)}_2(x_1^{-1},x_2^{-1},\dots x_n^{-1}),
\end{align}
where
\begin{align}
\Delta^{(n)}_1(x_1,x_2,\dots x_n)
&=
\prod_{1\le i\le n} \frac{(x_i^2;q)}
{(t_0 x_i,t_1 x_i,t_2\sqrt{q} x_i,t_3\sqrt{q} x_i;q)}
\prod_{1\le i<j\le n} \frac{(x_i x_j^{\pm 1};q)}{(t x_i x_j^{\pm 1};q)}\\
\Delta^{(n)}_2(x_1,x_2,\dots x_n)
&=
\prod_{1\le i\le n} \frac{(x_i^2;q)}
{(t_0\sqrt{q} x_i,t_1\sqrt{q} x_i,t_2 x_i,t_3 x_i;q)}
\prod_{1\le i<j\le n} \frac{(x_i x_j^{\pm 1};q)}{(t x_i x_j^{\pm 1};q)}.
\end{align}
The difference operators can then be expressed in the form
\begin{align}
(D^{(n)}(t_0,t_1;q,t)g)(x_1,x_2,\dots x_n)
&=
\sum_{\sigma\in \{\pm 1\}^n}
\prod_{1\le i\le n} R_{x_i}(\sigma_i)
\frac{
(\Delta^{(n)}_2 g)(\sqrt{q} x_1,\sqrt{q} x_2,\dots \sqrt{q} x_n)
}
{\Delta^{(n)}_1(x_1,x_2,\dots x_n)}
\\
(D^{(n)}(t_2,t_3;q,t)f)(x_1,x_2,\dots x_n)
&=
\sum_{\sigma\in \{\pm 1\}^n}
\prod_{1\le i\le n} R_{x_i}(\sigma_i)
\frac{
(\Delta^{(n)}_1 f)(\sqrt{q} x_1,\sqrt{q} x_2,\dots \sqrt{q} x_n)
}
{\Delta^{(n)}_2(x_1,x_2,\dots x_n)},
\end{align}
where as in the proof of Lemma \ref{lem:Dpoly}, $R_{x_i}(\pm 1)$ are
homomorphisms defined by $R_{x_i}(\pm 1)x_j=x_j$, $R_{x_i}(\pm 1) x_i =
x_i^{\pm 1}$.  Since $f$, $w$, and $d\T$ are preserved by these operators,
we find:
\begin{align}
\int
(D^{(n)}(t_0,t_1;q,t)g)f\,
w^{(n)}(;q,t&{};t_0,t_1,t_2\sqrt{q},t_3\sqrt{q})
d\T\notag\\*
&=
2^n
\int
(\Delta^{(n)}_1 f)(x_1^{-1},x_2^{-1},\dots x_n^{-1})
(\Delta^{(n)}_2 g)(\sqrt{q} x_1,\sqrt{q} x_2,\dots \sqrt{q} x_n)
d\T\\
&=
2^n
\int
(\Delta^{(n)}_1 f)(\sqrt{q} x_1,\sqrt{q} x_2,\dots \sqrt{q} x_n)
(\Delta^{(n)}_2 g)(x_1^{-1},x_2^{-1},\dots x_n^{-1})
d\T\\
&=
\int
(D^{(n)}(t_2,t_3;q,t)f)g\,
w^{(n)}(;q,t;t_2,t_3,t_0\sqrt{q},t_1\sqrt{q})
d\T.
\end{align}
\end{proof}

\begin{thm}\label{thm:koorn:ortho1}
Let $q$, $t$, $t_0$, $t_1$, $t_2$, $t_3$ be complex numbers of magnitude at
most 1.  Then the polynomials $K^{(n)}_\mu(;q,t;t_0,t_1,t_2,t_3)$ are
orthogonal with respect to the density $w^{(n)}_K(;q,t;t_0,t_1,t_2,t_3)$ on
the unit torus.
\end{thm}

\begin{proof}
On the one hand, the polynomials $K^{(n)}_\mu(;q,t;t_0,t_1,t_2,t_3)$
are eigenfunctions of the difference operator
\[
D^{(n)}(t_0,t_1;q,t)D^{(n)}(q^{-1/2} t_2,q^{-1/2} t_3;q,t)
,
\]
with generically distinct eigenvalues.  On the other hand, this difference
operator is self-adjoint with respect to the Koornwinder weight.
Orthogonality follows immediately.
\end{proof}

\begin{rem}
In particular, we have shown that our Koornwinder polynomials agree with
the usual Koornwinder polynomials, and thus the latter satisfy evaluation
symmetry.
\end{rem}

\bigskip
We now turn to the $q$-Racah case.  Following
\cite{vanDiejenJF/StokmanJV:1998}, suppose $t_0t_1=t^{1-n} q^{-m}$, and
define a function $\Delta^{qR}$ on partitions $\mu\subset m^n$ by
\[
\Delta^{qR}(\mu)
=
q^{-2n(\mu')}
t^{2n(\mu)}
(t^{2n-2} q t_0^2)^{-|\mu|}
\frac{C^0_\mu(t^{n-1} t_0 t_2,t^{n-1} t_0 t_3;q,t)C^0_{m^n-\mu}(t^{n-1} t_1 t_2,t^{n-1} t_1 t_3;q,t)
\langle P_\mu,P_\mu\rangle''_n
}
{
\bar{P}^{*(n)}_\mu(\mu;q,t,t_0)
\bar{P}^{*(n)}_{m^n-\mu}(m^n-\mu;q,t,t_1)
}.
\]
Aside from an overall constant, this is the weight function for the
multivariate $q$-Racah polynomials of \cite{vanDiejenJF/StokmanJV:1998}.
If we define a linear functional on $BC_n$-symmetric functions by
\[
\langle f\rangle_{qR}
=
\sum_{\mu\subset m^n} f(q^{\mu_i} t^{n-i} t_0) \Delta^{qR}(\mu),
\]
then
\begin{align}
\langle \bar{P}^{*(n)}_\kappa(;&q,t,t_0)
\bar{P}^{*(n)}_{m^n-\lambda}(;q,t,t_1)\rangle_{qR}
=
\sum_{\mu\subset m^n} \Delta^{qR}(\mu)
\bar{P}^*_\kappa(\mu;q,t,t_0)
\bar{P}^*_{m^n-\lambda}(m^n-\mu;q,t,t_1)\\
&=
\frac{\bar{P}^{*(n)}_\kappa(\kappa;q,t,t_0)\langle P_\lambda,P_\lambda\rangle''_n}
{\bar{P}^{*(n)}_\lambda(\lambda;q,t,t_0)}
\sum_{\mu\subset m^n}
\frac{(q/t_2t_3)^{|\mu|}
C^0_\mu(t^{n-1} t_0 t_2,t^{n-1} t_0 t_3;q,t)}
{C^0_\mu(t^{n-1} q t_0/t_2,t^{n-1} q t_0/t_3;q,t)}
\binomI{\lambda}{\mu}_{q,t,t_0 t^{n-1}}
\binomQ{\mu}{\kappa}_{q,t,t_0 t^{n-1}}.
\end{align}
This sums via Corollary \ref{cor:6vwp5} to give
\[
\langle \bar{P}^{*(n)}_\kappa(;q,t,t_0)
\bar{P}^{*(n)}_{m^n-\lambda}(;q,t,t_1)\rangle_{qR}
\propto
\frac{
k^0_\kappa(q,t,t^n;t_0,t_1,t_2,t_3)
P^{*(n)}_\kappa(\kappa;q,t;\hat{t}_0)}
{\hat{t}_0^{|\kappa|}
C^0_\kappa(t^n;q,t)C^0_\kappa(q^m;1/q,1/t)}
P_{\lambda/\kappa}([\frac{1-(q/t_2t_3)^k}{1-t^k}];q,t),
\]
where the constant of proportionality is independent of $\kappa$,
and $\hat{t}_0 = \sqrt{t_0t_1t_2t_3/q}$ as usual.
Now, from the binomial formula, we find
\begin{align}
\langle K^{(n)}_{m^n-\mu}(;q,t;t_0,t_1,t_2,t_3)&
\bar{P}^{*(n)}_{m^n-\lambda}(;q,t,t_1)\rangle_{qR}\notag\\*
&\propto
\sum_{\kappa\subset\lambda}
P^{*(n)}_\kappa(m^n-\mu;q,t,\hat{t}_0)
\frac{C^0_\lambda(t^n;q,t)C^0_\lambda(q^m;1/q,1/t)}
{C^0_\kappa(t^n;q,t)C^0_\kappa(q^m;1/q,1/t)}
P_{\lambda/\kappa}([\frac{\hat{t}_0^k-\hat{t}_1^k}{1-t^k}];q,t)\\
&\propto
P^{*(n)}_\lambda(\mu;q,t,\hat{t}_1),\label{eq:koorn:qRacah:tmp1}
\end{align}
where $\hat{t}_1=t^{1-n}q^{-m}/\hat{t}_0$, and \ref{eq:koorn:qRacah:tmp1}
follows from the connection coefficient identity for interpolation
polynomials.  In particular,
\[
\langle K^{(n)}_{m^n-\mu}(;q,t;t_0,t_1,t_2,t_3)
\bar{P}^{*(n)}_{m^n-\lambda}(;q,t,t_1)\rangle_{qR}
=0
\]
unless $\lambda\subset\mu$; it follows that the polynomials 
$K^{(n)}_\mu(;q,t;t_0,t_1,t_2,t_3)$ are orthogonal with respect to the
given inner product.  If we keep track of the constants when $\lambda=\mu$,
we obtain the following theorem.

\begin{thm}\cite{vanDiejenJF/StokmanJV:1998}
If $t_0t_1=t^{1-n} q^{-m}$ and $\lambda\subset m^n$, then
\[
\frac{
\langle K^{(n)}_{\lambda}(;q,t;t_0,t_1,t_2,t_3)
K^{(n)}_{\lambda}(;q,t;t_1,t_0,t_2,t_3)\rangle_{qR}}
{\langle 1\rangle_{qR}}
=
\delta_{\lambda\mu} N_\lambda(;q,t,t^n;t_0,t_1,t_2,t_3),
\]
where
\[
N_\lambda(;q,t,T;t_0,t_1,t_2,t_3)
=\frac{C^-_\lambda(q;q,t)C^+_\lambda(T^2 t_0t_1t_2t_3/t^3;q,t)
C^0_\lambda(T,Tt_0t_1t_2t_3/t^2;q,t)
\prod_{0\le i<j\le 3} C^0_\lambda(T t_i t_j/t;q,t)
}
{C^-_\lambda(t;q,t)C^+_\lambda(T^2 t_0t_1t_2t_3/qt^2;q,t)
C^0_{2\lambda^2}(T^2 t_0t_1t_2t_3/t^2;q,t)}
\]
\end{thm}

\begin{rems}
The proof of \cite{vanDiejenJF/StokmanJV:1998} was based on the usual
definition of Koornwinder polynomials, and involved showing that
Koornwinder's difference operator is self-adjoint with respect to the
$q$-Racah inner product.  One can presumably construct a similar proof
based on our difference operator.
\end{rems}

\begin{rems}
Except for the formula for the norm, we could have derived this from
orthogonality with respect to the Koornwinder weight, by reference to the
results of \cite{StokmanJV:2000a}.
\end{rems}

\begin{rems}
The normalization $\langle 1\rangle_{qR}$ can be computed using the
rectangle case of the generalized ${}_6\phi_5$ sum: take $\kappa=0$,
$\lambda=m^n$ above.  (This normalization was first computed in
\cite{vanDiejen:1997} by proving this case of the ${}_6\phi_5$ sum.)
\end{rems}

\begin{rems}
Similarly, we can express the inner product
\[
\langle
\bar{P}^{*(n)}_\kappa(;q,t,t_0)
\bar{P}^{*(n)}_{m^n-\lambda}(;q,t,t_1)
\prod_{1\le i\le n}
\frac{(u x_i,u/x_i;q)}{(v x_i,v/x_i;q)}
\rangle_{qR}
\]
as a generalized very-well-poised ${}_8\phi_7$; that is, in terms of the
special case
\[
(a;b,c,d,e)\mapsto (t^{2n-2} t_0^2;t^{n-1} t_0 t_2,t^{n-1} t_0 t_3,t^{n-1}
t_0 u,q t_0 t^{n-1}/v)
\]
of Theorem \ref{thm:8vwp7}.  Is there a similar interpretation of this
sum without the $q$-Racah constraint $t_0t_1=t^{1-n} q^{-m}$?
\end{rems}

If we define the {\it virtual Koornwinder integral}
\[
I^{(n)}_K(f;q,t;t_0,t_1,t_2,t_3)
:=
[K^{(n)}_0(;q,t;t_0,t_1,t_2,t_3)] f,
\]
for any $BC_n$-symmetric function $f$, the two orthogonality results imply
the following expressions.

\begin{cor}
For $q$, $t$, $t_0\dots t_3$ of magnitude $<1$, and all $BC_n$-symmetric
polynomials $f$,
\[
I^{(n)}_K(f;q,t;t_0,t_1,t_2,t_3)
=
Z^{-1}
\int f(x_1,x_2,\dots x_n)
w^{(n)}_K(x;q,t;t_0,t_1,t_2,t_3) d\T,
\]
where
\[
Z=\int w^{(n)}_K(x;q,t;t_0,t_1,t_2,t_3) d\T.
\]
Similarly, if $t_0t_1=t^{1-n} q^{-m}$, then
\[
I^{(n)}_K(f;q,t;t_0,t_1,t_2,t_3)
=
\frac{
\langle f\rangle_{qR}}
{\langle 1\rangle_{qR}}.
\]
\end{cor}

As the virtual integral is a rational function of the parameters, we
also conclude:

\begin{cor}
For all partitions $\lambda$, $\mu$,
\[
I^{(n)}_K(
K^{(n)}_\lambda(;q,t;t_0,t_1,t_2,t_3)
K^{(n)}_\mu(;q,t;t_0,t_1,t_2,t_3);q,t;t_0,t_1,t_2,t_3)
=
\delta_{\lambda\mu} N_\lambda(;q,t,t^n;t_0,t_1,t_2,t_3),
\]
with $N_\lambda$ as above.
\end{cor}

\begin{rem}
An alternative derivation of this (assuming evaluation symmetry) was given
in \cite{vanDiejen:1996}.
\end{rem}

\bigskip
The inversion formula for generalized binomial coefficients allows
us to invert the binomial formula.

\begin{thm}
For any partition $\lambda$, we have the expansion
\[
\bar{P}^{*(n)}_\lambda(;q,t,t_0)
=
\sum_{\mu\subset\lambda}
\binomI{\lambda}{\mu}_{q,t,t^{n-1} \sqrt{t_0t_1t_2t_3/q}}
\frac{k^0_\lambda(q,t,t^n;t_0{:}t_1,t_2,t_3)}
{k^0_\mu(q,t,t^n;t_0{:}t_1,t_2,t_3)}
K^{(n)}_\mu(;q,t;t_0,t_1,t_2,t_3)
\]
\end{thm}

This has two interesting consequences.  First, we obtain a connection
coefficient formula for Koornwinder polynomials, by expanding the
interpolation polynomial in the binomial formula using the inverse binomial
formula.

\begin{thm}
For any partitions $\kappa\subset\lambda$,
\begin{align}
\left[\frac{K^{(n)}_\kappa(;q,t;t_0,t_1,t_2,t_3)}{k^0_\kappa(q,t,t^n;t_0{:}t_1,t_2,t_3)}
\right]&
\frac{K^{(n)}_\lambda(;q,t;t_0,t'_1,t'_2,t'_3)}
{k^0_\lambda(q,t,t^n;t_0{:}t'_1,t'_2,t'_3)}
\notag\\&=
\sum_{\kappa\subset\mu\subset\lambda}
\binomQ{\lambda}{\mu}_{q,t,t^{n-1}\sqrt{t_0t'_1t'_2t'_3/q}}
\binomI{\mu}{\kappa}_{q,t,t^{n-1}\sqrt{t_0t_1t_2t_3/q}}
\frac{k^0_\mu(q,t,t^n;t_0{:}t_1,t_2,t_3)}
     {k^0_\mu(q,t,t^n;t_0{:}t'_1,t'_2,t'_3)}.
\end{align}
\end{thm}

\begin{rem}
If $t'_2=t_2$, $t'_3=t_3$, so only $t_1$ is changed, then the corresponding
sum for Askey-Wilson polynomials has a closed form evaluation; it turns out
that something similar is true for Koornwinder polynomials, in that the sum
can be evaluated in terms of the generalized binomial coefficients of
\cite{bctheta}.  See also Theorem \ref{thm:koorn:connt} below.
\end{rem}

Second, we obtain an integral formula generalizing Kadell's formula (and
the $q$-analogue) for the (normalized) integral of a Macdonald polynomial
over a Jacobi ensemble.

\begin{cor}
For any partition $\lambda$, one has the following virtual integral.
\[
I^{(n)}_K(\bar{P}^{*(n)}_\lambda(;q,t,t_0);q,t;t_0,t_1,t_2,t_3)
=
(-t_0 t^{n-1})^{-|\lambda|}
t^{2n(\lambda)}
q^{-n(\lambda')}
\frac{C^0_\lambda(t^n,t^{n-1} t_0 t_1,t^{n-1} t_0 t_2,t^{n-1} t_0 t_3;q,t)}
{C^-_\lambda(t;q,t)C^0_\lambda(t^{2n-2} t_0t_1t_2t_3;q,t)}
\]
\end{cor}

This allows us to prove the following integral representation for
${}_8W^{(n)}_7$ series.

\begin{thm}
Choose $t_0,t_1,t_2,t_3,t_4,q,t,u\in C$ such that $\max(|t_i|,|q|,|t|)<1$.
Then
\begin{align}
\label{eq:8W7int1}
I^{(n)}_K(\prod_{1\le i\le n} \frac{(u x_i,u/x_i;q)}{(t_4 x_i,t_4/x_i;q)}
;q,t;t_0,t_1,t_2,t_3)
&=
\prod_{0\le i<n}
\frac{(t^{-i} u' t'_0,t^{-i} u' t'_1,t^{-i} u' t'_2,t^{-i} t'_0 t'_1 t'_2 t'_4;q)}
{(t^{-i} t'_0 t'_4,t^{-i} t'_1 t'_4,t^{-i} t'_2 t'_4,t^{-i} u' t'_0 t'_1 t'_2;q)}\\*
&\phantom{{}={}}\quad
{}_8W^{(n)}_7(u' t'_0 t'_1 t'_2/q;t'_0 t'_1,t'_0
t'_2,t'_1 t'_2,u'/t'_3,u'/t'_4;q,t;t^{1-n} t'_3 t'_4),
\notag
\end{align}
where $t'_i=t^{(n-1)/2} t_i$, $u'=t^{(n-1)/2} u$.
\end{thm}

\begin{proof}
First, suppose that $t_4 = q^m u$.  Then
\begin{align}
\prod_{1\le i\le n} \frac{(u x_i,u/x_i;q)}{(t_4 x_i,t_4/x_i;q)}
&=
(-u)^{mn}
q^{n(n^m)}
\prod_{\substack{1\le i\le m\\1\le j\le n}}
(x_j+1/x_j-q^{i-1} u-q^{1-i}/u)\\
&=
(-u)^{mn}
q^{n(n^m)}
\bar{P}^{*(n)}_{m^n}(x;q,t,u)\\
&=
C^0_{m^n}(t^{n-1} t_0 u,u/t_0;q,t)
\sum_{\lambda\subset m^n}
\frac{q^{n(\lambda')}
(-t^{n-1} t_0 q)^{|\lambda|}
C^0_\lambda(q^{-m};q,t)\bar{P}^{*(n)}_{\lambda}(x;q,t,t_0)}
{C^0_\lambda(t^{n-1} t_0 u,t^{n-1} q^{1-m} t_0/u;q,t) C^-_\lambda(q;q,t)}
.
\end{align}
Integrating both sides, we find that
\begin{align}
I^{(n)}_K(\prod_{1\le i\le n} &\frac{(u x_i,u/x_i;q)}{(q^m u x_i,q^m u/x_i;q)}
;q,t;t_0,t_1,t_2,t_3)\notag\\*
&\qquad=
C^0_{m^n}(t^{n-1} t_0 u,u/t_0;q,t)
{}_4\Phi^{(n)}_3\left(
\genfrac{}{}{0pt}{}
{q^{-m},t^{n-1} t_0 t_1,t^{n-1} t_0 t_2,t^{n-1} t_0 t_3}
{t^{2n-2} t_0t_1t_2t_3,t^{n-1} t_0 u,t^{n-1} q^{1-m} t_0/u};q,t;q\right)\\
&\qquad=
\frac{C^0_{m^n}(t^{n-1} u t_0,t^{n-1} u t_1,t^{n-1} u t_2;q,t)}
{C^0_{m^n}(t^{2n-2} u t_0 t_1 t_2;q,t)}\notag\\*
&\qquad\phantom{{}={}}
{}_8W^{(n)}_7(t^{2n-2} u t_0 t_1 t_2/q;t^{n-1} t_0 t_1,t^{n-1} t_0
t_2,t^{n-1} t_1 t_2,u/t_3,q^{-m};q,t;q^m u t_3).
\end{align}
Thus the desired formula holds when $t_4=q^m u$.  The result follows from
the fact that both sides are manifestly analytic in the stated domain.
\end{proof}

If we exchange $t_2$ and $t_3$, the integral is clearly unchanged; we thus
obtain a transformation of ${}_8W^{(n)}_7$ series.  Similarly, permuting
the parameters of the ${}_8W^{(n)}_7$ leads to a transformation of the
integral.  The resulting symmetry groups are enlarged from $S_5$ to the
Weyl group $D_5$, and there are thus a total of three different expressions
for the integral, corresponding to the double cosets $S_5{\setminus}
D_5/S_5$.  The remaining two are:
\[
\begin{split}
\prod_{0\le i<n}&
\frac{
(t^{-i} u'/t'_4,t^{-i} u' t'_4,t^{-i} t'_0t'_1t'_2t'_4,t^{-i} t'_0t'_1t'_3t'_4,t^{-i} t'_0t'_2t'_3t'_4,t^{-i} t'_1t'_2t'_3t'_4;q)}
{
(t^{-i} t'_0t'_1t'_2t'_3,t^{-i} t'_0t'_4,t^{-i} t'_1t'_4,t^{-i} t'_2t'_4,t^{-i}
t'_3t'_4,t^{-i} t'_0t'_1t'_2t'_3t^{\prime 2}_4;q)}\\
&\qquad\qquad{}_8W^{(n)}_7(
t'_0t'_1t'_2t'_3t^{\prime 2}_4/q;
t'_0t'_4,
t'_1t'_4,
t'_2t'_4,
t'_3t'_4,
t'_0t'_1t'_2t'_3t'_4/u';
q,t;t^{1-n}u'/t'_4),
\label{eq:8W7int2}
\end{split}
\]
subject to the convergence condition $|t^{1-n}u'/t'_4|<1$, and
\[
\begin{split}
\prod_{0\le i<n}&
\frac{
(t^{-i} t'_0 u',t^{-i} t'_1 u',t^{-i} t'_2 u',t^{-i} t'_3 u',t^{-i} t'_4 u',t^{-i} t'_0t'_1t'_2t'_3t'_4/u';q)}
{
(t^{-i} t'_0 t'_4,t^{-i} t'_1 t'_4,t^{-i} t'_2 t'_4,t^{-i} t'_3 t'_4,t^{-i}
u^{\prime 2},t^{-i} t'_0t'_1t'_2t'_3;q)}\\
&\qquad\qquad{}_8W^{(n)}_7(
u^{\prime 2}/q;
u'/t'_0,u'/t'_1,u'/t'_2,u'/t'_3,u'/t'_4;q,t;t^{1-n} t'_0t'_1t'_2t'_3t'_4/u'),
\label{eq:8W7int3}
\end{split}
\]
subject to the convergence condition $|t^{1-n} t'_0t'_1t'_2t'_3t'_4/u'|<1$.
(Compare equations (6.3.7-9) of \cite{GasperG/RahmanM:1990}.)

\begin{cor}\label{cor:koorn:int_mn=nm}
Let $m$, $n$ be nonnegative integers.  Then
\[
\begin{split}
I^{(n)}_K(\prod_{1\le i\le n} &\frac{(t^{m/2} v x_i,t^{m/2}
v/x_i;q)}{(t^{-m/2} v x_i,t^{-m/2} v/x_i;q)}
;q,t;t^{m/2} t_0,t^{m/2} t_1,t^{m/2} t_2,t^{m/2} t_3)\\
&\phantom{t^{-m/2} v x_i}=
I^{(m)}_K(\prod_{1\le i\le m} \frac{(t^{n/2} v x_i,t^{n/2}
v/x_i;q)}{(t^{-n/2} v x_i,t^{-n/2} v/x_i;q)}
;q,t;t^{n/2} t_0,t^{n/2} t_1,t^{n/2} t_2,t^{n/2} t_3).
\end{split}
\]
\end{cor}

\begin{proof}
Expand the integrals via \eqref{eq:8W7int3}, and use the fact that
\[
{}_8W^{(n)}_7(a;b,c,d,e,t^m;q,t;z)=
{}_8W^{(m)}_7(a;b,c,d,e,t^n;q,t;z).
\]
\end{proof}

\begin{rem}
The Jacobi limit is implicitly used in 
\cite[Section 15.7]{ForresterBook}.
\end{rem}

\begin{cor}
Let $m$, $n$ be nonnegative integers, and set $u = t^{2n+2m-2} t_0t_1t_2t_3t_4$. Then
\[
\begin{split}
\prod_{0\le i<n}&
\frac{
(t^{i} t_0t_4,t^{i} t_1t_4,t^{i} t_2t_4,t^{i} t_3t_4;q)}
{
(t^{-i} u/t_0,t^{-i} u/t_1,t^{-i} u/t_2,t^{-i} u/t_3;q)}\\
&
\phantom{t^{-i} u/t_0,t^{-i} u/t_1}
I^{(n)}_K(\prod_{1\le i\le n} \frac{(t^{-m/2} u x_i,t^{-m/2} u/x_i;q)}{(t^{m/2} t_4
x_i,t^{m/2} t_4/x_i;q)}
;q,t;t^{m/2} t_0,t^{m/2} t_1,t^{m/2} t_2,t^{m/2} t_3)
\end{split}
\]
is symmetric in $m$, $n$.
\end{cor}

\begin{proof}
The same, but using \eqref{eq:8W7int2}.
\end{proof}

\begin{rems}
When $m=0$ and thus $u=t^{2n-2}t_0t_1t_2t_3t_4$, we conclude
\[
I^{(n)}_K(\prod_{1\le i\le n} \frac{(u x_i,u/x_i;q)}{(t_4 x_i,t_4/x_i;q)}
          ;q,t;t_0,t_1,t_2,t_3)
=
\prod_{0\le i<n}
\frac{
(t^{-i} u/t_0,t^{-i} u/t_1,t^{-i} u/t_2,t^{-i} u/t_3;q)}
{
(t^{i} t_0t_4,t^{i} t_1t_4,t^{i} t_2t_4,t^{i} t_3t_4;q)};
\]
aside from the normalization of the integral, this is 
Theorem 2.1 of \cite{GustafsonRA:1994}.
\end{rems}

\begin{rems}
Formally, there is an analogous result using \eqref{eq:8W7int1}; this is a
special case of equation \eqref{eq:koorn:intplethsym} below.
\end{rems}

\begin{rems}
If we were to ignore the constraints that certain ratios be integral powers
of $t$ in the above transformations, we would find that our integral has
symmetry group $D_6$.  It is unclear how to make this rigorous, however,
given the significant difficulties with convergence.
\end{rems}

\bigskip
We conclude with some miscellaneous results.

Combining the inversion identity for binomial coefficients, the duality
symmetry of binomial coefficients, and the Cauchy identity for
interpolation polynomials gives a Cauchy identity for Koornwinder
polynomials.

\begin{thm}\cite{MimachiK:2001}
For all integers $m,n\ge 0$,
\begin{align}
\sum_{\lambda\subset m^n}
(-1)^{mn-|\lambda|}
K^{(n)}_\lambda(x_1,\dots x_n;q,t;t_0,t_1,t_2,t_3)
K^{(m)}_{n^m-\lambda'}(y_1,\dots y_m;t,q;&t_0,t_1,t_2,t_3)\\*
&=
\prod_{1\le i\le n}
\prod_{1\le j\le m}
(x_i+1/x_i-y_j-1/y_j).\notag
\end{align}
\end{thm}

The action of the difference operators on the Koornwinder polynomials
(Lemma \ref{lem:koorn:ortho1}) is related via evaluation symmetry
to the following connection coefficient result.

\begin{thm}
For any partition $\lambda$,
\[
\frac{K^{(n)}_\lambda(;q,t;t_0,t_1,t_2,t_3)}
{k^0_\lambda(q,t,t^n;t_0{:}t_1,t_2,t_3)}
=
\sum_{\kappa\prec\lambda}
\psi^{(d)}_{\lambda/\kappa}(t^n t_0 t_1;q,t,t^{n-1}\sqrt{t_0t_1t_2t_3/q}) 
\frac{K^{(n)}_\kappa(;q,t;t_0,q t_1,t_2,t_3)}
{k^0_\kappa(q,t,t^n;t_0{:}q t_1,t_2,t_3)}
\]
\end{thm}

\begin{proof}
Apply the difference equation for binomial coefficients, with parameter
$u=t^n t_0 t_1$, and simplify the result using the fact that
\[
k^0_\mu(q,t,T;t_0{:}t_1,t_2,t_3)\psi^{(d)}_{\mu/\mu}(T t_0 t_1;q,t,(T/t)\sqrt{t_0t_1t_2t_3/q})
=
k^0_\mu(q,t,T;t_0{:}q t_1,t_2,t_3).
\]
\end{proof}

Dualizing (i.e., using \eqref{eq:koorn:duality}), or equivalently using the
integral equation, gives another special connection coefficient.

\begin{thm}\label{thm:koorn:connt}
For any partition $\lambda$,
\[
\frac{K^{(n)}_\lambda(;q,t;t_0,t_1 t,t_2,t_3)}
{k^0_\lambda(q,t,t^n;t_0{:}t_1 t,t_2,t_3)}
=
\sum_{\kappa'\prec\lambda'}
\psi^{(i)}_{\lambda/\kappa}(t^n t_0 t_1;q,t,t^{n-1}\sqrt{t_0t_1t_2t_3/q}) 
\frac{K^{(n)}_\kappa(;q,t;t_0,t_1,t_2,t_3)}
{k^0_\kappa(q,t,t^n;t_0{:}t_1,t_2,t_3)}
\]
\end{thm}

\begin{proof}
Take $u=t^n t_0 t_1$ in the integral equation, and note
\[
k^0_\mu(q,t,T;t_0{:}t t_1,t_2,t_3)\psi^{(i)}_{\mu/\mu}(T t_0 t_1;q,t,(T/t)\sqrt{t_0t_1t_2t_3/q})
=
k^0_\mu(q,t,T;t_0{:}t_1,t_2,t_3).
\]
\end{proof}

\begin{rem}
Note that this formula still depends nontrivially on $t$ in the
univariate case.  Since the Askey-Wilson polynomials are naturally
independent of $t$, we obtain the connection coefficient formula
\cite[Eqs.~6.4-5]{AskeyR/WilsonJ:1985} for
Askey-Wilson polynomials with only one parameter changed.
\end{rem}

We similarly obtain the following quasi-branching rule.

\begin{thm}\label{thm:koorn:brancht}
For any partition $\lambda$,
\[
\frac{K^{(n+1)}_\lambda(x_1,\dots x_n,t_0;q,t;t_0,t_1,t_2,t_3)}
{k^0_\lambda(q,t,t^{n+1};t_0{:}t_1,t_2,t_3)}
=
\sum_{\substack{\kappa'\prec\lambda'\\\ell(\kappa)\le n}}
\psi^{(i)}_{\lambda/\kappa}(t^{n+1};q,t,t^n\sqrt{t_0t_1t_2t_3/qt}) 
\frac{K^{(n)}_\kappa(x_1,\dots x_n;q,t;t_0 t,t_1,t_2,t_3)}
{k^0_\kappa(q,t,t^n;t_0 t{:}t_1,t_2,t_3)}
\]
\end{thm}

\begin{proof}
Applying Lemma \ref{lem:interp:dec_mn}, we find that it suffices to show
\[
\frac{\binomQ{\lambda}{\mu}_{q,t,t^n \sqrt{t_0t_1t_2t_3/q}}}
     {k^0_\mu(q,t,t^{n+1};t_0{:}t_1,t_2,t_3)}
=
\sum_{\substack{\kappa\subset\lambda\\\ell(\kappa)\le n}}
\psi^{(i)}_{\lambda/\kappa}(t^{n+1};q,t,t^{n-1}\sqrt{t_0t_1t_2t_3 t/q}) 
\frac{\binomQ{\kappa}{\mu}_{q,t,t^{n-1} \sqrt{t_0t_1t_2t_3 t/q}}}
     {k^0_\mu(q,t,t^n;t_0 t{:}t_1,t_2,t_3)},
\]
for $\ell(\mu)\le n$.  Since
\[
\psi^{(i)}_{\lambda/\kappa}(t^{n+1};q,t,t^{n-1}\sqrt{t_0t_1t_2t_3 t/q}) 
=
0
\]
for $\ell(\kappa)=n+1$, the restriction on the length of $\kappa$ can
be removed, at which point the integral equation may be applied.
\end{proof}

\begin{rems}
The special case of the integral equation used in the above proof
corresponds to the integral representation of \cite{OkounkovA:1998}.
\end{rems}

\begin{rems}
This implies via the Cauchy identity a quasi-Pieri identity of the
form
\[
\prod_{1\le i\le n} (x_i+1/x_i-t_0-1/t_0)
K^{(n)}_\mu(x_1,\dots x_n;q,t;qt_0,t_1,t_2,t_3)
=
\sum_{\substack{\lambda\succ\mu\\\ell(\lambda)\le n}}
c_{\lambda\mu}
K^{(n)}_\lambda(x_1,\dots x_n;q,t;t_0,t_1,t_2,t_3),
\]
for suitable coefficients $c_{\lambda\mu}$.
\end{rems}

If we combine the last two theorems, we obtain a special branching rule of the
form
\[
K^{(n+1)}_\lambda(x_1,\dots x_n,t_0;q,t;t_0,t_1,t_2,t_3)
=
\sum_{\mu'\prec \kappa'\prec\lambda'}
c_{\lambda/\kappa} d_{\kappa/\mu}
K^{(n)}_\mu(x_1,\dots x_n;q,t;t_0,t_1,t_2,t_3),
\]
for certain coefficients $c_{\lambda/\kappa}$, $d_{\kappa/\mu}$, thus
confirming the speculation in \cite{OkounkovA:1998} that the integral
representation is related to a branching rule for Koornwinder polynomials.
In general, it follows by combining the Pieri identities
\cite{vanDiejen:1996} and the Cauchy identity
that there exists a branching rule of the form
\[
K^{(n+1)}_\lambda(x_1,\dots x_n,u;q,t;t_0,t_1,t_2,t_3)
=
\sum_{\substack{\mu'\subset\lambda'\\\exists\kappa:\mu'\prec\kappa'\prec\lambda'}}
c_{\lambda/\mu}
K^{(n)}_\mu(x_1,\dots x_n;q,t;t_0,t_1,t_2,t_3),
\]
for suitable coefficients $c_{\lambda/\mu}\in \F(t_0,t_1,t_2,t_3)[u,1/u]$.
The expression one obtains for these coefficients is far from a closed
form, however.

\section{Symmetric functions from interpolation polynomials}\label{sec:sfinterp}

It turns out that the families of interpolation and Koornwinder polynomials
have natural lifts to families of symmetric functions; in each case, one
obtains two families with an additional algebraic parameter that reduce to
the given $BC_n$-symmetric polynomials when appropriately specialized.
We will first consider the interpolation polynomial case.

The first lifting involves inverting the natural projection from
$\Lambda$ to $\F[x_i^{\pm 1}]^{BC_n}$ given by
\[
f\mapsto f(x_1,1/x_1,x_2,1/x_2,\dots x_n,1/x_n).
\]
For any $n$, this map is surjective, but is quite far from injective, and
thus we have no way to define a unique lifting.  It turns out, however,
that if we introduce an algebraic parameter $T=t^n$, then there is a
unique lifting with coefficients in $\F(s,T)$.

The above homomorphism acts on power-sums by
\[
p_k\mapsto \sum_{1\le i\le n} x_i^k+x_i^{-k}.
\]
If we evaluate this at the partition $\lambda$, we have:
{\allowdisplaybreaks
\begin{align}
\sum_{1\le i\le n} x_i^k+x_i^{-k}
&=
\sum_{1\le i\le n} \left(q^{k\mu_i} t^{(n-i)k} s^k
+
q^{-k\mu_i} t^{-(n-i)k} s^{-k}\right)\\
&=
\sum_{1\le i\le \ell(\mu)}
\left((q^{k\mu_i}-1) t^{(n-i)k} s^k
+
(q^{-k\mu_i}-1) t^{-(n-i)k} s^{-k}\right)
+
\sum_{1\le i\le n}
\left(t^{(n-i)k} s^k + t^{-(n-i)k} s^{-k}\right)\\
&=
\sum_{1\le i\le \ell(\mu)}
\left((q^{k\mu_i}-1) t^{-ki} (sT)^k
+
(q^{-k\mu_i}-1) t^{ki} (sT)^{-k}\right)
+
s^k \frac{1-T^k}{1-t^k} + s^{-k} \frac{1-1/T^k}{1-1/t^k}.
\end{align}
}
This motivates the following definition.

\begin{defn}
The {\it lifted interpolation polynomials} are the unique family of
(inhomogeneous) symmetric functions $\tilde{P}^*_\lambda(;q,t,T;s)$ such
that
\begin{itemize}
\item $\tilde{P}^*_\lambda(;q,t,T;s)=m_\lambda+\text{dominated terms}$
\item $\tilde{P}^*_\lambda(\langle\mu\rangle_{q,t,T;s};q,t,T;s)=0$
for $\mu<\lambda$.
\end{itemize}
Here, for a symmetric function $f$,
$f(\langle\mu\rangle_{q,t,T;s})$
is its image under the homomorphism such that
\[
p_k(\langle\mu\rangle_{q,t,T;s})
=
\sum_{1\le i\le \ell(\mu)}
\left((q^{k\mu_i}-1) t^{-ki} (sT)^k
+
(q^{-k\mu_i}-1) t^{ki} (sT)^{-k}\right)
+
s^k \frac{1-T^k}{1-t^k} + s^{-k} \frac{1-1/T^k}{1-1/t^k}.
\]
\end{defn}

\begin{thm}
The lifted interpolation polynomials are well-defined, with coefficients
in $\F(s,T)$, and have the property that
\[
\tilde{P}^*_\lambda(x_1,x_1^{-1},x_2,x_2^{-1},\dots x_n,x_n^{-1};q,t,t^n;s)
=
\begin{cases}
\bar{P}^*_\lambda(x_1,\dots x_n;q,t,s) & \ell(\lambda)\le n\\
0&\ell(\lambda)>n.
\end{cases}
\]
Moreover,
\[
\tilde{P}^*_\lambda(\langle\mu\rangle_{q,t,T;s};q,t,T;s)=0
\]
unless $\lambda\subset\mu$, and
\[
\tilde{P}^*_\lambda(\langle\lambda\rangle_{q,t,T;s};q,t,T;s)=
(qsT/t)^{-|\lambda|} t^{n(\lambda)} q^{-2n(\lambda')}
C^-_\lambda(q;q,t)
C^+_\lambda((sT/t)^2;q,t)
\]
\end{thm}

\begin{proof}
Suppose that the claims are known for $\lambda<\kappa$.
Then we can write
\[
\tilde{P}^*_{\kappa}(;q,t,T;s)
=
m_\kappa
+
\sum_{\mu\le\kappa}
c_\mu
\tilde{P}^*_\mu(;q,t,T;s)
\]
for appropriate coefficients $c_\mu\in \F(s,T)$; the resulting equations
are triangular with nonzero diagonal by the inductive hypothesis, and thus
$\tilde{P}^*_\kappa$ is well-defined.  We cannot yet rule out a pole
at $T=t^n$, but can certainly conclude that only finitely many such
poles exist.

Now, for $n\ge \ell(\kappa)$ not hitting such a pole,
\[
\tilde{P}^*_\kappa(x_1,x_1^{-1},x_2,x_2^{-1},\dots x_n,x_n^{-1};q,t,t^n;s)
\]
is a $BC_n$-symmetric polynomial with leading monomial
$m_\kappa$ (since the natural projection is triangular on the monomial
functions) satisfying the necessary vanishing identities, so must
therefore equal the interpolation polynomial.

In particular, for any partition $\mu\not\supset\kappa$,
\[
\tilde{P}^*_\kappa(\langle\mu\rangle_{q,t,T;s};q,t,T;s)=0
\]
whenever $T=t^n$ for $n$ sufficiently large.  But then this identity
must in fact hold in $\F(s,T)$.  Similarly, the evaluation at $\mu=\kappa$
holds for sufficiently large $n$, and thus for all $T$.  In particular,
the diagonal coefficients in the equations for $c_\mu$ are nonzero in 
$\F(s)$ for any $T$ of the form $t^n$, and thus the coefficients have
no poles at such $T$.

Finally, we observe that for $n<\ell(\kappa)$,
\[
\tilde{P}^*_\kappa(x_1,x_1^{-1},x_2,x_2^{-1},\dots x_n,x_n^{-1};q,t,t^n;s)
\]
is a $BC_n$-symmetric polynomial vanishing at all partitions $\mu$ of
length at most $n$, and must therefore vanish.
\end{proof}

The next lemma indicates that the two parameters $s$ and $T$ give only one
degree of freedom, up to homomorphism.

\begin{lem}\label{lem:interp:homsT}
We have the plethystic identity
\[
\tilde{P}^*_\lambda(;q,t,T;s T')
=
\tilde{P}^*_\lambda([p_k+s^k \frac{T^k-T^{\prime k}}{1-t^k}+s^{-k}
\frac{1/T^k-1/T^{\prime k}}{1-1/t^k}];q,t,T';sT).
\]
\end{lem}

\begin{proof}
The right hand side has leading monomial $m_\lambda$ and vanishes at all
the appropriate partitions.
\end{proof}

\begin{cor}\label{cor:sfinterp:denoms}
The coefficients of $(sT)^{|\lambda|}\tilde{P}^*_\lambda(;q,t,T;s)$ lie in
$\F[s,T]$.
\end{cor}

\begin{proof}
The coefficients of $s^{|\lambda|}\bar{P}^*_\lambda(;q,t,t^n;s)$ must lie
in $\F[s]$ for $n$ sufficiently large, since the same is true of
$s^{|\lambda|}\bar{P}^{*(n)}_\lambda(;q,t,s)$.  In particular, the
coefficients of $(sT)^{|\lambda|}\bar{P}^*_\lambda(;q,t,t^n;s T/t^n)$ lie in
$\F[sT]$; applying the homomorphism
\[
p_k\mapsto p_k-s^k \frac{t^{nk}-T^k}{1-t^k}-s^{-k} \frac{t^{-nk}-T^{-k}}{1-t^{-k}}
\]
at most enlarges the coefficient ring to $\F[s,T]$, and produces
$(sT)^{|\lambda|}\tilde{P}^*_\lambda(;q,t,T;s)$.
\end{proof}

\begin{rem}
In fact, an examination of $\tilde{P}^*_\lambda(\langle\lambda\rangle)$ and
the homomorphism $\langle\lambda\rangle_{q,t,T;s}$ further shows that the
only possible denominator factors are $s$, $T$, $q$, $t$, and $(1-q^i t^j)$
for $i$, $j$ nonnegative integers, not both 0.
\end{rem}

\bigskip
We now define a slight modification of the Macdonald involution.  Recall
that the Macdonald involution is the homomorphism
\[
\omega_{q,t}:p_k\mapsto (-1)^{k-1} \frac{1-q^k}{1-t^k} p_k.
\]
We rescale this slightly, to give:
\[
\tilde{\omega}_{q,t}: p_k\mapsto (-1)^{k-1}
\frac{q^{k/2}-q^{-k/2}}{t^{k/2}-t^{-k/2}} p_k;
\]
thus this acts on ordinary Macdonald polynomials as
\[
\tilde{\omega}_{q,t} P_\mu(;q,t)
=
b_\mu(q,t)^{-1} (t/q)^{|\mu|/2} P_{\mu'}(;t,q).
\]

\begin{lem}
For any symmetric function $f$ and partition $\mu$,
\[
(\tilde{\omega}_{q,t}f)(\langle \mu\rangle_{t,q,1/T;-\sqrt{qt}/s})
=
f(\langle\mu'\rangle_{q,t,T;s}).
\]
\end{lem}

\begin{proof}
It suffices to prove this for the power-sum functions $p_k$.
Now, we find:
\begin{align}
(-1)^{k-1}
\frac{q^{k/2}-q^{-k/2}}{t^{k/2}-t^{-k/2}}
\sum_{1\le i\le \ell(\mu)}
(t^{k\mu_i}-1)q^{-ki} (-\sqrt{qt}/sT)^k
&=
\frac{1-q^k}{1-t^{-k}}
\sum_{1\le i\le \ell(\mu)}
(t^{k\mu_i}-1)(sT)^{-k}\\
&=
(1-q^k) \sum_{(i,j)\in\mu} t^{jk} (sT)^{-k}\\
&=
(1-q^k) \sum_{(j,i)\in\mu'} t^{jk} q^{-ki} (sT)^{-k}\\
&=
\sum_{1\le j\le \ell(\mu')} (q^{-\mu'_j}-1) t^{jk} (sT)^{-k}.
\end{align}
and
\[
(-1)^{k-1}
\frac{q^{k/2}-q^{-k/2}}{t^{k/2}-t^{-k/2}}
(-\sqrt{qt}/s)^k \frac{1-T^{-k}}{1-q^k}
=
s^{-k} \frac{1-T^{-k}}{1-t^{-k}}.
\]
The other two terms simplify analogously, and thus
\[
(\tilde{\omega}_{q,t}p_k)(\langle \mu\rangle_{t,q,1/T;-\sqrt{qt}/s})
=
p_k(\langle\mu'\rangle_{q,t,T;s})
\]
as required.
\end{proof}

This gives us the following result:

\begin{thm}\label{thm:interp:sfduality}
The interpolation polynomials satisfy
\[
\tilde{\omega}_{q,t} \tilde{P}_\mu(;q,t,T;s)
=
b_\mu(q,t)^{-1} (t/q)^{|\mu|/2} \tilde{P}_{\mu'}(;t,q,1/T;-\sqrt{qt}/s).
\]
\end{thm}

\begin{proof}
If we evaluate the left-hand side at
$\langle\nu\rangle_{t,q,1/T;-\sqrt{qt}/s}$,
the result is
\[
\tilde{P}_\mu(\langle\nu\rangle_{q,t,T;s};q,t,T;s),
\]
and thus vanishes unless $\nu\supset\mu$.  It follows that
\[
\tilde{\omega}_{q,t} \tilde{P}_\mu(;q,t,T;s)
\propto
\tilde{P}_{\mu'}(;t,q,1/T;-\sqrt{qt}/s),
\]
with some nonzero constant.  Now, by triangularity, we can write
\[
\tilde{P}_{\mu'}(;t,q,1/T;-\sqrt{qt}/s)
=
P_{\mu'}(;t,q)+\sum_{\nu<\mu'} c_{\nu} P_\nu(;t,q);
\]
similarly, applying $\tilde{\omega}_{qt}$ to $\tilde{P}_\mu$,
we have
\[
\tilde{P}_{\mu'}(;t,q,1/T;-\sqrt{qt}/s)
=
\sum_{\nu\le \mu} c'_{\nu} P_{\nu'}(;t,q).
\]
Since conjugation reverses dominance for partitions of the same size, we
find that the degree $|\mu|$ portion of $\tilde{P}_{\mu'}$ is $P_{\mu'}$, and
similarly for $\tilde{P}_\mu$.  The constant then follows immediately.
\end{proof}

\begin{rem}
In particular, we obtain a new proof of the leading term limit for
interpolation polynomials.
\end{rem}

\begin{cor}\label{cor:bc_duality2}
For any partitions $\mu$ and $\lambda$,
\[
\tilde{P}_\mu(\langle\lambda\rangle_{q,t,T;s};q,t,T;s)
=
b_\mu(q,t)^{-1} (t/q)^{|\mu|/2}
\tilde{P}_{\mu'}(\langle\lambda'\rangle_{t,q,1/T;-\sqrt{qt}/s};t,q,1/T;-\sqrt{qt}/s).
\]
In particular,
\begin{align}
\binomQ{\lambda}{\mu}_{q,t,s}
&=
\binomQ{\lambda'}{\mu'}_{t,q,1/\sqrt{qt}s}\\
\binomI{\lambda}{\mu}_{q,t,s}
&=
\binomI{\lambda'}{\mu'}_{t,q,1/\sqrt{qt}s}.
\end{align}
\end{cor}

\bigskip
Our second lifting of interpolation polynomials to symmetric functions
involves the observation that for $\lambda\subset m^n$,
\[
(\prod_{1\le i\le n} x_i^m)\bar{P}^{*(n)}_{m^n-\lambda}(x_1,\dots x_n;q,t,s)
\]
is an $S_n$-symmetric polynomial (without negative exponents).  Moreover,
we have the following:

\begin{lem}
For arbitrary positive integers $m,n$, and an arbitrary partition
$\lambda\subset m^n$,
\[
\lim_{x_n\to 0}
(\prod_{1\le i\le n} x_i^m)\bar{P}^{*(n)}_{m^n-\lambda}(x_1,\dots x_n;q,t,s)
=
\begin{cases}
(\prod_{1\le i\le n-1} x_i^m)\bar{P}^{*(n-1)}_{m^{n-1}-\lambda}(x_1,\dots
x_{n-1};q,t,s)
& \lambda_n=0\\
0& \lambda_n>0.
\end{cases}
\]
\end{lem}

\begin{proof}
We apply the branching rule; the only term that does not vanish in the
limit is the unique term in which the degree has been reduced by $m$,
namely
\[
\bar{P}^{*(n-1)}_{m^{n-1}-\lambda}(x_1,\dots
x_{n-1};q,t,s).
\]
\end{proof}

With this in mind, we define:

\begin{defn}
The {\it virtual interpolation polynomials} are the unique symmetric
functions $\hat{P}^*_{\lambda}(;q,t,Q;s)\in
\hat\Lambda$ with coefficients in $\F(Q,s)$ such that
\[
\hat{P}^*_{\lambda}(x_1,\dots x_n;q,t,Q;s)
=
\prod_{1\le i\le n}
\frac{x_i^m (s x_i,q^{m+1}x_i/s Q;q)}{(q x_i/s,sQ x_i/q^m;q)}
\bar{P}^{*(n)}_{m^n-\lambda}(x_1,\dots x_n;q,t,s Q/q^m)
\]
for all positive integers $m$, $n$ such that $\lambda\subset m^n$.
\end{defn}

We immediately see that these are well-defined and that the leading term
(now the term of smallest degree) is again the Macdonald polynomial.  The
two kinds of interpolation polynomials are related via the Cauchy identity:

\begin{thm}
We have the following identity of symmetric functions in $\Lambda_x\otimes
\hat\Lambda_y$:
\[
\sum_\lambda
(-1)^{|\lambda|}
\tilde{P}^*_{\lambda}(x;q,t,T;s)
\hat{P}^*_{\lambda'}(y;t,q,T;s)
=
\sum_\lambda
(-1)^{|\lambda|}
P_\lambda(x;q,t)
P_{\lambda'}(y;t,q)
\]
\end{thm}

\begin{proof}
Fix integers $m$, $n>0$.  Then we have:
{\allowdisplaybreaks
\begin{align}
\sum_\lambda
(-1)^{|\lambda|}
\tilde{P}^*_{\lambda}(x_1^{\pm 1},x_2^{\pm 1}&,\dots x_n^{\pm 1};q,t,t^n;s)
\hat{P}^*_{\lambda'}(y_1,y_2,\dots y_m;t,q,t^n;s)\notag\\*
&=
\sum_{\lambda\subset m^n}
(-1)^{|\lambda|}
\bar{P}^*_{\lambda}(x_1,x_2,\dots x_n;q,t,s)
\prod_{1\le j\le m} y_j^n
\bar{P}^*_{n^m-\lambda'}(y_1,y_2,\dots y_m;t,q,s)\\
&=
\prod_{1\le j\le m} y_j^n
\prod_{\substack{1\le i\le n\\1\le j\le m}}
(y_j+1/y_j-x_i-1/x_i)\\
&=
\prod_{\substack{1\le i\le n\\1\le j\le m}}
(1-y_j x_i)(1-y_j/x_i)\\
&=
\sum_\lambda
(-1)^{|\lambda|}
P_\lambda(x_1^{\pm 1},x_2^{\pm 1},\dots x_n^{\pm 1};q,t)
P_{\lambda'}(y_1,\dots y_m;t,q)
\end{align}
}
The desired identity follows by rationality of coefficients.
\end{proof}

\begin{cor}
\[
\tilde{\omega}_{q,t} \hat{P}_{\lambda}(;q,t,Q;s)
=
b_\lambda(q,t)^{-1} (t/q)^{|\lambda|/2}
\hat{P}_{\lambda'}(;t,q,1/Q;-\sqrt{qt}/s)
\]
\end{cor}

\begin{proof}
In the sum
\[
\sum_\lambda
(-1)^{|\lambda|}
\tilde{P}^*_{\lambda'}(x;t,q,Q;s)
\hat{P}^*_{\lambda}(y;q,t,Q;s)
=
\sum_\lambda
(-1)^{|\lambda|}
P_{\lambda'}(x;t,q)
P_\lambda(y;q,t),
\]
apply $\tilde{\omega}_{q,t}$ to the $y$ variables and
$\tilde{\omega}_{t,q}$ to the $x$ variables.  We find:
\[
\sum_\lambda
(-1)^{|\lambda|}
P_{\lambda'}(x;q,t)
P_\lambda(y;t,q)
=
\sum_\lambda
(-1)^{|\lambda|}
b_{\lambda'}(t,q)^{-1} (q/t)^{|\lambda|/2}
\tilde{P}^*_\lambda(x;q,t,1/Q;-\sqrt{qt}/s)
(\tilde{\omega}_{q,t}\hat{P}^*_\lambda(y;q,t,Q;s)),
\]
and thus, taking coefficients of
$\tilde{P}^*_{\lambda}(x;q,t,1/Q;-\sqrt{qt}/s)$ on both sides:
\[
\hat{P}^*_{\lambda'}(y;t,q,1/Q;-\sqrt{qt}/s)
=
b_{\lambda'}(t,q)^{-1} (q/t)^{|\lambda|/2}
(\tilde{\omega}_{q,t}\hat{P}^*_{\lambda}(y;q,t,Q;s)).
\]
The desired result is immediate.
\end{proof}

Applying the involution to only one of the sets of variables gives
another Cauchy identity:

\begin{cor}
\[
\sum_\lambda
b_\lambda(q,t) (q/t)^{|\lambda|/2}
\tilde{P}^*_{\lambda}(x;q,t,T;s)
\hat{P}^*_{\lambda}(y;q,t,1/T;\sqrt{qt}/s)
=
\sum_\lambda
(q/t)^{|\lambda|/2}
P_\lambda(x;q,t)
Q_\lambda(y;q,t)
\]
\end{cor}

\begin{cor}\label{cor:interp:expmac2}
The virtual interpolation polynomials can be expanded in terms of Macdonald
polynomials as follows:
\[
\hat{P}^*_\mu(;q,t,Q;s)
=
\sum_\lambda
(-1)^{|\lambda|-|\mu|}
\left([\tilde{P}^*_{\mu'}(;t,q,Q;s)] P_{\lambda'}(;t,q)\right)
P_\lambda(;q,t)
\]
\end{cor}

\bigskip
The bulk branching rule lifts to the following result.

\begin{thm}
For any partition $\lambda$,
\[
\tilde{P}^{*}_\lambda
    ([p_k+\frac{u^k-v^k}{1-t^k}+\frac{u^{-k}-v^{-k}}{1-t^{-k}}]
     ;q,t,T v/u;s)
=
\sum_{\mu\subset\lambda}
\psi^{(B)}_{\lambda/\mu}(u,v;q,t;s T)
\tilde{P}^{*}_\mu(;q,t,T;s),
\]
where $\psi^{(B)}_{\lambda/\mu}(u,v;q,t;s)$ is as in Theorem
\ref{thm:interp:bulk_branch}.
Similarly, for any partition $\lambda$, we have the following identity in
$\hat\Lambda$ with coefficients in $\F(s)[[u,v]]$.
\[
\hat{P}^*_\lambda([p_k+\frac{u^k-v^k}{1-t^k}];q,t,Q;s)
=
\hat{P}^*_0([\frac{u^k-v^k}{1-t^k}];q,t,Q;s)
\sum_{\mu\subset\lambda}
\hat{\psi}^{(B)}_{\lambda/\mu}(u,v;q,t;s Q)
\hat{P}^*_\mu(;q,t,Q;s),
\]
where
\begin{align}
\hat{P}^*_0([\frac{u^k-v^k}{1-t^k}];q,t,Q;s)
&=
\prod_{j\ge 0}
\frac{(t^j u s,t^j qu/s Q,t^j v s Q,t^j qv/s;q)}
{(t^j v s,t^j qv/s Q,t^j u s Q,t^j qu/s;q)}\\
\hat{\psi}^{(B)}_{\lambda/\mu}(u,v;q,t;s)
&=
\frac{C^0_\mu(qu/ts;q,t)C^0_\mu(sv/q;1/q,1/t)}
{C^0_\lambda(qv/ts;q,t)C^0_\lambda(su/q;1/q,1/t)}
P_{\lambda/\mu}([(u^k-v^k)/(1-t^k)];q,t).
\end{align}
\end{thm}

Corollary \ref{cor:interp:rect1} (evaluation at a ``constant'') has the
following especially pleasing lift.

\begin{cor}
For any partition $\lambda$,
\begin{align}
(-\sqrt{xyzt})^{|\lambda|}
\tilde{P}^*_\lambda
    ([\frac{(x^{k/2}-x^{-k/2})(y^{k/2}-y^{-k/2})(z^{k/2}-z^{-k/2})}
           {t^{k/2}-t^{-k/2}}];q,t,1;&\sqrt{t/xyz})
=\notag\\
&t^{-2n(\lambda)} q^{n'(\lambda)}
\frac{C^0_\lambda(x,y,z;1/q,1/t)}{C^-_\lambda(1/t;1/q,1/t)}
\end{align}
\end{cor}

\begin{rem}
Corollary \ref{cor:interp:rect1} is the special case $y=t^n$, $z=1/(t^{n-1}
x s^2)$.
\end{rem}

We also note the lifted versions of the connection coefficient identity.

\begin{thm}\label{thm:sfinterp:conn}
For any partitions $\lambda$, $\mu$,
\begin{align}
[\tilde{P}^*_\mu(;q,t,T;s)]
\tilde{P}^*_\lambda(;q,t,T;s')
&=
\frac{
C^0_\lambda(T;q,t)
C^0_\lambda(t/Ts s';1/q,1/t)}
{
C^0_\mu(T;q,t)
C^0_\mu(t/Ts s';1/q,1/t)}
P_{\lambda/\mu}([\frac{s^k-s^{\prime k}}{1-t^k}];q,t).\\
[\hat{P}^*_\lambda(;q,t,Q;s)]\hat{P}^*_\mu(;q,t,Q;s')
&=
\frac{
C^0_\lambda(Q;1/q,1/t)
C^0_\lambda(t/Qs s';q,t)}
{
C^0_\mu(Q;1/q,1/t)
C^0_\mu(t/Qs s';q,t)}
Q_{\lambda/\mu}([\frac{s^k-s^{\prime k}}{1-t^k}];q,t)
\end{align}
\end{thm}

We also have a lift of the bulk Pieri identity, albeit only to the lifted
polynomials (for the virtual polynomials, there are convergence problems,
even formally).

\begin{thm}
For any partition $\mu$, the following identity holds in $(\Lambda\otimes
\F(s))[[u,v]]$.
\begin{align}
\Bigl(
\sum_\kappa Q_\kappa(\left[\frac{u^k-v^k}{1-t^k}\right];q,t)&P_\kappa(;q,t)
\Bigr)
\tilde{P}^{*}_\mu(;q,t,T;s)\\*
&=
\Bigl(
\sum_\kappa Q_\kappa(\left[\frac{u^k-v^k}{1-t^k}\right];q,t)
P_\kappa(\langle\mu\rangle_{q,t,T;s};q,t)
\Bigr)
\sum_{\lambda\supset\mu}
\psi^{(P)}_{\lambda/\mu}(u,v;q,t;s T)
\tilde{P}^{*}_\lambda(;q,t,T;s),\notag
\end{align}
where $\psi^{(P)}_{\lambda/\mu}(u,v;q,t;s)$ is as in Theorem
\ref{thm:interp:bulk_Pieri}.
\end{thm}

\bigskip
We close with the following refinement of the fact that the leading terms
of the interpolation polynomials are Macdonald polynomials.

\begin{thm}\label{thm:PsPtriangular}
The lifted and virtual interpolation polynomials are triangular in the
Macdonald polynomial basis, with respect to the inclusion partial order.
That is, we have the expansions
\begin{align}
\tilde{P}^*_\lambda(;q,t,T;s)
&=
P_\lambda(;q,t)+
\sum_{\mu\subset\lambda} c_{\lambda/\mu} P_\mu(;q,t)\\
\hat{P}^*_\mu(;q,t,T;s)
&=
P_\mu(;q,t)+
\sum_{\lambda\supset\mu} \hat{c}_{\lambda/\mu} P_\lambda(;q,t)
\end{align}
for suitable constants $c$.
\end{thm}

\begin{proof}
We first note that by Corollary \ref{cor:interp:expmac2}, the
first claim implies the second.  Now, suppose the first claim is false,
and let $\lambda$ be an inclusion-minimal partition such that
$\tilde{P}^*_\lambda(;q,t,T;s)$ is not triangular.  Furthermore, let
$\mu$ be an inclusion-maximal partition not contained in $\lambda$ such
that
\[
[P_\mu(;q,t)]\tilde{P}^*_\lambda(;q,t,T;s)\ne 0.
\]
We will show that this coefficient is independent of $T$ and $s$; that it
is 0 (giving a contradiction) will be shown in the proof of Theorem
\ref{thm:KPtriangular} below.

By Theorem \ref{thm:sfinterp:conn}, for any fixed $T\ne 0$, the lifted
interpolation polynomials are mutually triangular with respect to the
inclusion ordering, and thus the minimality of $\lambda$ implies that the
coefficient is independent of $s$.  On the other hand,
\[
[P_\mu(;q,t)] \tilde{P}^*_\lambda(;q,t,T;s)
=
[P_\mu(;q,t)] 
\tilde{P}^*_\lambda([p_k+s^k \frac{T^k-1}{1-t^k}+s^{-k}
\frac{1/T^k-1}{1-1/t^k}];q,t,1;sT);
\]
since this homomorphism is triangular in the Macdonald basis, the
maximality of $\mu$ implies that
\[
[P_\mu(;q,t)] \tilde{P}^*_\lambda(;q,t,T;s)
=
[P_\mu(;q,t)] \tilde{P}^*_\lambda(;q,t,1;sT),
\]
and is thus independent of $s$ and $T$ as required.
\end{proof}

\begin{cor}
For any integer $n\ge 0$,
\[
\tilde{P}^*_{1^n}(;q,t,T;s)=
(e_n-e_{n-2})([
p_k-s^k \frac{1-(T/t^{n-1})^k}{1-t^k}-
     s^{-k} \frac{1-(t^{n-1}/T)^k}{1-1/t^k}]),
\]
where $e_{-1}=e_{-2}=0$.
\end{cor}

\begin{proof}
It suffices to prove this in the case $T=t^{n-1}$.  Triangularity
then tells us that
\[
\tilde{P}^*_{1^n}(;q,t,t^{n-1};s)
=
e_n+\sum_{0\le m<n} c_m e_m
\]
for suitable coefficients $c_m$.  On the other hand, we know that
\[
\tilde{P}^*_{1^n}(x_1^{\pm 1},x_2^{\pm 1},\dots x_{n-1}^{\pm
1};q,t,t^{n-1};s)=0.
\]
The only algebraic relation satisfied by the quantities
\[
e_m(x_1^{\pm 1},x_2^{\pm 1},\dots x_{n-1}^{\pm 1})
\]
for $0\le m\le n$ is that $e_n=e_{n-2}$, and thus $e_n-e_{n-2}$ is
the only symmetric function satisfying both requirements.
\end{proof}

\section{Symmetric functions from Koornwinder polynomials}\label{sec:sfkoorn}

Via the binomial formula, the lifted interpolation polynomials lead
immediately to a lifting for Koornwinder polynomials.

\begin{defn}
The {\it lifted Koornwinder polynomials} are defined by the expansion
\[
\tilde{K}_\lambda(;q,t,T;t_0,t_1,t_2,t_3)
=
\sum_{\mu\subset\lambda}
\binomQ{\lambda}{\mu}_{q,t,(T/t)\sqrt{t_0t_1t_2t_3/q}}
\frac{k^0_\lambda(q,t,T;t_0{:}t_1,t_2,t_3)}
     {k^0_\mu(q,t,T;t_0{:}t_1,t_2,t_3)}
\tilde{P}^*_\mu(;q,t,T;t_0).
\]
\end{defn}

\begin{thm}
For any integer $n>0$ and partition $\lambda$, and for generic values
of the parameters,
\[
\tilde{K}_\lambda(x_1^{\pm 1},\dots x_n^{\pm 1};q,t,t^n;t_0,t_1,t_2,t_3)
=
\begin{cases}
K^{(n)}_\lambda(x_1,\dots x_n;q,t;t_0,t_1,t_2,t_3)
&
\ell(\lambda)\le n\\
0
&
\text{otherwise.}
\end{cases}
\]
\end{thm}

\begin{proof}
The claim when $\ell(\lambda)\le n$ is immediate from the binomial formula
for ordinary Koornwinder polynomials.  Thus assume $\ell(\lambda)>n$,
and consider the term of the lifted binomial formula corresponding to
a partition $\mu$.  The only factors that can lead to a zero or pole
at $T=t^n$ are the factors $C^0_\lambda(T;q,t)$ of $k^0_\lambda$,
$C^0_\mu(T;q,t)$ of $k^0_\mu(T;q,t)$, and the lifted interpolation
polynomial itself.  Now, when $\ell(\mu)\le n$, $C^0_\mu(t^n;q,t)\ne 0$
while $C^0_\lambda(t^n;q,t)=0$, and thus the $\mu$ term vanishes.  On the
other hand, when $\ell(\mu)>n$, we find that
\[
\lim_{T\to t^n} C^0_\lambda(T;q,t)/C^0_\mu(T;q,t)
\]
is well-defined and nonzero.  In this case, however, the interpolation
polynomial itself vanishes.  Thus all terms in the expansion vanish,
as required.
\end{proof}

\begin{rem}
The genericity hypothesis is necessary when $\ell(\lambda)>n$:
\[
\lim_{T\to t} \tilde{K}_{1^2}(x_1,1/x_1;q,t,T;1,-1,\sqrt{t},-\sqrt{t})
=
1,
\]
not 0.  However, as long as $T$ is specialized before any of the other
parameters, this will not be a problem.  Indeed, by the following
corollary, the only possible problems (for generic $q$ and $t$) arise when
\[
C^+_\lambda(t^{2n-2} t_0t_1t_2t_3/q;q,t)=0.
\]
\end{rem}

\begin{cor}
For any partition $\lambda$, the only possible factors of the denominators
of the coefficients of the symmetric function
\[
t^{2|\lambda|+3n(\lambda)}
C^+_\lambda((T/t)^2 t_0t_1t_2t_3/q;q,t)
\tilde{K}_\lambda(;q,t,T;t_0,t_1,t_2,t_3)
\]
are $t$ and binomials of the form $1-q^i t^j$ for $i,j\ge 0$.
\end{cor}

\begin{proof}
If we specialize $T=t^n$ for any sufficiently large $n$, then we obtain an
ordinary Koornwinder polynomial, in which the only denominator factors that
can appear are of the form $1-q^i t^j$, $i,j\ge 0$.  Since this is
true for all sufficiently large $n$, the result follows.
\end{proof}

\begin{rem}
We conjecture that, in fact,
\[
t^{2|\lambda|+3n(\lambda)}
C^-_\lambda(t;q,t)
C^+_\lambda((T/t)^2 t_0t_1t_2t_3/q;q,t)
\tilde{K}_\lambda([p_k/(1-t)];q,t,T;t_0,t_1,t_2,t_3)
\]
has coefficients in $\Z[q,t,t_0,t_1,t_2,t_3,T]$.  (This is true for the
leading terms, as in that case it reduces to the corresponding integrality
result for Macdonald polynomials.)
\end{rem}

\begin{prop}
For any pair of partitions $\mu$, $\lambda$,
\[
\frac{
\tilde{K}_\lambda(\langle\mu\rangle_{q,t,T;t_0};q,t,T;t_0,t_1,t_2,t_3)
}{
k^0_\lambda(q,t,T;t_0{:}t_1,t_2,t_3)
}
=
\frac{
\tilde{K}_\mu(\langle\lambda\rangle_{q,t,T;\hat{t}_0};q,t,T;\hat{t}_0,\hat{t}_1,\hat{t}_2,\hat{t}_3)
}{
k^0_\mu(q,t,T;\hat{t}_0{:}\hat{t}_1,\hat{t}_2,\hat{t}_3),
}
\]
where
\[
\hat{t}_0=\sqrt{t_0t_1t_2t_3/q}\text{;\qquad}
\hat{t}_i=t_0t_i/\hat{t}_0,\text{ $i\in \{1,2,3\}$}.
\]
\end{prop}

The symmetries of ordinary Koornwinder polynomials lift; in addition,
we obtain new symmetries involving $T$.

\begin{prop}\label{prop:sfkoorn:trivsim}
For any partition $\lambda$,
\begin{align}
\tilde{K}_\lambda(;q,t,T;t_0,t_1,t_2,t_3)
&=
\tilde{K}_\lambda(;1/q,1/t,1/T;1/t_0,1/t_1,1/t_2,1/t_3)\\
\tilde{K}_\lambda(;q,t,T;t_0,t_1,t_2,t_3)
&=
(-1)^{|\lambda|} \tilde{K}_\lambda([(-1)^k p_k];q,t,T;-t_0,-t_1,-t_2,-t_3)\\
\tilde{K}_\lambda(;q,t,T;t_0,t_1,t_2,t_3)
&=
\tilde{K}_\lambda(
[p_k+\frac{(t/t_0)^k+(t/t_1)^k-t_0^k-t_1^k}{(1-t^k)}]
;q,t,T t_0 t_1/t;t/t_1,t/t_0,t_2,t_3)\label{eq:koorn:plethsym}\\
\tilde{\omega}_{q,t}(\tilde{K}_\lambda(;q,t,T;t_0,t_1,t_2,t_3))
&=
b_\lambda(q,t)^{-1}
(t/q)^{|\lambda|/2}
\tilde{K}_{\lambda'}(;t,q,1/T;\frac{-\sqrt{qt}}{t_0},\frac{-\sqrt{qt}}{t_1},\frac{-\sqrt{qt}}{t_2},\frac{-\sqrt{qt}}{t_3}).\label{eq:koorn:duality}
\end{align}
Furthermore, $\tilde{K}_\lambda$ is invariant under permutations of $t_0$,
$t_1$, $t_2$, $t_3$.
\end{prop}

\begin{rem}
We note three particularly nice special cases of \eqref{eq:koorn:plethsym}:
\begin{align}
\tilde{K}_\lambda(;q,t,T;\sqrt{t},-\sqrt{t},t_2,t_3)
&=
\tilde{K}_\lambda(;q,t,-T;-\sqrt{t},\sqrt{t},t_2,t_3)\\
\tilde{K}_\lambda(;q,t,T;t,\sqrt{t},t_2,t_3)
&=
\tilde{K}_\lambda(
[p_k+1]
;q,t,T \sqrt{t};\sqrt{t},1,t_2,t_3)\\
\tilde{K}_\lambda(;q,t,T;-t,-\sqrt{t},t_2,t_3)
&=
\tilde{K}_\lambda(
[p_k+(-1)^k]
;q,t,T \sqrt{t};-\sqrt{t},-1,t_2,t_3).
\end{align}
\end{rem}

\begin{defn}
The {virtual Koornwinder integral} $I_K(;q,t,T;t_0,t_1,t_2,t_3)$ is
the linear functional on symmetric functions defined by
\[
I_K(f;q,t,T;t_0,t_1,t_2,t_3)
=
[\tilde{K}_0(;q,t,T;t_0,t_1,t_2,t_3)] f.
\]
\end{defn}

In particular, we note that when $T=t^n$, this reduces to the virtual
Koornwinder integral defined above:
\[
I^{(n)}_K(;q,t;t_0,t_1,t_2,t_3)=I_K(;q,t,t^n;t_0,t_1,t_2,t_3)
\]
Again, the specialization of $T$ must occur before any other
specialization.  For instance, the identity
\[
I^{(1)}_K(e_2;q,t;t_0,t_1,t_2,t_3) = 1
\]
holds for all values of the parameters.  On the other hand,
\[
\lim_{T\to t} I_K(e_2;q,t,T;\pm 1,\pm \sqrt{t}) = 0.
\]
Again, as long as $q$ and $t$ are generic and $C^+_\lambda(t^{2n-2}
t_0t_1t_2t_3/q;q,t)\ne 0$, there is no problem.

Since the set of parameters with $T=t^n$ is Zariski dense, the
orthogonality of ordinary Koornwinder polynomials lifts.

\begin{prop}
\[
I_K(\tilde{K}_\lambda(;q,t,T;t_0,t_1,t_2,t_3)
\tilde{K}_\mu(;q,t,T;t_0,t_1,t_2,t_3);q,t,T;t_0,t_1,t_2,t_3)
=
\delta_{\lambda\mu} N_\lambda(;q,t,T;t_0,t_1,t_2,t_3),
\]
with $N_\lambda$ as above.
\end{prop}

The symmetries of Proposition \ref{prop:sfkoorn:trivsim} carry over to the
virtual integral.

\begin{cor}
For any partition $\lambda$ and any symmetric function $f$,
\begin{align}
I_K(f;1/q,1/t,1/T;1/t_0,1/t_1,1/t_2,1/t_3)
&=
I_K(f;q,t,T;t_0,t_1,t_2,t_3)
\\
I_K(f;q,t,T;-t_0,-t_1,-t_2,-t_3)
&=
I_K(f([(-1)^k p_k]);q,t,T;t_0,t_1,t_2,t_3)
\\
I_K(f;
q,t,T t_0 t_1/t;t/t_1,t/t_0,t_2,t_3)
&=
I_K(f([p_k-\frac{t_0^k+t_1^k-(t/t_0)^k-(t/t_1)^k}{(1-t)}]);q,t,T;t_0,t_1,t_2,t_3)\label{eq:koorn:intplethsym}
\\
I_K(f;t,q,1/T;\frac{-\sqrt{qt}}{t_0},\frac{-\sqrt{qt}}{t_1},\frac{-\sqrt{qt}}{t_2},\frac{-\sqrt{qt}}{t_3})
&=
I_K(\tilde{\omega}_{q,t} f;q,t,T;t_0,t_1,t_2,t_3).
\end{align}
Furthermore, $I_K$ is symmetric in $t_0$, $t_1$, $t_2$, $t_3$.
\end{cor}

The next several results are immediate lifts of the analogous results
above for ordinary Koornwinder polynomials.

\begin{prop}
For any partition $\lambda$, we have the expansion
\[
\tilde{P}^*_\lambda(;q,t,T;t_0)
=
\sum_{\mu\subset\lambda}
\binomI{\lambda}{\mu}_{q,t,(T/t)\sqrt{t_0t_1t_2t_3/q}}
\frac{k^0_\lambda(q,t,T;t_0{:}t_1,t_2,t_3)}
{k^0_\mu(q,t,T;t_0{:}t_1,t_2,t_3)}
\tilde{K}_\mu(;q,t,T;t_0,t_1,t_2,t_3)
\]
\end{prop}

\begin{prop}
For any partitions $\kappa\subset\lambda$,
\begin{align}
\left[\frac{\tilde{K}_\kappa(;q,t,T;t_0,t_1,t_2,t_3)}{k^0_\kappa(q,t,T;t_0{:}t_1,t_2,t_3)}
\right]&
\frac{\tilde{K}_\lambda(;q,t,T;t_0,t'_1,t'_2,t'_3)}
{k^0_\lambda(q,t,T;t_0{:}t'_1,t'_2,t'_3)}
\notag\\&=
\sum_{\kappa\subset\mu\subset\lambda}
\binomQ{\lambda}{\mu}_{q,t,(T/t)\sqrt{t_0t'_1t'_2t'_3/q}}
\binomI{\mu}{\kappa}_{q,t,(T/t)\sqrt{t_0t_1t_2t_3/q}}
\frac{k^0_\mu(q,t,T;t_0{:}t_1,t_2,t_3)}
     {k^0_\mu(q,t,T;t_0{:}t'_1,t'_2,t'_3)}.
\end{align}
\end{prop}

\begin{prop}
For any partition $\lambda$, one has the following virtual integral.
\[
I_K(\tilde{P}^*_\lambda(;q,t,T;t_0);q,t,T;t_0,t_1,t_2,t_3)
=
(-t_0 T/t)^{-|\lambda|}
t^{2n(\lambda)}
q^{-n(\lambda')}
\frac{C^0_\lambda(T,T t_0 t_1/t,T t_0 t_2/t,T t_0 t_3/t;q,t)}
{C^-_\lambda(t;q,t)C^0_\lambda(T^2 t_0t_1t_2t_3/t^2;q,t)}
\]
\end{prop}

\begin{prop}
For any partition $\lambda$,
\[
\frac{\tilde{K}_\lambda(;q,t,T;t_0,t_1,t_2,t_3)}
{k^0_\lambda(q,t,T;t_0{:}t_1,t_2,t_3)}
=
\sum_{\kappa\prec\lambda}
\psi^{(d)}_{\lambda/\kappa}(T t_0 t_1;q,t,(T/t)\sqrt{t_0t_1t_2t_3/q}) 
\frac{\tilde{K}_\kappa(;q,t,T;t_0,q t_1,t_2,t_3)}
{k^0_\kappa(q,t,T;t_0{:}q t_1,t_2,t_3)}
\]
\end{prop}

\begin{prop}\label{prop:sfkoorn:connt}
For any partition $\lambda$,
\[
\frac{\tilde{K}_\lambda(;q,t,T;t_0,t_1 t,t_2,t_3)}
{k^0_\lambda(q,t,T;t_0{:}t_1 t,t_2,t_3)}
=
\sum_{\kappa'\prec\lambda'}
\psi^{(i)}_{\lambda/\kappa}(T t_0 t_1;q,t,(T/t)\sqrt{t_0t_1t_2t_3/q}) 
\frac{\tilde{K}_\kappa(;q,t,T;t_0,t_1,t_2,t_3)}
{k^0_\kappa(q,t,T;t_0{:}t_1,t_2,t_3)}
\]
\end{prop}

\begin{prop}\label{prop:sfkoorn:brancht}
For any partition $\lambda$,
\[
\frac{\tilde{K}_\lambda([p_k+t_0^k+t_0^{-k}];q,t,t T;t_0,t_1,t_2,t_3)}
{k^0_\lambda(q,t,t T;t_0{:}t_1,t_2,t_3)}
=
\sum_{\kappa'\prec\lambda'}
\psi^{(i)}_{\lambda/\kappa}(t T;q,t,(T/t)\sqrt{t_0t_1t_2t_3 t/q})
\frac{\tilde{K}_\kappa(;q,t,T;t_0 t,t_1,t_2,t_3)}
{k^0_\kappa(q,t,T;t_0 t{:}t_1,t_2,t_3)}
\]
\end{prop}

\bigskip
We next turn to the Cauchy identities.

\begin{defn}
The {\it virtual Koornwinder polynomials} are the basis of $\hat\Lambda$
given by
\[
\hat{K}_\lambda(;q,t,Q;t_0,t_1,t_2,t_3)
=
\sum_{\mu\supset\lambda}
(-1)^{|\mu|-|\lambda|}
\binomI{\mu}{\lambda}_{q,t,1/Q\hat{t}_0}
\frac{k^0_{\mu'}(t,q,Q;t_0{:}t_1,t_2,t_3)}
{k^0_{\lambda'}(t,q,Q;t_0{:}t_1,t_2,t_3)}
\hat{P}^*_\mu(;q,t,Q;t_0),
\]
where $\hat{t}_0=\sqrt{t_0t_1t_2t_3/q}$.
\end{defn}

The name is justified by the following result, a straightforward
verification from the binomial formula.

\begin{thm}
If $\lambda\subset m^n$, then
\[
\hat{K}_\lambda(x_1,\dots x_n;q,t,q^m;t_0,t_1,t_2,t_3)
=
\prod_{1\le i\le n} x_i^m
K^{(n)}_{m^n-\lambda}(x_1,\dots x_n;q,t;t_0,t_1,t_2,t_3)
\]
\end{thm}

The Cauchy identities follow immediately from the inversion formula for
binomial coefficients and the Cauchy identities for interpolation
polynomials.

\begin{thm}
We have the following identities in $\Lambda_x\otimes\hat\Lambda_y$.
\[
\sum_\lambda
(-1)^{|\lambda|}
\tilde{K}_\lambda(x;q,t,T;t_0,t_1,t_2,t_3)
\hat{K}_{\lambda'}(y;t,q,T;t_0,t_1,t_2,t_3)
=
\sum_\lambda
(-1)^{|\lambda|}
P_\lambda(x;q,t)
P_{\lambda'}(y;t,q)\]
\begin{align}
\sum_\lambda
b_\lambda(q,t)
(\frac{q}{t})^{|\lambda|/2}
\tilde{K}_{\lambda}(x;q,t,T;t_0,t_1,t_2,t_3)
\hat{K}_{\lambda}(y;q,t,1/T;&\frac{\sqrt{qt}}{t_0},\frac{\sqrt{qt}}{t_1},\frac{\sqrt{qt}}{t_2},\frac{\sqrt{qt}}{t_3})\notag\\*
&=
\sum_\lambda
(\frac{q}{t})^{|\lambda|/2}
P_\lambda(x;q,t)
Q_\lambda(y;q,t)
\end{align}
\end{thm}

The special case $\hat{K}_{0}$ is of particular interest, as
a result of the following integrals (obtained by integrating
the Cauchy identities):

\begin{cor}
\begin{align}
I^{(n)}_K(
\prod_{\substack{1\le i\le n\\1\le j\le m}} (1-y_j x_i)(1-y_j/x_i)
;q,t;t_0,t_1,t_2,t_3)
&=
\hat{K}_{0}(y_1,\dots y_m;t,q,t^n;t_0,t_1,t_2,t_3)\\
I^{(n)}_K(
\prod_{\substack{1\le i\le n\\1\le j\le m}} 
\frac{
(\sqrt{qt} y_j x_i^{\pm 1};q)}
{
(\sqrt{q/t} y_j x_i^{\pm 1};q)}
;q,t;t_0,t_1,t_2,t_3)
&=
\hat{K}_{0}(y_1,\dots y_m;q,t,t^{-n};\frac{\sqrt{qt}}{t_0},\frac{\sqrt{qt}}{t_1},\frac{\sqrt{qt}}{t_2},\frac{\sqrt{qt}}{t_3})\label{eq:koorn:intK0b}
\end{align}
\end{cor}

\begin{cor}
The following identity holds for all sufficiently small $y$.
\[
\begin{split}
\hat{K}_{0}(y;q,t,Q;t_0,t_1,t_2,t_3)
&=
\frac{(q y/t_0Q,q y/t_1Q,q y/t_2Q,q^2 y/t_0 t_1 t_2 Q;q)}
{(q y/t_0,q y/t_1,q y/t_2,q^2 y/t_0t_1t_2 Q^2;q)}\\*
&\phantom{{}={}}\qquad\qquad
{}_8W_7(\frac{qy}{t_0t_1t_2 Q^2};
\frac{q}{t_0t_1 Q},
\frac{q}{t_0t_2 Q},
\frac{q}{t_1t_2 Q},
y t_3,
1/Q;q,q y/t_3)
\end{split}
\]
\end{cor}

\begin{proof}
Applying Corollary \ref{cor:koorn:int_mn=nm} to equation
\eqref{eq:koorn:intK0b}, we see that the identity holds whenever
$Q=t^{-n}$, and thus for all $Q$.
\end{proof}

\bigskip
The case $T=0$ of the lifted and virtual Koornwinder polynomials turns out
to be especially nice.  This is somewhat surprising, considering that the
binomial formula behaves very badly in that case: a rather large amount of
cancellation is required to eliminate the apparent singularities at $T=0$.
The key to dealing with the $T=0$ case turns out to be the inner product.

\begin{defn}
Let $\mu$ and $\sigma$ be sequences of complex numbers, such that
$\Re(\sigma_j)>0$ for all $j$.  Then the {Gaussian functional}
$I_G(;\mu;\sigma)$ is the linear functional on symmetric functions defined
by
\[
I_G(f;\mu;\sigma)
=
\int_{\R^{\deg(f)}}
f
\prod_{1\le j\le \deg(f)}
(2\pi\sigma_j)^{-1/2}
e^{-(p_j-\mu_j)^2/2\sigma_j}
dp_j;
\]
in particular, $I_G(1;\mu;\sigma)=1$.
\end{defn}

This has a probabilistic interpretation: $I_G(f;\mu;\sigma)$ is the
expected value of $f$ if the power sum functions $p_k$ are independent and
normally distributed random variables with mean $\mu_k$ and variance
$\sigma_k$.  In particular, $I_G(f;\mu;\sigma)$ is polynomial in $\mu$
and $\sigma$; we extend it to arbitrary $\mu$ and $\sigma$ accordingly.

Our reason for introducing Gaussian functionals is the following theorem.

\begin{thm}
For any symmetric function $f$,
\[
I_K(f;q,t,0;t_0,t_1,t_2,t_3)
=
I_G(f;\mu;\sigma),
\]
where
\begin{align}
\mu_{2k-1}&=\frac{t_0^{2k-1}+t_1^{2k-1}+t_2^{2k-1}+t_3^{2k-1}}{1-t^{2k-1}}\\
\mu_{2k}&=\frac{t_0^{2k}+t_1^{2k}+t_2^{2k}+t_3^{2k}-1-t^k-q^k-(qt)^k}{1-t^{2k}}\\
\sigma_k &= \frac{1-t^k}{1-q^k}.
\end{align}
\end{thm}

\begin{proof}
Since $I_K(f)$ is a rational function of the parameters for any $f$,
it suffices to consider the limit
\begin{align}
\lim_{n\to\infty}
Z_n^{-1}
\int f(x_i^{\pm 1})
w^{(n)}_K(x;q,t;t_0,t_1,t_2,t_3) d\T
\end{align}
after expanding the integrand as a formal power series in the parameters.

Now, if $q=t=t_0=t_1=t_2=t_3=0$, $w^{(n)}_K$ reduces to the distribution
function for the eigenvalues of a (Haar-distributed) random symplectic
matrix.  It follows from Theorem 6 of \cite{DiaconisP/ShahshahaniM:1994}
(see also
\cite[Section 8]{BaikJ/RainsEM:2001}) that
\[
\frac{\int f(x_i^{\pm 1})
w^{(n)}_K(x;0,0;0,0,0,0) d\T} {\int w^{(n)}_K(x;0,0;0,0,0,0) d\T} =
I_G(f;\mu^{(0)},\sigma^{(0)}),
\]
where $\mu^{(0)}_{2k-1}=0$, $\mu^{(0)}_{2k}=-1$, $\sigma^{(0)}_k=k$, for
all sufficiently large $n$.  Thus for any power series $f$ with
coefficients in $\Lambda$, we have the formal limit
\[
\lim_{n\to \infty} I^{(n)}_K(f;0,0;0,0,0,0)
=
I_G(f;\mu^{(0)},\sigma^{(0)}).
\]

Now, we have
\begin{align}
\frac{w^{(n)}_K(x;q,t;t_0,t_1,t_2,t_3)}
{w^{(n)}_K(x;0,0;0,0,0,0)}
=&\\
(q;q)^{-n}&
\prod_{1\le k}
\exp\left(
-
\frac{q^k+t^k}{2k(1-q^k)}
p_{2k}(x_i^{\pm 1})
+
\frac{t_0^k+t_1^k+t_2^k+t_3^k}{k(1-q^k)}
p_k(x_i^{\pm 1})
-
\frac{q^k-t^k}{2k(1-q^k)}
p_k^2(x_i^{\pm 1})\right).\notag
\end{align}
This combines with the Gaussian density from $I_G(f;\mu^{(0)},\sigma^{(0)})$ to give
the desired Gaussian density.
\end{proof}

In the sequel, we will evaluate the Gaussian integral via the following
expansions.

\begin{lem}\label{lem:sfkoorn:gaussian}
We have the following integrals:
\begin{align}
I_K(\prod_{j,k} \frac{(t x_j y_k;q)}{(x_j y_k;q)};q,t,0;t_0,t_1,t_2,t_3)_y
&=
\prod_{j<k} \frac{(t x_j x_k;q)}{(x_j x_k;q)}
\prod_j \frac{(t x_j^2;q)}{(t_0 x_j,t_1 x_j,t_2 x_j,t_3 x_j;q)}\label{eq:sfkoorn:gaussian1}\\
I_K(\prod_{j,k} (1+x_j y_k);q,t,0;t_0,t_1,t_2,t_3)_y
&=
\prod_{j<k} \frac{(q x_j x_k;t)}{(x_j x_k;t)}
\prod_j
\frac{(-t_0 x_j,-t_1 x_j,-t_2 x_j,-t_3 x_j;t)}{(x_j^2;t)}.\label{eq:sfkoorn:gaussian2}
\end{align}
\end{lem}

\begin{proof}
Complete the square in the Gaussian integrals.
\end{proof}

The lifted Koornwinder polynomials for $T=0$ are given by the following
generating function.

\begin{thm}
We have the following identity in $\hat\Lambda_x\otimes \Lambda_y$.
\begin{align}
\sum_\lambda
Q_\lambda(x;q,t)
\tilde{K}_\lambda(y;q,t,0;t_0,t_1,t_2,t_3)
=
\prod_{j,k} \frac{(t x_j y_k;q)}{(x_j y_k;q)}
\prod_{j<k} \frac{(x_j x_k;q)}{(t x_j x_k;q)}
\prod_j \frac{(t_0 x_j,t_1 x_j,t_2 x_j,t_3 x_j;q)}{(t x_j^2;q)}
\end{align}
\end{thm}

\begin{proof}
It is equivalent to show that the polynomials
\[
K'_\lambda(y)
=
[Q_\lambda(x;q,t)]
\prod_{j,k} \frac{(t x_j y_k;q)}{(x_j y_k;q)}
\prod_{j<k} \frac{(x_j x_k;q)}{(t x_j x_k;q)}
\prod_j \frac{(t_0 x_j,t_1 x_j,t_2 x_j,t_3 x_j;q)}{(t x_j^2;q)}
\]
are orthogonal with respect to the inner product
$I_K(;q,t,0;t_0,t_1,t_2,t_3)$.  Evaluating the integral
\[
I_K(
\prod_{j,k} \frac{(t x_j y_k,t z_j y_k;q)}{(x_j y_k,z_j y_k;q)}
\prod_{j<k} \frac{(x_j x_k,z_j z_k;q)}{(t x_j x_k,t z_j z_k;q)}
\prod_j \frac{(t_0 x_j,t_1 x_j,t_2 x_j,t_3 x_j,t_0 z_j,t_1 z_j,t_2 z_j,t_3 z_j;q)}{(t x_j^2,t z_j^2;q)}
;q,t,0;t_0,t_1,t_2,t_3)_y
\]
using Lemma \ref{lem:sfkoorn:gaussian} produces the result
\[
\prod_{j,k} \frac{(t x_j z_k;q)}{(x_j z_k;q)}
=
\sum_\lambda Q_\lambda(x;q,t) P_\lambda(z;q,t);
\]
it follows that
\[
I_K(K'_\lambda K'_\mu;q,t,0;t_0,t_1,t_2,t_3)
=
\delta_{\lambda\mu} b_\lambda(q,t)^{-1}.
\]
\end{proof}

\begin{rems}
Note, in particular, that
\[
N_\lambda(q,t,0;t_0,t_1,t_2,t_3) = b_\lambda(q,t)^{-1}
\]
as required.
\end{rems}

\begin{rems}
The special cases
\begin{align}
(q,t;t_0,t_1,t_2,t_3)&\mapsto (q,q;1,-1,\sqrt{q},-\sqrt{q})\\
(q,t;t_0,t_1,t_2,t_3)&\mapsto (q,q;q,-q,\sqrt{q},-\sqrt{q})
\end{align}
correspond to the Cauchy identities for orthogonal and symplectic
characters, respectively; in those cases, the lifted
Koornwinder polynomials are independent of $T$, and equal to the
``virtual characters'' surveyed in section 5 of \cite{SundaramS:1990}.
Similar comments apply to the $T=0$ results below.
\end{rems}

By duality, we obtain:

\begin{cor}
We have the following identity in $\hat\Lambda_x\otimes\Lambda_y$.
\[
\sum_\lambda
P_{\lambda'}(x;t,q)
\tilde{K}_\lambda(y;q,t,0;t_0,t_1,t_2,t_3)
=
\prod_{j,k}(1+x_j y_k)
\prod_{j<k} \frac{(x_j x_k;t)}{(q x_j x_k;t)}
\prod_j
\frac{(x_j^2;t)}{(t_0 x_j,t_1 x_j,t_2 x_j,t_3 x_j;t)}.
\]
\end{cor}

Comparing this with the Cauchy identity, we find:

\begin{cor}
The virtual Koornwinder polynomials are, when $T=0$, given by
the formula
\[
\hat{K}_{\lambda}(x;q,t,0;t_0,t_1,t_2,t_3)
=
\prod_j \frac{(t_0 x_j,t_1 x_j,t_2 x_j,t_3 x_j;q)}{(x_j^2;q)}
\prod_{j<k} \frac{(t x_j x_k;q)}{(x_j x_k;q)}
P_\lambda(x;q,t)
\]
\end{cor}

\begin{proof}
Here we use the fact that if $x$ is a formal variable, then
$(x;q)(x/q;1/q) = 1$.
\end{proof}

The generating function also gives us the following branching rule.

\begin{thm}
For any partition $\lambda$, we have the following identity in 
$\Lambda_x\otimes \Lambda_y$.
\[
\tilde{K}_\lambda(x,y;q,t,0;t_0,t_1,t_2,t_3)
=
\sum_{\mu\subset\lambda}
P_{\lambda/\mu}(x;q,t)
\tilde{K}_\mu(y;q,t,0;t_0,t_1,t_2,t_3).
\]
\end{thm}

\begin{proof}
We have
\[
\sum_\lambda Q_\lambda(z;q,t) \tilde{K}_\lambda(x,y;q,t,0;t_0,t_1,t_2,t_3)
=
\sum_\mu Q_\mu(z;q,t) P_\mu(x)
\sum_\lambda Q_\lambda(z;q,t) \tilde{K}_\lambda(y;q,t,0;t_0,t_1,t_2,t_3).
\]
Taking coefficients of $Q_\kappa(z;q,t)$ and
$\tilde{K}_\nu(y;q,t,0;t_0,t_1,t_2,t_3)$ gives
\begin{align}
[\tilde{K}_\nu(y;q,t,0;t_0,t_1,t_2,t_3)]
\tilde{K}_\kappa(x,y;q,t,0;t_0,t_1,t_2,t_3)
&=
\sum_\mu [Q_\kappa(z;q,t)](Q_\mu(z;q,t)Q_\nu(z;q,t))
P_\mu(x)\\
&=P_{\kappa/\nu}(x).
\end{align}
\end{proof}

Setting $y=0$ gives us the following expansion.

\begin{cor}
For any partition $\lambda$,
\[
\tilde{K}_\lambda(;q,t,0;t_0,t_1,t_2,t_3)
=
\sum_{\mu\subset\lambda}
P_{\lambda/\mu}(;q,t)
\tilde{K}_\mu(0;q,t,0;t_0,t_1,t_2,t_3).
\]
\end{cor}

The following is then immediate from the fact that
$[P_\kappa]P_{\lambda/\mu}=0$ unless $\kappa\subset\lambda$.

\begin{cor}\label{cor:T0:KPexp}
When $T=0$, the lifted and virtual Koornwinder polynomials are triangular
in the Macdonald polynomial basis, with respect to the inclusion partial
order.
\end{cor}

This allows us to finish the proof of Theorem \ref{thm:PsPtriangular}, as
well as prove a corresponding result for Koornwinder polynomials.

\begin{thm}\label{thm:KPtriangular}
The lifted and virtual Koornwinder polynomials are triangular in the
Macdonald polynomial basis, with respect to the inclusion partial order.
\end{thm}

\begin{proof}
Let $\lambda$ and $\mu$ be chosen as in the proof of Theorem
\ref{thm:PsPtriangular}.  Recall that we had shown there that the coefficient
\[
[P_\mu(;q,t)]\tilde{P}^*_\lambda(;q,t,T;s)
\]
is independent of $s$ and $T$, and still need to show that it is 0.  From
the binomial formula, it follows that
\[
[P_\mu(;q,t)]\tilde{K}_\lambda(;q,t,T;t_0,t_1,t_2,t_3)
=
[P_\mu(;q,t)]\tilde{P}^*_\lambda(;q,t,T;t_0)
\]
It thus suffices to show that the left-hand-side is 0 for some value of the
parameters.  Taking $T=0$ suffices, by Corollary \ref{cor:T0:KPexp}, and
thus Theorem \ref{thm:PsPtriangular} holds.  The present result follows
from the binomial formulas.
\end{proof}

The branching rule also implies the following expansion.

\begin{cor}
For any partition $\lambda$,
\[
P_\lambda(;q,t)
=
\sum_{\mu\subset\lambda}
I_K(P_{\lambda/\mu}(x;q,t);q,t,0;t_0,t_1,t_2,t_3)
\tilde{K}_\mu(;q,t,0;t_0,t_1,t_2,t_3)
\]
\end{cor}

\begin{proof}
We observe that
\begin{align}
\sum_{\mu\subset\lambda}
P_{\lambda/\mu}(x;q,t)
\tilde{K}_\mu(y;q,t,0;t_0,t_1,t_2,t_3)
&=
\tilde{K}_\mu(x,y;q,t,0;t_0,t_1,t_2,t_3)\\
&=
\sum_{\mu\subset\lambda}
P_{\lambda/\mu}(y;q,t)
\tilde{K}_\mu(x;q,t,0;t_0,t_1,t_2,t_3).
\end{align}
Integrating the $y$ variables gives the desired result.
\end{proof}

The fact that the $T=0$ branching rule is independent of $t_0$, $t_1$,
$t_2$, $t_3$ is related to the following plethystic symmetry, which follows
easily from the generating function.

\begin{thm}
For any partition $\lambda$,
\[
\tilde{K}_\lambda(;q,t,0;t'_0,t'_1,t'_2,t'_3)
=
\tilde{K}_\lambda(\left[p_k+\frac{t_0^k+t_1^k+t_2^k+t_3^k-t_0^{\prime
k}-t_1^{\prime k}-t_2^{\prime k}-t_3^{\prime
k}}{1-t^k}\right];q,t,0;t_0,t_1,t_2,t_3).
\]
In particular,
\[
\tilde{K}_\lambda(;q,t,0;t'_0,t_1,t_2,t_3)
=
\sum_{\mu\subset\lambda}
P_{\lambda/\mu}(\left[\frac{t^k_0-t^{\prime k}_0}{1-t^k}\right];q,t)
\tilde{K}_\mu(;q,t,0;t_0;t_1,t_2,t_3).
\]
\end{thm}

\begin{rem}
The special case $t'_0=t_0 t$ gives
\[
\tilde{K}_\lambda(;q,t,0;t_0 t,t_1,t_2,t_3)
=
\sum_{\mu\subset\lambda}
\psi_{\lambda/\mu}(q,t) t_0^{|\lambda/\mu|}
\tilde{K}_\mu(;q,t,0;t_0;t_1,t_2,t_3),
\]
a limiting case of Theorem \ref{thm:koorn:connt}.  Similarly,
\begin{align}
\tilde{K}_\lambda([p_k+t_0^k+t_0^{-k}];q,t,0;t_0,t_1,t_2,t_3)
&=
\sum_{\mu\subset\lambda}
\psi_{\lambda/\mu}(q,t) t_0^{-|\lambda/\mu|}
\tilde{K}_\mu([p_k+t_0^k];q,t,0;t_0,t_1,t_2,t_3)\\
&=
\sum_{\mu\subset\lambda}
\psi_{\lambda/\mu}(q,t) t_0^{-|\lambda/\mu|}
\tilde{K}_\mu(;q,t,0;t_0 t,t_1,t_2,t_3),
\end{align}
a limiting case of Theorem \ref{thm:koorn:brancht}.  (In fact, these
identities are what originally suggested that something like those theorems
should be true.)  Compare also the plethystic symmetry of lifted
Koornwinder polynomials (equation \eqref{eq:koorn:plethsym}).
\end{rem}

The Pieri identities are also nice when $T=0$.

\begin{thm}
For any partition $\lambda$ and integer $n\ge 0$,
\begin{align}
(\sum_n u^n g_n)
\tilde{K}_\lambda(;q,t,0;&t_0,t_1,t_2,t_3)\\*
&=
\frac{(tu^2;q)}{(t_0 u,t_1 u,t_2 u,t_3 u;q)}
\sum_{\substack{\nu\prec\lambda\\\mu\succ\nu}}
u^{|\lambda|-2|\nu|+|\mu|}
\psi_{\lambda/\nu}(q,t)
\varphi_{\mu/\nu}(q,t)
\tilde{K}_\mu(;q,t,0;t_0,t_1,t_2,t_3)
\notag
\end{align}
\end{thm}

\begin{proof}
Let $f_\lambda$ denote the left-hand side, and consider the generating
function
\[
\sum_{\lambda,\mu}
([\tilde{K}_\mu(;q,t,0;t_0,t_1,t_2,t_3)]f_\lambda)
Q_\lambda(x;q,t) P_\mu(y;q,t).
\]
Writing
\[
[\tilde{K}_\mu(;q,t,0;t_0,t_1,t_2,t_3)]f_\lambda
=
b_\mu(q,t)^{-1}
I_K(\tilde{K}_\mu(;q,t,0;t_0,t_1,t_2,t_3)f_\lambda;q,t,0;t_0,t_1,t_2,t_3),
\]
we find that the generating function can be written as a Gaussian integral.
Applying Lemma \ref{lem:sfkoorn:gaussian}, we find
\begin{align}
\sum_{\lambda,\mu}
([\tilde{K}_\mu(;q,t,0;t_0,t_1,t_2,t_3)]f_\lambda)
&Q_\lambda(x;q,t) P_\mu(y;q,t)\notag\\*
&=
\frac{(t u^2;q)}{(t_0 u,t_1 u,t_2 u,t_3 u;q)}
\prod_{j} \frac{(t u x_j,t u y_j;q)}{(u x_j,u y_j;q)}
\sum_\nu Q_\nu(x;q,t) P_\nu(y;q,t).\\
&=
\frac{(t u^2;q)}{(t_0 u,t_1 u,t_2 u,t_3 u;q)}
\sum_\nu
\sum_{\lambda\succ\nu}
\sum_{\mu\succ\nu}
\psi_{\lambda/\nu}(q,t) Q_\lambda(x;q,t) \varphi_{\mu/\nu}(q,t) P_\mu(y;q,t).
\end{align}
Comparing coefficients of $Q_\lambda(x;q,t)P_\mu(y;q,t)$ on both sides
gives the desired result.
\end{proof}

Similarly,

\begin{thm}
For any partition $\lambda$ and integer $n\ge 0$,
\begin{align}
(\sum_n u^n e_n)
\tilde{K}_\lambda(;q,t,0&;t_0,t_1,t_2,t_3)\\*
&=
\frac{(-t_0 u,-t_1 u,-t_2 u,-t_3 u;t)}{(u^2;t)}
\sum_{\substack{\nu'\prec\lambda'\\\mu'\succ\nu'}}
u^{|\lambda|-2|\nu|+|\mu|}
\varphi'_{\lambda/\nu}(q,t)
\psi'_{\mu/\nu}(q,t)
\tilde{K}_\mu(;q,t,0;t_0,t_1,t_2,t_3).\notag
\end{align}
\end{thm}

\section{Vanishing conjectures}\label{sec:vanish}

If we substitute $(t_0,t_1,t_2,t_3)\mapsto (\pm \sqrt{t},\pm \sqrt{qt})$
in equation \eqref{eq:sfkoorn:gaussian1}, the right-hand side becomes
\[
\prod_{j<k} \frac{(t x_j x_k;q)}{(x_j x_k;q)}.
\]
This is the right-hand side of a generalized Littlewood identity due to
Macdonald \cite[Ex. VI.7.4]{MacdonaldIG:1995}; we thus obtain the following
proposition.

\begin{prop}
For any partition $\lambda$,
\[
I_K(P_\lambda(;q,t);q,t,0;\pm\sqrt{t},\pm\sqrt{qt})=0
\]
unless $\lambda$ is of the form $\mu^2$, in which case
\[
I_K(P_{\mu^2}(;q,t);q,t,0;\pm\sqrt{t},\pm\sqrt{qt})=
\frac{C^-_\mu(qt;q,t^2)}{C^-_\mu(t^2;q,t^2)}.
\]
\end{prop}

The dual identity is also related to an integral.

\begin{prop}
For any partition $\lambda$,
\[
I_K(P_\lambda(;q,t);q,t,0;\pm 1,\pm\sqrt{t})=0
\]
unless $\lambda$ is of the form $2\mu$, in which case
\[
I_K(P_{2\mu}(;q,t);q,t,0;\pm 1,\pm\sqrt{t})=
\frac{C^-_\mu(q;q^2,t)}{C^-_\mu(t;q^2,t)}.
\]
\end{prop}

This leads us to formulate the following conjectures.

\begin{conj}\label{conj:vanish_U/Sp_T}
For any partition $\lambda$,
\[
I_K(P_\lambda(;q,t);q,t,T;\pm \sqrt{t},\pm\sqrt{qt})=0
\]
unless $\lambda$ is of the form $\mu^2$, in which case
\[
I_K(P_{\mu^2}(;q,t);q,t,T;\pm \sqrt{t},\pm\sqrt{qt})
=
\frac{C^0_\mu(T^2;q,t^2) C^-_\mu(qt;q,t^2)}
{C^0_\mu(qT^2/t;q,t^2) C^-_\mu(t^2;q,t^2)}
\]
\end{conj}

\begin{conj}\label{conj:vanish_U/O_T}
For any partition $\lambda$,
\[
I_K(P_\lambda(;q,t);q,t,T;\pm 1,\pm \sqrt{t})=0
\]
unless $\lambda$ is of the form $2\mu$, in which case
\[
I_K(P_{2\mu}(;q,t);q,t,T;\pm 1,\pm \sqrt{t})
=
\frac{C^0_\mu(T^2;q^2,t) C^-_\mu(q;q^2,t)}
{C^0_\mu(qT^2/t;q^2,t) C^-_\mu(t;q^2,t)}
\]
\end{conj}

(Note that these conjectures are equivalent, by duality.)  The vanishing
of these integrals is, of course, a natural conjecture; the nonzero values
will be justified below.

\begin{prop}
Conjecture \ref{conj:vanish_U/Sp_T} is equivalent to the following claim: For
all integers $n\ge 0$ and partitions $\lambda$ with $\ell(\lambda)\le 2n$,
\[
I^{(n)}_K(P_\lambda(x_1^{\pm 1},\dots x_n^{\pm
1};q,t);q,t;\pm\sqrt{t},\pm\sqrt{qt})=0
\]
unless $\lambda$ is of the form $\mu^2$.
\end{prop}

\begin{proof}
The given claim is the specialization $T=t^n$ of the vanishing part of
Conjecture \ref{conj:vanish_U/Sp_T}, and thus by rationality is equivalent
to the vanishing claim.  It thus remains to show that vanishing implies
\[
I^{(n)}_K(P_{\mu^2}(x_1^{\pm 1},\dots x_n^{\pm
1};q,t);q,t;\pm\sqrt{t},\pm\sqrt{qt})
=
\frac{C^0_\mu(t^n;q^2,t) C^-_\mu(q;q^2,t)}
{C^0_\mu(qt^{n-1};q^2,t) C^-_\mu(t;q^2,t)}.
\]
Suppose this is true for $\mu\subsetneq\lambda$, and choose
$\mu\subset\lambda$ so that $|\lambda/\mu|=1$.  Let $\nu$ be the unique
partition such that $\mu^2\subsetneq \nu\subsetneq \lambda^2$.  Now,
$e_1-e_{2n-1}$ is in the kernel of the homomorphism $f\mapsto f(x_1^{\pm
1},\dots x_n^{\pm 1})$, and thus
\[
I^{(n)}_K((e_1 P_\nu(;q,t))(x_1^{\pm 1},\dots x_n^{\pm
1});q,t;\pm\sqrt{t},\pm\sqrt{qt})
=
I^{(n)}_K((e_{2n-1} P_\nu(;q,t))(x_1^{\pm 1},\dots x_n^{\pm
1};q,t);q,t;\pm\sqrt{t},\pm\sqrt{qt}),
\]
Now, if we expand each side via the Pieri identity, only one term on each
side has a nonvanishing integral; we thus find
\[
\psi'_{\lambda^2/\nu}
I^{(n)}_K(P_{\lambda^2}(x_1^{\pm 1},\dots x_n^{\pm
1};q,t);q,t;\pm\sqrt{t},\pm\sqrt{qt})
=
\psi'_{(1^{2n}+\mu^2)/\nu}
I^{(n)}_K(P_{\mu^2}(x_1^{\pm 1},\dots x_n^{\pm
1};q,t);q,t;\pm\sqrt{t},\pm\sqrt{qt}).
\]
Solving this for
\[
I^{(n)}_K(P_{\lambda^2}(;q,t)(x_1^{\pm 1},\dots x_n^{\pm
1});q,t;\pm\sqrt{t},\pm\sqrt{qt})
\]
gives the desired result.
\end{proof}

\begin{rem}
Compare the computation of the nonzero value of the Macdonald inner product
in \cite[Section VI.9]{MacdonaldIG:1995}, and the computation of the nonzero
value of the Koornwinder inner product in \cite{vanDiejen:1996}.
\end{rem}

For the other version of the conjecture, we have the following refinement.

\begin{prop}\label{prop:conj_equiv_O1}
Let $m$ be a nonnegative integer.  Then the following claims are
equivalent.
\begin{itemize}
\item Conjecture \ref{conj:vanish_U/O_T} holds for all partitions $\lambda$
with $\lambda_1\le m$.
\item For all integers $n\ge 0$ and partitions $\lambda\subset m^{2n}$,
\[
I^{(n)}_K(P_\lambda(x_1^{\pm 1},\dots x_n^{\pm
1};q,t);q,t;\pm 1,\pm\sqrt{t})
+
I^{(n-1)}_K(P_\lambda(x_1^{\pm 1},\dots x_{n-1}^{\pm
1},1,-1;q,t);q,t;\pm t,\pm\sqrt{t})=0
\]
unless $\lambda$ is of the form $2\mu$.
\item For all integers $n\ge 0$ and partitions $\lambda\subset m^{2n+1}$,
\[
I^{(n)}_K(P_\lambda(x_1^{\pm 1},\dots x_n^{\pm
1},1;q,t);q,t;t,-1,\pm\sqrt{t})
+
I^{(n)}_K(P_\lambda(x_1^{\pm 1},\dots x_{n-1}^{\pm
1},-1;q,t);q,t;1,-t,\pm\sqrt{t})=0
\]
unless $\lambda$ is of the form $2\mu$.
\end{itemize}
\end{prop}

\begin{proof}
First, assume the vanishing portion of Conjecture \ref{conj:vanish_U/O_T}.
We cannot directly specialize $T\to t^n$, since this is a case in which
the order of specialization is important.  A consideration of possible
poles shows that the only partitions for which 
\[
\tilde{K}_\lambda(x_1^{\pm 1},\dots x_n^{\pm 1};q,t,t^n;\pm 1,\pm \sqrt{t})
\]
is ill-defined are the partitions $\lambda=1^k$ for $n+1\le k\le 2n$.
In this case, we find (by the Cauchy identity, say) that
\[
\tilde{K}_{1^k}(;q,t,T;\pm 1,\pm \sqrt{t})=e_k
\]
for all $k$, independent of $T$.  We thus find that
\[
I^{(n)}_K(\tilde{K}_{1^k}(x_1^{\pm 1},\dots x_n^{\pm 1};q,t,T;\pm 1,\pm
\sqrt{t})=0
\]
unless $k=2n$.  It follows that for $\ell(\lambda)<2n$,
\[
I^{(n)}_K(P_\lambda(x_1^{\pm 1},\dots x_n^{\pm
1};q,t);q,t;\pm 1,\pm\sqrt{t}) 
=
\lim_{T\to t^n}
I_K(P_\lambda(;q,t);q,t,T;\pm 1,\pm\sqrt{t}).
\]
Similarly, by the plethystic symmetry of Koornwinder integrals,
\begin{align}
I^{(n)}_K(P_\lambda(x_1^{\pm 1},\dots x_{n-1}^{\pm
1},1,-1;q,t);q,t;\pm t,\pm\sqrt{t})
&=
\lim_{T\to t^n}
I_K(P_\lambda([p_k+1+(-1)^k];q,t);q,t,T;\pm t,\pm\sqrt{t})\\
&=
\lim_{T\to t^n}
I_K(P_\lambda(p_k;q,t);q,t,T;\pm 1,\pm\sqrt{t}),
\end{align}
while for $\ell(\lambda)<2n+1$,
\begin{align}
I^{(n)}_K(P_\lambda(x_1^{\pm 1},\dots x_n^{\pm
1},1;q,t);q,t;t,-1,\pm\sqrt{t})
&=
\lim_{T\to t^{n+1/2}}
I_K(P_\lambda(;q,t);q,t,T;\pm 1,\pm\sqrt{t})\\
I^{(n)}_K(P_\lambda(x_1^{\pm 1},\dots x_n^{\pm
1},-1;q,t);q,t;1,-t,\pm\sqrt{t})
&=
\lim_{T\to t^{n+1/2}}
I_K(P_\lambda(;q,t);q,t,T;\pm 1,\pm\sqrt{t}).
\end{align}
More generally, if $\lambda=k^{2n}+\mu$ with $\ell(\mu)<2n$,
\begin{align}
I^{(n)}_K(P_\lambda(x_1^{\pm 1},\dots x_n^{\pm
1};q,t);q,t;\pm 1,\pm\sqrt{t}) 
&=
\lim_{T\to t^n}
I_K(P_\mu(;q,t);q,t,T;\pm 1,\pm\sqrt{t})\\
I^{(n)}_K(P_\lambda(x_1^{\pm 1},\dots x_{n-1}^{\pm
1},1,-1;q,t);q,t;\pm t,\pm\sqrt{t}) 
&=
(-1)^k
\lim_{T\to t^n}
I_K(P_\mu(;q,t);q,t,T;\pm 1,\pm\sqrt{t}),
\end{align}
and if $\lambda=k^{2n+1}+\mu$ with $\ell(\mu)<2n+1$,
\begin{align}
I^{(n)}_K(P_\lambda(x_1^{\pm 1},\dots x_n^{\pm
1},1;q,t);q,t;t,-1,\pm\sqrt{t}) 
&=
\lim_{T\to t^{n+1/2}}
I_K(P_\mu(;q,t);q,t,T;\pm 1,\pm\sqrt{t})\\
I^{(n)}_K(P_\lambda(x_1^{\pm 1},\dots x_{n-1}^{\pm
1},-1;q,t);q,t;1,-t,\pm\sqrt{t}) 
&=
(-1)^k
\lim_{T\to t^{n+1/2}}
I_K(P_\mu(;q,t);q,t,T;\pm 1,\pm\sqrt{t}).
\end{align}
Thus the vanishing portion of Conjecture \ref{conj:vanish_U/O_T} is
equivalent to the other two claims.

It remains to consider the nonzero values of the integrals.  We consider
the $2n$ version; the other is analogous.  It suffices to consider
\[
I^{(n)}_K(P_\lambda(x_1^{\pm 1},\dots x_{n-1}^{\pm
1},1,-1;q,t);q,t;\pm t,\pm\sqrt{t})
\]
when $\lambda$ is an even partition with $\lambda_1\le m$,
$\ell(\lambda)<2n$.  Let $\mu\subset\lambda$ be any partition such that
$\mu_1<m$ and $\lambda/\mu$ is a horizontal strip.  There is then a unique
even partition $\nu$ such that $\mu/\nu$ is a horizontal strip.  Now,
the function
\[
(e_{|\lambda/\mu|}+e_{2n-|\lambda/\mu|}) P_\mu(;q,t).
\]
is annihilated by the integral, since
\[
(e_{|\lambda/\mu|}+e_{2n-|\lambda/\mu|})(x_1^\pm 1, \dots x_{n-1}^{\pm
1},1,-1) = 0.
\]
On the other hand, every term $P_\kappa(;q,t)$ in the Pieri expansion has
$\kappa_1\le m$, and thus the only nonvanishing terms are
\[
\psi_{\lambda/\mu} P_\lambda(;q,t) +
\psi_{(1^{2n}+\nu)/\mu} P_{1^{2n}+\nu}(;q,t)
\]
(plus an additional factor of 2 if $|\mu/\nu|=n$).  Solving the resulting
recurrence gives the desired nonzero values of the integral.
\end{proof}

We can show these conjectures in a number of special cases.  We begin
with the Schur case.

\begin{thm}
Conjectures \ref{conj:vanish_U/Sp_T} and \ref{conj:vanish_U/O_T} hold when
$q=t$.
\end{thm}

\begin{proof}
We apply the propositions, and observe that the resulting integrals
can be expressed in terms of integrals over classical groups.  To
be precise, the $q=t$ cases of the conjectures are equivalent to the
claims
\begin{align}
\int_{U\in Sp(2n)} s_\lambda(U) dU &= 0\quad\text{unless $\lambda=\mu^2$}\\
\int_{U\in O(n)} s_\lambda(U) dU &= 0\quad\text{unless $\lambda=2\mu$},
\end{align}
where the integrals are with respect to the corresponding Haar measures.
These follow from the theory of zonal polynomials \cite[Chapter
VI]{MacdonaldIG:1995}, or equivalently from the theory of symmetric spaces
\cite{HelgasonS:1984} (specifically $U(2n)/Sp(2n)$ and $U(n)/O(n)$); this
also shows that the integrals are 1 when nonzero, agreeing with the general
formula.
\end{proof}

\begin{rem} In fact, our original motivation for the above conjectures was
to generalize these results to Macdonald polynomials.  The connection to
generalized Littlewood identities was then suggested by the results of
\cite[Section 5]{BaikJ/RainsEM:2001}.
\end{rem}

The theory of zonal polynomials gives us another special case.

\begin{thm}
Conjecture \ref{conj:vanish_U/Sp_T} holds in the case $(q,t)\mapsto
\lim_{q\to 1} (q^2,q)$.
Conjecture \ref{conj:vanish_U/O_T} holds in the case $(q,t)\mapsto
\lim_{q\to 1} (q,q^2)$.
\end{thm}

\begin{proof}
The expression of zonal polynomials in terms of Jack polynomials gives,
for partitions $\mu$ with $\ell(\mu)\le n$,
\begin{align}
\int_{U\in Sp(2n)} s_{\mu^2}(A U) &\propto \lim_{q\to 1} P_{\mu}(A J
A^t;q,q^2)\\
\int_{U\in O(2n)} s_{2\mu}(A U) &\propto \lim_{q\to 1} P_{\mu}(A A^t
A^t;q^2,q),
\end{align}
where $J$ is the symplectic inner product; the integrals vanish on
Schur functions not of the stated form.  Now, consider the integral
\[
\int_{U\in Sp(n)} \int_{U'\in O(2n)}
s_\lambda(U U').
\]
This vanishes unless $\lambda$ has both the form $2\kappa$ and the
form $\kappa^2$; that is, unless $\lambda=2\mu^2$ for some $\mu$.
Thus the integrals
\[
\int_{U\in Sp(2n)} \lim_{q\to 1} P_\lambda(U U^t;q^2,q)
\quad\text{and}\quad
\int_{U\in O(2n)} \lim_{q\to 1} P_\lambda(U J U^t;q,q^2)
\]
vanish unless $\lambda$ has the form $\mu^2$ and $2\mu$ respectively.  As
these can be written as integrals over the spaces $Sp(2n)/U(n)$ and
$O(2n)/U(n)$, we can express them as limiting cases of Koornwinder
integrals (more precisely, Jacobi integrals); the theorem follows.
\end{proof}

\begin{rem}
Since the theory of zonal polynomials extends via quantum groups to
the cases $P_\lambda(;q,q^2)$, $P_\lambda(;q^2,q)$ (without limits)
\cite{NoumiM:1996}, it should be possible to extend this argument
accordingly.
\end{rem}

We next turn to special cases with generic parameters, but for which
$\lambda$ has been constrained.

\begin{thm}
Conjecture \ref{conj:vanish_U/Sp_T} holds if $\lambda_1\le 1$;
Conjecture \ref{conj:vanish_U/O_T} holds if $\ell(\lambda)\le 1$.
\end{thm}

\begin{proof}
We consider the second claim; the first will follow by duality.
We thus need to show, for $u$ a formal variable:
\[
I_K(\sum_k u^k g_k;q,t,T;\pm 1,\pm \sqrt{t})
=
\sum_k
u^{2k}
\frac{(T;q^2)_k (qt;q^2)_k}
{(qT/t;q^2)_k (q^2;q^2)_k}
=
{}_2\phi_1(T^2,qt;qT^2/t;q^2,u^2).
\]
In fact, we claim that, more generally,
\[
I_K(\sum_k u^k g_k;q,t,T;\pm a,\pm \sqrt{t})
=
{}_2\phi_1(T^2,qt/a^2;T^2 q a^2/t;q^2,u^2 a^2).
\]
By rationality, we may assume $T=t^n$ for some integer $n\ge 0$.  Then
\begin{align}
I_K(\sum_k u^k g_k;q,t,T;\pm a,\pm \sqrt{t})
&=
I^{(n)}_K(\prod_j \frac{(t u x^{\pm 1}_j;q)}{(u x^{\pm 1}_j;q)};q,t;\pm
a,\pm \sqrt{t})\\
&=
I^{(1)}_K(\frac{(t^{(n+1)/2} u x^{\pm 1};q)}{(t^{-(n-1)/2} u x^{\pm 1};q)}
;q,t;\pm t^{(n-1)/2} a,\pm t^{n/2})\\
&=
\frac{
(T^2 a^2 u^2,T^2 u^2 t;q^2) (u^2 t,T a^2/t;q)}
{
(a^2 u^2,u^2 t;q^2)
(t T u^2, a^2 T^2/t;q)}\\*
&\phantom{{}={}}
{}_8W_7(
t T u^2/q;
tu/a,-tu/a,\sqrt{t}u,-\sqrt{t}u,T;q,T a^2b^2/t^2),
\end{align}
as long as $|a^2b^2T/t^2|<1$.  The claim follows by quadratic
transformation (equation (3.5.4) of \cite{GasperG/RahmanM:1990}).
\end{proof}

\begin{rem}
The conjectures can thus be viewed as multivariate analogues of quadratic
transformations.
\end{rem}

The same argument shows that Conjecture \ref{conj:vanish_U/O_T} holds whenever
$\lambda_1\le 1$.  It turns out that we can show something much stronger;
to do so, we will need yet another equivalent form of the conjectures.

\begin{thm}
Let $m$ be a nonnegative integer.  Then the following statements are
equivalent.
\begin{itemize}
\item Conjecture \ref{conj:vanish_U/O_T} holds for all partitions $\lambda$
with $\lambda_1\le m$.
\item For all integers $n\ge 0$, we have
\[
[P_\lambda(x_1,\dots x_m;q,t)]
(x_1\dots x_m)^{n/2} P^{D_m}_{n\omega_m}(x_1^{-1},\dots x_m^{-1};q,t)
=
0
\]
unless $\lambda=\mu^2$ for some $\mu$, where $P^{D_m}_{n\omega_m}(;q,t)$ is
the $D_m$-type Macdonald polynomial associated to the weight
\[
n\omega_m = (n/2,n/2,\dots n/2).
\]
\end{itemize}
\end{thm}

\begin{proof}
We will show that the second claim is equivalent to the second and third
claims of Proposition \ref{prop:conj_equiv_O1}.  By orthogonality of
Koornwinder polynomials, these may be written in the form
\begin{align}
[K^{(n)}_0(x_1,\dots x_n;q,t;\pm 1,\pm\sqrt{t})]&
P_\lambda(x_1^{\pm 1},\dots x_n^{\pm 1};q,t)\notag\\*
&+
[K^{(n-1)}_0(x_1,\dots x_{n-1};q,t;\pm t,\pm\sqrt{t})]
P_\lambda(x_1^{\pm 1},\dots x_{n-1}^{\pm 1},1,-1;q,t)=0\\
[K^{(n)}_0(x_1,\dots x_n;q,t;t,-1,\pm\sqrt{t})]&
P_\lambda(x_1^{\pm 1},\dots x_n^{\pm 1},1;q,t)\notag\\*
&+
[K^{(n)}_0(x_1,\dots x_n;q,t;1,-t,\pm\sqrt{t})]
P_\lambda(x_1^{\pm 1},\dots x_{n}^{\pm 1},-1;q,t)=0
\end{align}
(unless $\lambda=2\mu$).  By the Cauchy identity for Koornwinder
polynomials, we compute
{\allowdisplaybreaks
\begin{align}
[K^{(n)}_0(x_1,\dots x_n&;q,t;\pm 1,\pm\sqrt{t})]
P_\lambda(x_1^{\pm 1},\dots x_n^{\pm 1};q,t)
\notag\\*
&=(-1)^{|\lambda|}
[P_{\lambda'}(y_1,\dots y_m;t,q)]
\prod_{1\le j\le m} y_j^n
K^{(m)}_{n^m}(y_1,\dots y_m;t,q;\pm 1,\pm\sqrt{t})\\
[K^{(n-1)}_0(x_1,\dots x_{n-1}&;q,t;\pm t,\pm\sqrt{t})]
P_\lambda(x_1^{\pm 1},\dots x_{n-1}^{\pm 1},1,-1;q,t)
\notag\\*&=
(-1)^{|\lambda|}
[P_{\lambda'}(y_1,\dots y_m;t,q)]
\prod_{1\le j\le m} y_j^n (y_j^{-1}-y_j)
K^{(m)}_{(n-1)^m}(y_1,\dots y_m;t,q;\pm t,\pm\sqrt{t})\\
[K^{(n)}_0(x_1,\dots x_n&;q,t;t,-1,\pm\sqrt{t})]
P_\lambda(x_1^{\pm 1},\dots x_n^{\pm 1},1;q,t)
\\*&=
(-1)^{|\lambda|}
[P_{\lambda'}(y_1,\dots y_m;t,q)]
\prod_{1\le j\le m} y_j^{n+1/2} (y_j^{-1/2}-y_j^{1/2})
K^{(m)}_{n^m}(y_1,\dots y_m;t,q;t,-1,\pm\sqrt{t})\notag\\
[K^{(n)}_0(x_1,\dots x_n&;q,t;1,-t,\pm\sqrt{t})]
P_\lambda(x_1^{\pm 1},\dots x_n^{\pm 1},-1;q,t)
\\*&=
(-1)^{|\lambda|}
[P_{\lambda'}(y_1,\dots y_m;t,q)]
\prod_{1\le j\le m} y_j^{n+1/2} (y_j^{-1/2}+y_j^{1/2})
K^{(m)}_{n^m}(y_1,\dots y_m;t,q;1,-t,\pm\sqrt{t}).\notag
\end{align}
}
By the considerations of \cite[Section 5.4]{vanDiejen:1995}, we find
\[
K^{(m)}_{n^m}(y_1,\dots y_m;t,q;\pm 1,\pm\sqrt{t})
+
\prod_{1\le j\le m} (y_j^{-1}-y_j)
K^{(m)}_{n^m}(y_1,\dots y_m;t,q;\pm t,\pm\sqrt{t})
=
P^{D_m}_{2n\omega_m}(y_1^{-1},\dots y_m^{-1};t,q),
\]
and similarly
\begin{align}
\prod_{1\le j\le m} (y_j^{-1/2}-y_j^{1/2})&
K^{(m)}_{n^m}(y_1,\dots y_m;t,q;t,-1,\pm\sqrt{t})\\*
&+
\prod_{1\le j\le m} (y_j^{-1/2}+y_j^{1/2})
K^{(m)}_{n^m}(y_1,\dots y_m;t,q;-t,1,\pm\sqrt{t})
=
P^{D_m}_{(2n+1)\omega}(y_1^{-1},\dots y_m^{-1};t,q).\notag
\end{align}
The theorem follows.
\end{proof}

\begin{cor}
Conjecture \ref{conj:vanish_U/Sp_T} holds whenever $\ell(\lambda)\le 4$;
Conjecture \ref{conj:vanish_U/O_T} holds whenever $\lambda_1\le 4$.
\end{cor}

\begin{proof}
By the theorem, we must show that
\[
[P_\lambda(x_1,\dots x_4;q,t)]
(x_1\dots x_4)^{n/2} P^{D_4}_{n\omega_4}(x_1^{-1},\dots x_4^{-1};q,t)
=
0
\]
unless $\lambda$ is of the form $\mu^2$.  Now, the triality automorphism of
$D_4$ (which still applies in the Macdonald setting) implies the identity
\begin{align}
(x_1x_2x_3x_4)^{n/2}
P^{D_4}_{n\omega_4}(x_1^{-1},\dots x_4^{-1};q,t)
&=
u^n P^{D_4}_n(u,x_1x_2/u,x_1x_3/u,x_1x_4/u;q,t),\\
&=
u^n K^{(4)}_n(u,x_1x_2/u,x_1x_3/u,x_1x_4/u;q,t;\pm 1,\pm\sqrt{q}),
\end{align}
where $u=\sqrt{x_1x_2x_3x_4}$.  By triangularity, this is a linear
combination of the polynomials
\[
u^n P_k((u)^{\pm 1},(x_1x_2/u)^{\pm 1},(x_1x_3/u)^{\pm 1},(x_1x_4/u)^{\pm 1};q,t)
\]
for $k\le n$; by symmetry, only those $k$ having the same parity as $n$
occur.  Since for $\ell(\lambda)\le 4$,
\[
u^{2l} P_\lambda(x_1,\dots x_4) = P_{l^4+\lambda}(x_1,\dots x_4),
\]
preserving the constraint $\lambda=\mu^2$, we find that it suffices to
show that
\[
[P_\lambda(x_1,x_2,\dots x_4;q,t)]
u^k P_k((u)^{\pm 1},(x_1x_2/u)^{\pm 1},(x_1x_3/u)^{\pm 1},(x_1x_4/u)^{\pm 1};q,t)
\]
for $\lambda\ne \mu^2$.  Now, we have the generating function
\begin{align}
\sum_k v^k u^k P_k((u)^{\pm 1},(x_1x_2/u)^{\pm 1},{}&(x_1x_3/u)^{\pm 1},(x_1x_4/u)^{\pm 1};q,t)\notag\\*
&=
\frac
{(tvx_1x_2x_3x_4,tv,tvx_1x_2,tvx_1x_3,tvx_1x_4,tvx_2x_3,tvx_2x_4,tvx_3x_4;q)}
{(vx_1x_2x_3x_4,v,vx_1x_2,vx_1x_3,vx_1x_4,vx_2x_3,vx_2x_4,vx_3x_4;q)}\\
&=
\frac
{(tvx_1x_2x_3x_4,tv;q)}
{(vx_1x_2x_3x_4,v;q)}
\sum_{\ell(\mu)\le 2}
\frac{C^-_\mu(t;q,t^2)}{C^-_\mu(q;q,t^2)}
P_{\mu^2}(x_1,x_2,x_3,x_4;q,t),
\end{align}
by Macdonald's generalized Littlewood conjecture.  The factors out front
have no effect on the vanishing requirement; the corollary follows.
\end{proof}

\begin{rem}
In particular, the conjectures hold if $|\lambda|\le 5$.
\end{rem}

\bigskip
We observed above that the case $q=t$ of Conjectures
\ref{conj:vanish_U/Sp_T} and \ref{conj:vanish_U/O_T} follows from the
theory of symmetric spaces, specifically the spaces $U(2n)/Sp(2n)$ and
$U(n)/O(n)$.  It is thus natural to wonder whether one can formulate
similar conjectures for other symmetric spaces.  This indeed is the case;
for instance, the analogous ``conjecture'' for spaces of the form $G\times
G/G$ results is simply the orthogonality of the Macdonald polynomials for
the associated root system.  For the other classical symmetric spaces,
the situation turns out to be more complicated, as we shall see below.

One approach to generating such conjectures is simply to make an educated
guess based on the form of the integral for $q=t$.  For the Grassmannian
$U(m+n)/U(m)\times U(n)$ with $m\le n$, the Schur case is
\[
\int_{U_1\in U(m),U_2\in U(n)} s_\lambda(U_1\oplus U_2) dU_1 dU_2
=
0,
\]
unless the dominant weight $\lambda$ of $U(m+n)$ satisfies
\begin{align}
\lambda_i+\lambda_{m+n+1-i}&=0,\quad 1\le i\le m\\
\lambda_i&=0,\quad m+1\le i\le n-m,
\end{align}
in which case the integral is $1$.  This condition can be stated more
concisely as $\lambda=\mu\overline{\mu}$ for $\ell(\mu)\le m$, where
$\mu\overline{\nu}$ denotes the dominant weight of $U(m+n)$ with positive
part $\mu$ and negative part $0^{m+n}-\nu$.  This immediately suggests the
following conjecture.  Here and for the remainder of this section,
we take the convention that a factor $1/Z$ in front of an integral of
a Macdonald or Koornwinder polynomial over a weight function is the
constant that makes the integral 1 when the polynomial is trivial.

\begin{conj}\label{conj:vanish_Ugrass1}
Let $m$ and $n$ be integers with $0\le m\le n$.  Then
for a dominant weight $\mu\overline{\nu}$ of $U(m+n)$,
\[
\frac{1}{Z}
\int
P_{\mu\overline{\nu}}(x_1,\dots x_m,y_1,\dots y_n)
\prod_{1\le i\ne j\le m} \frac{(x_i/x_j;q)}{(tx_i/x_j;q)}
\prod_{1\le i\ne j\le n} \frac{(y_i/y_j;q)}{(ty_i/y_j;q)}
\prod_{1\le i\le m} \frac{dx_i}{2\pi \sqrt{-1} x_i}
\prod_{1\le i\le n} \frac{dy_i}{2\pi \sqrt{-1} y_i}
=
0
\]
unless $\mu=\nu$ and $\ell(\mu)\le m$, in which case the integral is
\[
\frac{
C^-_\mu(q;q,t)
C^+_\mu(t^{m+n-2} q;q,t)
C^0_\mu(t^n,t^m;q,t)
}{
C^-_\mu(t;q,t)
C^+_\mu(t^{m+n-2} t;q,t)
C^0_\mu(q t^{n-1},q t^{m-1};q,t)
}.
\]
\end{conj}

\begin{rems}
Unlike Conjectures \ref{conj:vanish_U/Sp_T} and \ref{conj:vanish_U/O_T} (as
well as the other conjectures below), the nonzero value here has not been
computed via Pieri identities, but has merely been guessed from low-order
examples.
\end{rems}

\begin{rems}
There is an obvious analogue for Koornwinder polynomials (related to the
Grassmannians $O(m+n)/O(m)\times O(n)$ and $Sp(2m+2n)/Sp(2m)\times Sp(2n)$)
but we have not been able to test it enough to justify making a formal
conjecture.
\end{rems}

The reason why it was relatively easy to formulate conjectures for the
spaces $U(n)/O(n)$, $U(2n)/Sp(2n)$ is that in those cases, the rank of the
smaller group is about half the rank of the bigger group.  This, for
instance, is what allowed us to compute the nonzero values via Pieri
identities.  In the remaining cases, the rank differs, if at all, by only
1, and thus the vanishing condition is not enough to determine the weight
function.  In a number of cases, however, the small group is most naturally
taken to be disconnected, and while the rank of the identity component is
indeed large, the effective rank of the nonidentity component is often much
smaller.

The simplest example of this is the case $U(2n)/U(n)\times U(n)$.  As the
stabilizer group of a symmetric space, $U(n)\times U(n)$ is the subgroup
preserved by an involution acting on $U(2n)$; to be precise, it is the
centralizer of the element
\[
\begin{pmatrix} I_n&0\\0&-I_n\end{pmatrix}.
\]
Now, the element
\[
\begin{pmatrix} 0&I_n\\I_n&0\end{pmatrix},
\]
while not preserved by the involution, is at least preserved up to sign;
furthermore, it normalizes $U(n)\times U(n)$, acting by switching the two
unitary groups.  If we integrate a Schur function over the corresponding
coset of $U(n)\times U(n)$, the integral vanishes on the same weights, and
evaluates to $\pm 1$ where nonzero.  To extend this to the Macdonald case,
we observe that the eigenvalues of an element of this coset come in $\pm$
pairs; we thus wish an integral of the form
\[
\int P_{\mu\overline{\nu}}(\pm \sqrt{x_1},\pm \sqrt{x_2},\dots
\pm \sqrt{x_n};q,t) w(x)
\prod_{1\le j\le n} \frac{dx_j}{2\pi \sqrt{-1} x_j},
\]
vanishing unless $\mu=\nu$.  Since
\[
e_1(\pm\sqrt{x_1},\pm\sqrt{x_2},\dots \pm\sqrt{x_n})=0,
\]
the Pieri identity argument applies here, and thus the weight function (if
it exists) is unique.  By examining low-rank cases, we are led to the
following conjecture.

\begin{conj}\label{conj:vanish_U_grass2}
For any integer $n\ge 0$, and any dominant weight $\mu\overline\nu$ of
$U(2n)$,
\[
\frac{1}{Z}
\int
P_{\mu\overline{\nu}}(\pm\sqrt{x_1},\pm\sqrt{x_2},\dots \pm\sqrt{x_n};q,t)
\prod_{1\le i\ne j\le n} \frac{(x_i/x_j;q^2)}{(t^2 x_i/x_j;q^2)}
\prod_{1\le i\le n} \frac{dx_i}{2\pi\sqrt{-1}x_i}
=
0
\]
unless $\mu=\nu$, when the integral is
\[
\frac{(-1)^{|\mu|}C^-_\mu(q;q,t) C^+_\mu(t^{2n-2} q;q,t) C^0_\mu(t^n,-t^n;q,t)}
{C^-_\mu(t;q,t) C^+_\mu(t^{2n-2} t;q,t) C^0_\mu(q t^{n-1},-q t^{n-1};q,t)}
\]
\end{conj}

By analogy with Proposition \ref{prop:conj_equiv_O1}, we would have
expected the nonzero values of this integral to differ from the nonzero
values in Conjecture \ref{conj:vanish_Ugrass1} by only a sign factor.
It turns out that these values are (conjecturally) attained by another
nice integral.

\begin{conj}
For any integer $n\ge 0$, and any dominant weight $\mu\overline\nu$ of
$U(2n)$,
\[
\frac{1}{Z}
\int
P_{\mu\overline{\nu}}(x_1,\dots x_n,y_1,\dots y_n;q,t)
\prod_{1\le i,j\le n}
\frac{(q x_i/y_j,q y_i/x_j;q^2)}{(t x_i/y_j,t y_i/x_j;q^2)}
\prod_{1\le i\ne j\le n}
\frac{(x_i/x_j,y_i/y_j;q^2)}{(q t x_i/x_j,q t y_i/y_j;q^2)}
\prod_{1\le i\le n} \frac{dx_i}{2\pi\sqrt{-1}x_i}
=
0
\]
unless $\mu=\nu$, when the integral is
\[
\frac{C^-_\mu(q;q,t) C^+_\mu(t^{2n-2} q;q,t) C^0_\mu(t^n,-t^n;q,t)}
{C^-_\mu(t;q,t) C^+_\mu(t^{2n-2} t;q,t) C^0_\mu(q t^{n-1},-q t^{n-1};q,t)}
\]
\end{conj}

\begin{rem}
Note that the weight function in this case is not of a standard Macdonald
or Koornwinder form.  The associated orthogonal polynomials may be of
interest.
\end{rem}

For the spaces $O(2n)/O(n)\times O(n)$ and $Sp(2n)/U(n)$, we have the
following conjecture.  Here $[2p_{k/2}]$ represents the homomorphism such
that $p_{2k+1}\to 0$, $p_{2k}\to 2p_k$; this is just the infinite variable
analogue of the specialization
$\pm\sqrt{x_1},\pm\sqrt{x_2},\dots \pm\sqrt{x_n}$.

\begin{conj}\label{conj:vanish_O_grass}
For any partition $\lambda$,
\[
I_K(\tilde{K}_\lambda([2p_{k/2}];q,t,T;a,-a,b,-b);q^2,t^2,T;-1,-t,a^2,b^2)
=0
\]
unless $\lambda$ is of the form $2\mu$, in which case the integral is
\[
\frac{
(-1)^{|\mu|}
C^-_\mu(q;q^2,t)
C^+_\mu(a^2 b^2 T^2/t^3;q^2,t)
C^0_\mu(T,-a^2 T/t,-b^2 T/t,a^2 b^2 T/t^2;q^2,t)
}{
C^-_\mu(t;q^2,t)
C^+_\mu(a^2 b^2 T^2/qt^2;q^2,t)
C^0_{2\mu}(a^2 b^2 T^2/t^3;q^2,t^2)
}.
\]
\end{conj}

Similarly, for the spaces $O(4n)/U(2n)$ and $Sp(4n)/Sp(2n)\times Sp(2n)$,

\begin{conj}\label{conj:vanish_Sp_grass}
For any partition $\lambda$,
\[
I_K(\tilde{K}_\lambda([2p_{k/2}];q,t,T;a,-a,b,-b);q^2,t^2,T;-t,-qt,a^2,b^2)
=0
\]
unless $\lambda$ is of the form $\mu^2$, in which case the integral is
\[
\frac{
(-1)^{|\mu|}
C^-_\mu(qt;q,t^2)
C^+_\mu(a^2 b^2 T^2/t^4;q,t^2)
C^0_\mu(T,-a^2 T/t,-b^2 T/t,a^2 b^2 T/t^2;q,t^2)
}{
C^+_\mu(a^2 b^2 T^2/qt^3;q,t^2)
C^-_\mu(t^2;q,t^2)
C^0_{\mu^2}(a^2 b^2 T^2 q/t^2;q^2,t^2)
}.
\]
\end{conj}

In these cases, we have no conjectured weight function corresponding to the
values with the sign factors removed; the problem is that the two
Schur cases associated to each integral not only break the $BC_n$ symmetry,
but do so in different ways.

It turns out that Propositions \ref{prop:sfkoorn:connt} and
\ref{prop:sfkoorn:brancht} produce integrals associated to orthogonal group
Grassmannians.  For $O(2n)/O(1)\times O(2n-1)$:

\begin{prop}
For any partition $\lambda$,
\[
I_K(\tilde{K}_\lambda([p_k+t_0^k+t_0^{-k}];q,t,t T;t_0,t_1,t_2,t_3)
;q,t,T;t_0 t,t_1,t_2,t_3)
=
0
\]
unless $\ell(\lambda)\le 1$, in which case the integral is
\[
t_0^{-\lambda_1}
\frac{(T t_0 t_1,T t_0 t_2,T t_0 t_3,T t_0t_1t_2t_3/t;q)_{\lambda_1}}
{(q^{\lambda_1-1} T^2 t_0t_1t_2t_3,T^2 t_0t_1t_2t_3/t;q)_{\lambda_1}}.
\]
\end{prop}

Similarly, for $O(2n+1)/O(1)\times O(2n)$,

\begin{prop}
For any partition $\lambda$,
\[
I_K(\tilde{K}_\lambda(;q,t,T;t_0 t,t_1,t_2,t_3)
;q,t,T;t_0,t_1,t_2,t_3)
=
0
\]
unless $\ell(\lambda)\le 1$, in which case the integral is
\[
t_0^{\lambda_1}
\frac
{(T,T t_1 t_2/t,T t_1 t_3/t,T t_2t_3/t;q)_{\lambda_1}}
{(q^{\lambda_1-1} T^2 t_0t_1t_2t_3/t,T^2 t_0t_1t_2t_3/t^2;q)_{\lambda_1}}.
\]
\end{prop}

Less trivially, for the nonidentity component of $O(2n+1)/O(2)\times O(2n-1)$,

\begin{thm}
For any partition $\lambda$,
\[
I_K(\tilde{K}_\lambda([p_k+a^k+(-a)^k+a^{-k}+(-a)^{-k}]
;q,t,t^2 T;a,-a,b,-b);q,t,T;a t,-a t,b,-b)
=
0
\]
unless $\ell(\lambda)\le 2$ and $|\lambda|$ is even, in which case
the integral is generically nonzero and admits a factorization into
$q$-symbols.
\end{thm}

\begin{proof}
By two consecutive applications of the quasi-branching rule, the integral
evaluates to
\[
\sum_{0\prec\kappa'\prec\lambda'}
\psi^{(i)}_{\lambda/\kappa}(t T;q,t,T\sqrt{a^2b^2 t/q})
\psi^{(i)}_{\kappa/0}(t T;q,t,T\sqrt{a^2b^2/q})
\frac{k^0_{\kappa}(q,t,t T;a{:}-at,b,-b)
k^0_\lambda(q,t,t^2 T;-a{:}a,b,-b)}
{k^0_\kappa(q,t,t T;-at{:}a,b,-b)}.
\]
This is clearly nonzero unless $\ell(\lambda)\le 2$, in which case
$\ell(\kappa)\le 1$.  The sum turns out to be proportional to a terminating
very-well-poised ${}_8W_7$, summable by Equation II.16 of
\cite{GasperG/RahmanM:1990}.
\end{proof}

Finally, for the nonidentity component of $O(2n+2)/O(2)\times O(2n)$,
a similar calculation gives

\begin{thm}
For any partition $\lambda$,
\[
I_K(\tilde{K}_\lambda(;q,t,T;a t,-a t,b,-b);q,t,T;a,-a,b,-b)
=
0
\]
unless $\ell(\lambda)\le 2$ and $|\lambda|$ is even, in which case the
integral is generically nonzero, and admits a factorization into
$q$-symbols.
\end{thm}

It is likely that there are a number of other ``nice'' integrals satisfying
appropriate vanishing conditions, but a more systematic method of
construction will likely be needed to find them.  It is, however, unclear
to what extent our existing conjectures can be systematized, especially
given the multiple (untwisted) integrals associated to $U(2n)/U(n)\times
U(n)$.

\end{document}